\documentclass[aos,MSNbibl,seceqn,dvips]{arximspdf}
\usepackage{graphicx}
\doi{10.1214/13-AOS1121} %kopijuoti is PTS
\volume{41}
\issue{3}
\pubyear{2013}
\firstpage{1593}
\lastpage{1641}

\makeatletter

\newcommand{\rright}{\right}
\newcommand{\lleft}{\left}
\newcommand{\rrVert}{\Vert}
\newcommand{\rrvert}{\vert}
\newcommand{\llVert}{\Vert}
\newcommand{\llvert}{\vert}
\newtheorem{thmm}{Theorem}[section]
\newtheorem{lem}[thmm]{Lemma}
\newtheorem{prop}[thmm]{Proposition}
\newtheorem{cor}[thmm]{Corollary}
\newproclaim{rem}[thmm]{Remark}
\newproclaim{ass}[thmm]{Assumption}

\newcommand{\iint}{\int\!\!\!\int}

\newcommand{\mca}{\mathcal{A}}
\newcommand{\mcc}{\mathcal{C}}\newcommand{\mcd}{\mathcal{D}}
\newcommand{\mce}{\mathcal{E}}\newcommand{\mcf}{\mathcal{F}}
\newcommand{\mcg}{\mathcal{G}}
\newcommand{\mcj}{\mathcal{J}}
\newcommand{\mcl}{\mathcal{L}}
\newcommand{\mcn}{\mathcal{N}}
\newcommand{\mcx}{\mathcal{X}}

\newcommand{\mbbg}{\mathbb{G}}
\newcommand{\mbbm}{\mathbb{M}}\newcommand{\mbbn}{\mathbb{N}}
\newcommand{\mbbq}{\mathbb{Q}}\newcommand{\mbbr}{\mathbb{R}}
\newcommand{\mbbv}{\mathbb{V}}
\newcommand{\mbby}{\mathbb{Y}}\newcommand{\mbbz}{\mathbb{Z}}
\newcommand{\mbbrp}{\mathbb{R}_{+}}
\newcommand{\mbbzp}{\mathbb{Z}_{+}}

\newcommand{\mbF}{\mathbf{F}}

\newcommand{\al}{\alpha} \newcommand{\ep}{\varepsilon}
\newcommand{\vp}{\varphi} \newcommand{\del}{\delta}
\newcommand{\D}{\Delta} \newcommand{\Sig}{\Sigma}
\newcommand{\Lam}{\Lambda} \newcommand{\Gam}{\Gamma}

\newcommand{\p}{\partial}
\newcommand{\cil}{\to^{\mcl}} % <- Convergence in law
\newcommand{\argmax}{\mathop{\operatorname{argmax}}}

\def\nn{\nonumber}

\makeatother

\begin{document}
\begin{frontmatter}

\title{Convergence of Gaussian quasi-likelihood random fields
for ergodic L\'evy driven SDE observed at~high frequency}
\runtitle{Estimation of ergodic L\'evy driven SDE}

\begin{aug}
\author[A]{\fnms{Hiroki} \snm{Masuda}\corref{}\thanksref{t1}\ead[label=e1]{hiroki@imi.kyushu-u.ac.jp}}
\thankstext{t1}{Supported in part by JSPS KAKENHI Grant Number 20740061
and 23740082.}
\runauthor{H. Masuda}
\affiliation{Kyushu University}
\address[A]{Institute of Mathematics for Industry\\
Kyushu University\\
744 Motooka Nishi-ku Fukuoka 819-0395\\
Japan\\
\printead{e1}} %adresu isvedimo komanda gale!
\end{aug}

% HISTORY:
\received{\smonth{2} \syear{2012}}
\revised{\smonth{4} \syear{2013}}

% ABSTRACT
%
\begin{abstract}
This paper investigates the Gaussian quasi-likelihood estimation of
an exponentially ergodic multidimensional Markov process,
which is expressed as a solution to a L\'evy driven stochastic
differential equation
whose coefficients are known except for the finite-dimensional
parameters to be estimated,
where the diffusion coefficient may be degenerate or even null.
We suppose that the process is discretely observed under the rapidly
increasing experimental design
with step size $h_{n}$.
By means of the polynomial-type large deviation inequality,
convergence of the corresponding statistical random fields is derived
in a mighty mode,
which especially leads to the asymptotic normality at rate $\sqrt {nh_{n}}$
 for all the target parameters,
and also to the convergence of their moments.
As our Gaussian quasi-likelihood solely looks at the local-mean and
local-covariance structures,
efficiency loss would be large in some instances.
Nevertheless, it has the practically important advantages:
first, the computation of estimates does not require any fine tuning,
and hence it is straightforward;
second, the estimation procedure can be adopted without full
specification of the L\'evy measure.
\end{abstract}

% KEYWORDS

\end{frontmatter}

%s1 #&#
\section{Introduction}\label{sec1}

Let $X=(X_{t})_{t\in\mbbrp}$ be a solution to the stochastic
differential equation (SDE)
%
%
%e1.1 #&#
\begin{equation}
dX_{t}=a(X_{t},\al)\,dt+b(X_{t},
\beta)\,dW_{t}+c(X_{t-},\beta)\,dJ_{t}, \label{SDE}
\end{equation}
where the ingredients involved are as follows:
\begin{itemize}
\item the finite-dimensional unknown parameter
\[
\theta=(\al,\beta)\in\Theta_{\al}\times\Theta_{\beta}=:\Theta,
\]
where, for simplicity, the parameter spaces
$\Theta_{\al}\subset\mbbr^{p_{\al}}$ and $\Theta_{\beta}\subset
\mbbr
^{p_{\beta}}$ are supposed to be
bounded convex domains;
the parameter $\al$ (resp., $\beta$) affects local trend (resp., local
dispersion);

\item an $r'$-dimensional standard Wiener process $W$ and an $r''$-dimensional
centered pure-jump L\'evy process $J$, whose L\'evy measure is denoted
by $\nu$;

\item the initial variable $X_{0}$ independent of $(W,J)$,
with $\eta:=\mcl(X_{0})$ possibly depending on $\theta$;

\item the measurable functions $a\dvtx \mbbr^{d}\times\Theta_{\al}\to
\mbbr^{d}$,
$b\dvtx \mbbr^{d}\times\Theta_{\beta}\to\mbbr^{d}\otimes\mbbr^{r'}$,
and $c\dvtx \mbbr^{d}\times\Theta_{\beta}\to\mbbr^{d}\otimes\mbbr^{r''}$.
\end{itemize}
Incorporation of the jump part extends the continuous-path diffusion
parametric model,
which are nowadays widely used in many application fields.
We denote by $P_{\theta}$ the image measure of a solution process $X$
associated with $\theta\in\Theta\subset\mbbr^{p}$, where
$p:=p_{\al
}+p_{\beta}$.
Suppose that the true parameter $\theta_{0}=(\al_{0},\beta_{0})\in
\Theta
$ does exist,
with $P_{0}$ denoting the shorthand for the true image measure
$P_{\theta_{0}}$,
and that $X$ is not completely (continuously) observed but only
discretely at high frequency
under the condition for the rapidly increasing experimental design:
we are given a sample $(X_{t_{0}},X_{t_{1}},\ldots,X_{t_{n}})$,
where $t_{j}=t^{n}_{j}=jh_{n}$ for some $h_{n}>0$ such that
%
%
%e1.2 #&#
\begin{equation}
T_{n}:=nh_{n}\to\infty\quad\mbox{and}\quad nh_{n}^{2}
\to0 \label{RIED}
\end{equation}
for $n\to\infty$. The main objective of this paper is to estimate
$\theta_{0}$
under the exponential ergodicity of $X$;
the equidistant sampling assumption can be weakened to some extent
as long as the long-term and high-frequency framework is concerned;
however, it is just a technical extension making the presentation
notationally messy, and
hence we do not deal with it in the main context to make the
presentation more clear.

It is common knowledge that the maximum likelihood estimation is
generally infeasible,
since the transition probability is most often unavailable in a closed form.
This implies that the conventional statistical analyses based on the
genuine likelihood have no utility.
For this reason, we have to resort to some other feasible estimation
procedure, which could be a lot of things.
Among several possibilities, we are concerned here with
the \textit{Gaussian quasi-likelihood (\textit{GQL})} function defined
as if the conditional distributions of $X_{t_{j}}$ given $X_{t_{j-1}}$
are Gaussian
with approximate but explicit mean vector and covariance matrix; see
(\ref{lGa}) below.

The terminology ``quasi-likelihood'' has originated as the pioneering
work of Wedderburn~\cite{Wed74},
the concept of which formed a basis of the generalized linear regression.
The GQL-based estimation has been known to have
the advantage of computational simplicity and robustness for
misspecification of the noise distribution,
and is well-established as a fundamental tool in estimating possibly
non-Gaussian and dependent statistical models.
Just to be a little more precise,
consider a time-series $Y_{1},\ldots,Y_{n}$ in $\mbbr$ with a fixed~$Y_{0}$,
and\vadjust{\goodbreak}
denote by $m_{j-1}(\theta)\in\mbbr$ and $v_{j-1}(\theta)>0$ the
conditional mean and conditional variance of $Y_{j}$
given $(Y_{0},\ldots,Y_{j-1})$, where $\theta$ is an unknown parameter
of interest.
Then, the \textit{Gaussian quasi maximum likelihood estimator (GQMLE)}
is defined to be a maximizer of the function
\[
\theta\mapsto\sum_{j=1}^{n}\log \biggl\{
\frac{1}{\sqrt{2\pi
v_{j-1}(\theta
)}}\exp \biggl( -\frac{(Y_{j}-m_{j-1}(\theta))^{2}}{2v_{j-1}(\theta)} \biggr) \biggr\}.
\]
Namely, we compute the likelihood of $(Y_{1},Y_{2},\ldots,Y_{n})$ as if
the conditional law of $Y_{j}$ given $(Y_{1},\ldots,Y_{j-1})$ is
Gaussian with mean $m_{j-1}(\theta)$
and variance $v_{j-1}(\theta)$, so that only the structures of the
conditional mean and variance do matter.
Although it is not asymptotically efficient in general, it can serve as
a widely applicable estimation procedure.
One can consult Heyde~\cite{Hey97} for an extensive and systematic
account of statistical inference based on the GQL.
The GQL has been a quite popular tool for (semi)parametric estimation,
and especially there exists a vast amount of literature concerning
asymptotics of the GQL
for time series models with possibly non-Gaussian error sequence;
among others, we refer to Straumann and Mikosch~\cite{StrMik06} for
a class of conditionally heteroscedastic time series models,
and Bardet and Wintenburger~\cite{BarWin09} for multidimensional causal
time series,
as well as the references therein.

Let us return to our framework. On one hand, for the diffusion case
(where $c\equiv0$),
the GQL-estimation issue has been solved under some regularity
conditions, especially the GQL, which leads to an asymptotically
efficient estimator, where
the crucial point is that the optimal rates of convergence for
estimating $\al$ and $\beta$ are different
and given by $\sqrt{T_{n}}$ and $\sqrt{n}$, respectively;
see Gobet~\cite{Gob02} for the local asymptotic normality of the
corresponding statistical experiments.
For how to construct an explicit contrast function,
we refer to Yoshida~\cite{Yos92} and Kessler~\cite{Kes97} as well as
the references therein;
specifically, they employed a discretized version of the
continuous-observation likelihood process,
and a higher order local-Gauss approximation of the transition density,
respectively.
S\o rensen~\cite{Sor08} includes an extensive bibliography of many
existing results, including
explicit martingale estimating functions for discretely observed
diffusions (not necessarily at high frequency).
On the other hand, the issue has not been addressed enough in the
presence of jumps (possibly of infinite variation).
The question we should then ask is what will occur when one adopts the
GQL function.
In this paper, we will provide sufficient conditions under which
the GQL random field associated with our statistical experiments
converges in a mighty mode; see Section~\ref{sec_MC}.
We will apply Yoshida~\cite{Yos11} to derive the mighty convergence
with the limit being shifted Gaussian.
As results, we will obtain an asymptotically normally distributed estimator
at rate $\sqrt{T_{n}}$ for both $\al$ and $\beta$ and also, very importantly,
the convergence of their moments to the corresponding ones of the limit
centered Gaussian distribution.
Different from the diffusion case, the GQL does not lead to an
asymptotically efficient estimator
in the presence of jumps, and is not even rate-efficient for $\beta$:
for instance, in the case where $X$ is a diffusion with
compound-Poisson jumps,
the information loss in the GQMLE of $\al$ can be large if the jump
part is much larger than the diffusion part; see Section~\ref{sec_discussion_eff}.
That is to say, the performance of our GQMLE may strongly depend on the
structure of the jump part and
its relation to the possibly nondegenerate diffusion one,
which may be considered as a possible major drawback of our estimation
procedure.
Nevertheless, it has the practically important advantages:
first, the computation of estimates does not require any fine tuning,
hence is straightforward;
second, the estimation procedure can be adopted without full
specification of the L\'evy measure $\nu$.
Further, our numerical experiments in Section~\ref{sec_sim} reveal
that, when the diffusion part is absent,
it can happen that the finite-sample performance of $\hat{\theta}_{n}$
becomes as good as the diffusion case if~$J$ ``distributionally'' close
to the Wiener process.

We should mention that the convergence of moments especially serves as
a fundamental tool when analyzing
asymptotic behavior of the expectations of statistics depending on the
estimator, for example,
asymptotic bias and mean squared prediction error, model-selection
devices (information criteria) and
remainder estimation in higher-order inference.
In the past, several authors have investigated such a strong mode of
convergence of estimators;
see Bhansali and Papangelou~\cite{BhaPap91}, Chan and Ing \cite
{ChaIng11}, Findley and Wei~\cite{FinWei02},
Inagaki and Ogata~\cite{InaOga75}, Jeganathan~\cite{Jeg82,Jeg89},
Ogata and Inagaki~\cite{OgaIna77}, Sieders and Dzhaparidze \cite
{SieDzh87} and Uchida~\cite{Uch10},
as well as Ibragimov and Has'minski~\cite{IbrHas81}, Kutoyants \cite
{Kut84,Kut04} and Yoshida~\cite{Yos11}.
See also the recent paper Uchida and Yoshida~\cite{UY_adapt} for an
adaptive parametric estimation
of diffusions with moment convergence of estimators under the sampling
design $nh_{n}^{k}\to0$
for arbitrary integer $k\ge2$.

The rest of this paper is organized as follows.
Section~\ref{sec_GQMLE} introduces our GQL random field and presents
its asymptotic behavior,
together with a small numerical example for observing finite-sample
performance of the GQMLE.
Section~\ref{sec_MC} provides a somewhat general result concerning the
mighty convergence, based on which
we prove our main result in Section~\ref{sec_proofs}. In Section~\ref{sec_proof_prop_ergo},
we prove a fairly simple criterion for the exponential ergodicity
assumption in dimension one,
only in terms of the coefficient $(a,b,c)$ and the L\'evy measure $\nu(dz)$.

Throughout this paper, asymptotics are taken for $n\to\infty$ unless
otherwise mentioned,
and the following notation is used:

\begin{itemize}
\item$I_{r}$ denotes the $r\times r$-identity matrix;

\item given a multilinear form $M=\{M^{(i_{1}i_{2}\cdots i_{K})}\dvtx i_{k}=1,\ldots,d_{k}; k=1,\ldots,K\}
\in\mbbr^{d_{1}}\otimes\cdots\otimes\mbbr^{d_{K}}$
and variables $u_{k}=\{u_{k}^{(i)}\}_{i\le d_{k}}\in\mbbr^{d_{k}}$,
we write
\[
M[u_{1},\ldots,u_{K}]=\sum_{i_{1}=1}^{d_{1}}
\cdots \sum_{i_{K}=1}^{d_{K}}M^{(i_{1}i_{2}\cdots i_{K})}u_{1}^{(i_{1})}
\cdots u_{K}^{(i_{K})}.
\]
The correspondences of indices of $M$ and $u_{k}$ will be clear from
each context.
Some of $u_{k}$ may be missing in ``$M[u_{1},\ldots,u_{K}]$''
so that the resulting form again defines a multilinear form,
for example, $M[u_{3},\ldots,u_{K}]\in\mbbr^{d_{1}}\otimes\mbbr^{d_{2}}$.
In particular, given two multilinear forms
$M^{(j)}=\{M^{(i_{1}i_{2}\cdots i_{K(j)})}\}$, $j=1,2$, we often use
the notation $M^{(1)}\otimes M^{(2)}$ for the tensor product
\begin{eqnarray*}
&& \bigl(M^{(1)}\otimes M^{(2)}\bigr)[u_{1},
\ldots,u_{K(1)},v_{1},\ldots,v_{K(2)}]
\\
&&\qquad:=\bigl(M^{(1)}[u_{1},\ldots,u_{K(1)}]\bigr)
\bigl(M^{(2)}[v_{1},\ldots,v_{K(2)}]\bigr).
\end{eqnarray*}
%
%(M^{(1)}\otimes M^{(2)})[u_{1},\cdots,u_{K(1)},v_{1},\cdots,v_{K(2)}]
%:=(M^{(1)}[u_{1},\cdots,u_{K(1)}])(M^{(2)}[v_{1},\cdots,v_{K(2)}]).
When $K\le2$, identifying $M$ as a vector or matrix,
we write $M^{\otimes2}=MM^{\top}$ with $\top$ denoting the transpose;
furthermore, $|M|$ denotes either, depending on the context, $\operatorname{
det}(M)$ when $d_{1}=d_{2}$,
or any matrix norm of $M$.

\item$\p_{a}^{m}$ stands for the bundled $m$th partial differential
operator with respect to $a=\{a^{(i)}\}$.

\item$C$ denotes generic positive constant possibly varying from line
to line,
and we write $x_{n}\lesssim y_{n}$ if $x_{n}\le Cy_{n}$ a.s. for every
$n$ large enough.
\end{itemize}

%%%%%
%%%%%

%s2 #&#
\section{Gaussian quasi-likelihood estimation}\label{sec_GQMLE}

We denote by $(\Omega,\mcf,\mbF=\break(\mcf_{t})_{t\in\mbbrp},P)$ a complete
filtered probability space
on which the process $X$ given by (\ref{SDE}) is defined:
the initial variable $X_{0}$ being $\mcf_{0}$-measurable, and $(W,J)$
being $\mbF$-adapted.

%%%

%s2.1 #&#
\subsection{Assumptions}

%as2.1 #&#
\begin{ass}[(Moments)]\label{A_J}
$E[J_{1}]=0$, $E[J_{1}^{\otimes2}]=I_{r''}$, and
$E[|J_{1}|^{q}]<\infty
$ for all $q>0$.
\end{ass}

We introduce the function $V\dvtx \mbbr^{d}\times\Theta_{\beta}\to\mbbr
^{d}\otimes\mbbr^{d}$ by
\[
V=b^{\otimes2}+c^{\otimes2}.
\]
For each $\theta$, the function $x\mapsto V(x,\beta)$ can be viewed as
an approximate local covariance matrix of
the law of $h_{n}^{-1/2}(X_{h_{n}}-x)$ under $P_{\theta}[\cdot| X_{0}=x]$.

Let $\overline{\Theta}$ denote the closure of $\Theta$.

%as2.2 #&#
\begin{ass}[(Smoothness)]\label{A_coeff}
(a)
The coefficient $(a,b,c)$ has the extension in $\mcc(\mbbr^{d}\times
\overline{\Theta})$,
and has partial derivatives such that $(\p_{\al}a,\p_{\beta}b,\p
_{\beta}c)$
admits the extension in $\mcc(\mbbr^{d}\times\overline{\Theta})$, that
\[
\sup_{(x,\theta)\in\mbbr^{d}\times\Theta} \bigl\{\bigl|\p_{x}a(x,\al )\bigr|+\bigl|\p
_{x}b(x,\beta)\bigr|+\bigl|\p_{x}c(x,\beta)\bigr| \bigr\}<\infty,
\]
and that, for each $k\in\{0,1,2\}$ and $l\in\{0,1,\ldots,5\}$,
there exists a constant $C(k,l)\ge0$ for which
\[
\sup_{(x,\theta)\in\mbbr^{d}\times\Theta}\bigl(1+|x|\bigr)^{-C(k,l)} \bigl\{ \bigl|
\p_{x}^{k}\,\p_{\al}^{l}a(x,\al)\bigr|+\bigl|
\p_{x}^{k}\,\p_{\beta
}^{l}b(x,\beta )\bigr|+\bigl|
\p_{x}^{k}\,\p_{\beta}^{l}c(x,\beta)\bigr|
\bigr\}<\infty.
\]
\begin{longlist}[(a)]
\item[(b)] $V(x,\beta)$ is invertible for each $(x,\beta)$, and there
exists a constant $C(V)\ge0$ such that
\[
\sup_{(x,\beta)\in\mbbr^{d}\times\Theta_{\beta
}}\bigl(1+|x|\bigr)^{-C(V)}\bigl|V^{-1}(x,\beta)\bigr|<
\infty.
\]
\end{longlist}
\end{ass}

When considering large-time asymptotics, the stability property of $X$
much affects the statistical analysis in essential ways. A typical
situation to be considered is that $X$ is ergodic.
We impose here a stronger stability condition.
Let $(P_{t})$ denote the transition semigroup of $X$.
Given a function $\rho\dvtx \mbbr^{d}\to\mbbrp$ and a signed measure $m$ on
the $d$-dimensional Borel space,
we define
\[
\|m\|_{\rho}=\sup \bigl\{\bigl|m(f)\bigr|\dvtx \mbox{$f$ is $\mbbr$-valued and
measurable, and fulfils that $|f|\le\rho$} \bigr\}.
\]

%as2.3 #&#
\begin{ass}[(Stability)]\label{A_ergo}
(a)
There exists a probability measure $\pi_{0}$ such that for every $q>0$
we can find a constant $a>0$ for which
%
%
%e2.1 #&#
\begin{equation}
\sup_{t\in\mbbrp}e^{at}\bigl\|P_{t}(x,\cdot)-
\pi_{0}(\cdot)\bigr\| _{g}\lesssim g(x),\qquad x\in
\mbbr^{d}, \label{exp_ergo_def}
\end{equation}
where $g(x):=1+|x|^{q}$.
\begin{longlist}[(b)]
\item[(b)] For every $q>0$,
%
%
%e2.2 #&#
\begin{equation}
\sup_{t\in\mbbrp}E_{0}\bigl[|X_{t}|^{q}
\bigr]<\infty. \label{g_moment_bound}
\end{equation}
\end{longlist}
\end{ass}

Here and in the sequel, $E_{0}$ denotes the expectation operator with
respect to $P_{0}$.
Condition (\ref{exp_ergo_def}) with $g$ replaced by the constant $1$
is the \textit{exponential ergodicity}, which in particular entails the
ergodic theorem:
the limit $\pi_{0}$ is a unique invariant distribution such that, for
every $f\in L^{1}(\pi_{0})$,
%
%
%e2.3 #&#
\begin{equation}
\frac{1}{T_{n}}\int_{0}^{T_{n}}f(X_{t})\,dt
\to^{p}\int f(x)\pi_{0}(dx), \label{LLN}
\end{equation}
where $\to^{p}$ stands for the convergence in $P_{0}$-probability;
we see that %$n^{-1}\sum_{j=1}^{n}f(X_{t_{j-1}})\to^{p}\int f(x)
\[
\frac{1}{n}\sum_{j=1}^{n}f(X_{t_{j-1}})
\to^{p}\int f(x)\pi_{0}(dx)
\]
for continuously differentiable $f$ with $\p f$ at most polynomial
order, since
%
%
%e2.4 #&#
\begin{eqnarray}\label{LLN_gap}
&& E_{0} \Biggl[ \Biggl|\frac{1}{T_{n}}\int_{0}^{T_{n}}f(X_{t})\,dt
-\frac{1}{n}\sum_{j=1}^{n}f(X_{t_{j-1}})
\Biggr| \Biggr]
\nonumber
\\[-8pt]
\\[-8pt]
\nonumber
&&\qquad\lesssim\frac{1}{n}\sum_{j=1}^{n}
\sup_{t_{j-1}\le s\le t_{j}} \sqrt{E_{0}\bigl[|X_{s}-X_{t_{j-1}}|^{2}
\bigr]}\to0.
\end{eqnarray}
We also note that Assumption~\ref{A_ergo} entails the \textit
{exponential absolute regularity},
also referred to as the \textit{exponential $\beta$-mixing property}.
This means that $\beta_{X}(t)=O(e^{-at})$ as $t\to\infty$ for some
$a>0$, where
$\beta_{X}$ denotes the $\beta$-mixing coefficient
\[
\beta_{X}(t):=\sup_{s\in\mbbrp}\int\bigl\|P_{t}(x,
\cdot)-\eta P_{s+t}(\cdot )\bigr\|\eta P_{s}(dx),
\]
where $\eta P_{t}:=\mcl(X_{t})$ and $\|m\|:=\|m\|_{1}$.
Let us recall that the exponential absolute regularity implies the
\textit{exponential strong-mixing property},
which plays an essential role in Yoshida~\cite{Yos11}, Lemma 4,
which we will apply in the proof of Theorem~\ref{thmm_main}.

Several sufficient conditions for Assumption~\ref{A_ergo} are known;
for diffusion processes,
see the references of Masuda~\cite{Mas07,Mas08} for some details. In
the presence of the jump component,
verification of (\ref{exp_ergo_def}) can become much more involved.
Especially if the coefficients are nonlinear and the L\'evy process $J$
is of infinite variation,
the verification may be far from being a trivial matter.
We refer to Kulik~\cite{Kul09,Kul11}, Maruyama and Tanaka~\cite{MarTan59},
Menaldi and Robin~\cite{MenRob99}, Meyn and Tweedie~\cite{mt3}
and Wang~\cite{Wan08} as well as Masuda~\cite{Mas07,Mas08}
for some general results concerning the (exponential) ergodicity.
For the sake of convenience, focusing on the univariate case and
setting ease of verification above generality,
we will provide in Proposition~\ref{prop_ergo} sufficient conditions
for Assumption~\ref{A_ergo},
in a form enabling us to deal with cases of nonlinear coefficients and
infinite-variation $J$;
see also Remark~\ref{rem_ergo_2}.

Define $\mbbg_{\infty}(\theta)=(\mbbg_{\infty}^{\al}(\theta
),\mbbg
_{\infty}^{\beta}(\beta))\in\mbbr^{p}$ by
%
%
%e2.5 #&#
%e2.6 #&#
\begin{eqnarray}
\mbbg_{\infty}^{\al}(\theta)&=&\int\p_{\al}a(x,\al)
\bigl[V^{-1}(x,\beta) \bigl[a(x,\al_{0})-a(x,\al) \bigr]
\bigr]\pi_{0}(dx), \label{Ga_infty}
\\
\mbbg_{\infty}^{\beta}(\beta)&=&\int \bigl\{V^{-1}(
\p_{\beta}V)V^{-1}(x,\beta) \bigr\} \bigl[V(x,
\beta_{0})-V(x,\beta) \bigr]\pi_{0}(dx).
\label{Gb_infty}
\end{eqnarray}
[In (\ref{Gb_infty}), we regarded ``$V^{-1}(\p_{\beta
}V)V^{-1}(x,\beta
)$'' as
a bilinear form with dimensions of indices being $p_{\beta}$ and $d^{2}$.]
Further, let $\mbbg_{\infty}'(\theta_{0}):=\operatorname{ diag}\{\mbbg
_{\infty
}^{\prime\al}(\theta_{0}),\break
\mbbg_{\infty}^{\prime\beta}(\theta_{0})\}\in\mbbr^{p}\otimes
\mbbr^{p}$,
where, for each $v'_{1}, v'_{2}\in\mbbr^{p_{\al}}$ and $v''_{1},
v''_{2}\in\mbbr^{p_{\beta}}$,
%
%
%e2.7 #&#
%e2.8 #&#
\begin{eqnarray}
 \label{G'a}&&\mbbg_{\infty}^{\prime\al}(\theta_{0})\bigl[v'_{1},v'_{2}
\bigr]
\nonumber
\\[-8pt]
\\[-8pt]
\nonumber
&&\qquad=-\int V^{-1}(x,\beta_{0}) \bigl[
\p_{\al}a(x,\al_{0})\bigl[v'_{1}
\bigr], \p _{\al
}a(x,\al_{0})\bigl[v'_{2}
\bigr] \bigr]\pi_{0}(dx),
\\
\label{G'b}&&\mbbg_{\infty}^{\prime\beta}(\theta_{0})\bigl[v''_{1},v''_{2}
\bigr]
\nonumber
\\[-8pt]
\\[-8pt]
\nonumber
&&\qquad=-\int\operatorname{ trace} \bigl[ \bigl\{\bigl(V^{-1}\p_{\beta}V
\bigr)\otimes \bigl(V^{-1}\p _{\beta}V\bigr) \bigr\}(x,
\beta_{0}) \bigl[v''_{1},v''_{2}
\bigr] \bigr]\pi_{0}(dx).
\end{eqnarray}

%as2.4 #&#
\begin{ass}[(Identifiability)]\label{A_iden}
There exist positive constants $\chi_{\al}=\chi_{\al}(\theta_{0})$ and
$\chi_{\beta}=\chi_{\beta}(\theta_{0})$
such that $|\mbbg_{\infty}^{\al}(\theta)|^{2}\ge\chi_{\al}|\al
-\al
_{0}|^{2}$ and
$|\mbbg_{\infty}^{\beta}(\beta)|^{2}\ge\chi_{\beta}|\beta-\beta
_{0}|^{2}$ for every $\theta\in\Theta$.
\end{ass}

%as2.5 #&#
\begin{ass}[(Nondegeneracy)]
Both $\mbbg_{\infty}^{\prime\al}(\theta_{0})$ and $\mbbg_{\infty
}^{\prime\beta}(\theta_{0})$ are invertible.
\label{A_nond}
\end{ass}

Assumptions~\ref{A_iden} and~\ref{A_nond} are quite typical in
statistical estimation.
In Lemma~\ref{lem_iden_nond} below,
both assumptions are implied by a kind of uniform nonsingularity.
Define two bilinear forms $\bar{A}(\al',\al'',\beta')$ and $\bar
{B}(\beta',\beta'')$ by,
just like (\ref{G'a}) and (\ref{G'b}),
\begin{eqnarray*}
&&\bar{A}\bigl(\al',\al'',
\beta'\bigr)\bigl[v'_{1},v'_{2}
\bigr] =\int V^{-1}\bigl(x,\beta'\bigr) \bigl[
\p_{\al}a\bigl(x,\al'\bigr)\bigl[v'_{1}
\bigr], \p_{\al
}a\bigl(x,\al ''\bigr)
\bigl[v'_{2}\bigr] \bigr]\pi_{0}(dx),
\\
&&\bar{B}\bigl(\beta',\beta''\bigr)
\bigl[v''_{1},v''_{2}
\bigr]\\
&&\qquad =\int\operatorname{ trace} \bigl[ \bigl\{\bigl(V^{-1}(
\p_{\beta}V)V^{-1}\bigr) \bigl(x,\beta '\bigr)
\otimes\p_{\beta}V\bigl(x,\beta''\bigr)
\bigl[v''_{1},v''_{2}
\bigr] \bigr\} \bigr]\pi_{0}(dx).
\end{eqnarray*}

%le2.6 #&#
\begin{lem}\label{lem_iden_nond}
Suppose that $\bar{A}(\al',\al'',\beta')$ and $\bar{B}(\beta
',\beta'')$ are
nonsingular uniformly in $\al',\al''\in\Theta_{\al}$ and $\beta
',\beta
''\in\Theta_{\beta}$.
Then both Assumptions~\ref{A_iden} and~\ref{A_nond} hold true.
\end{lem}

\begin{pf}
It is obvious that Assumption~\ref{A_nond} follows.
The mean-value theorem applied to (\ref{Ga_infty}) and (\ref{Gb_infty})
leads to
$\mbbg^{\al}_{\infty}(\theta)=\bar{A}(\al,\tilde{\al},\beta
)[\al_{0}-\al
]$ for some $\tilde{\al}$
lying the segment connecting $\al$ and $\al_{0}$, with a similar form
for $\mbbg^{\beta}_{\infty}(\beta)$;
recall that $\Theta_{\al}$ and $\Theta_{\beta}$ are presupposed to
be convex.
Since $\inf_{\al',\al'',\beta'} \Vert\bar{A}(\al',\al'',\break\beta')\Vert>0$
and $\inf_{\beta',\beta''} \Vert\bar{B}(\beta',\beta'')\Vert>0$ under
the assumption, the matrices
$\bar{A}^{\otimes2}$ and $\bar{B}^{\otimes2}$ are uniformly
positive definite,
hence Assumption~\ref{A_iden} follows.
\end{pf}

%%%

%s2.2 #&#
\subsection{Asymptotics: Main results}

In what follows, we write
\[
\D_{j}Y=Y_{t_{j}}-Y_{t_{j-1}}
\]
for any process $Y$, and
\[
f_{j-1}(a)=f(X_{t_{j-1}},a)
\]
for a variable $a$ in some set $A$ and a measurable function $f$ on
$\mbbr^{d}\times A$.
The Euler approximation for SDE (\ref{SDE}) is formally
\[
X_{t_{j}}\approx X_{t_{j-1}}+a_{j-1}(\al)h_{n}+b_{j-1}(
\beta)\D _{j}W+c_{j-1}(\beta)\D_{j}J
\]
under $P_{\theta}$, which leads us to consider the local-Gauss
distribution approximation
%
%
%e2.9 #&#
\begin{equation}
\mcl(X_{t_{j}}|X_{t_{j-1}})\approx\mcn_{d} \bigl(
X_{t_{j-1}}+a_{j-1}(\al)h_{n}, h_{n}V_{j-1}(
\beta) \bigr). \label{lGa}
\end{equation}
Put
\[
\chi_{j}(\al)=\D_{j}X-h_{n}a_{j-1}(
\al).\vadjust{\goodbreak}
\]
Based on (\ref{lGa}), we define our GQL by
%
%
%e2.10 #&#
\begin{equation}
\mbbq_{n}(\theta)=-\sum_{j=1}^{n}
\biggl\{ \log\bigl|V_{j-1}(\beta)\bigr|+\frac{1}{h_{n}}V_{j-1}^{-1}(
\beta) \bigl[\chi _{j}(\al)^{\otimes2} \bigr] \biggr\},
\label{GQL}
\end{equation}
and the corresponding GQMLE by any element
\[
\hat{\theta}_{n}=(\hat{\al}_{n},\hat{\beta}_{n})
\in \mathop{\operatorname{ argmax}}_{\theta\in\overline{\Theta}}\mbbq_{n}(\theta).
\]

Under Assumption~\ref{A_J} we have $\int z^{(k)}z^{(l)}\nu(dz)=\del_{kl}$
for $k,l\in\{1,\ldots,r''\}$. We need some further notation in this direction.
For $i_{1},\ldots,i_{m}\in\{1,\ldots,r''\}$ with \mbox{$m\ge3$},
we write $\nu(m)$ for the $m$th mixed moments of $\nu$,
\[
\nu(m)=\bigl\{\nu_{i_{1}\cdots i_{m}}(m)\bigr\}_{i_{1},\ldots,i_{m}}:= \biggl\{\int
z^{(i_{1})}\cdots z^{(i_{m})}\nu(dz) \biggr\} _{i_{1},\ldots,i_{m}}.
\]
Let $c^{(\cdot k)}(x,\beta)\in\mbbr^{d}$ denote the $k$th column of
$c(x,\beta)$.
We introduce the matrix
%
%
%e2.11 #&#
\begin{equation}
\mbbv(\theta_{0}):=\pmatrix{
\mbbg^{\prime\al}_{\infty}(\theta_{0}) &
\mbbv_{\al\beta}
\vspace*{2pt}\cr
\mbbv_{\al\beta}^{\top} & \mbbv_{\beta\beta}}, \label{mbbv_def}
\end{equation}
where, for each $v'\in\mbbr^{p_{\al}}$ and $v_{1}'', v''_{2}\in
\mbbr
^{p_{\beta}}$,
\begin{eqnarray*}
\mbbv_{\al\beta}\bigl[v',v''_{1}
\bigr]&:=&-\int\sum_{k',l',s'}\nu_{k'l's'}(3)
V^{-1}(x,\beta_{0}) \bigl[\p_{\al}a(x,
\al_{0})\bigl[v'\bigr],c^{(\cdot
s')}(x,\beta
_{0}) \bigr]
\\
&&\hspace*{42pt}{} \times \bigl\{\p_{\beta}V^{-1}(x,\beta_{0})
\bigr\} \bigl[v''_{1},c^{(\cdot k')}(x,
\beta_{0}),c^{(\cdot l')}(x,\beta_{0})\bigr]\pi
_{0}(dx),
\\
\mbbv_{\beta\beta}\bigl[v''_{1},v''_{2}
\bigr] &:=&\int\sum_{s,t,s',t'}\nu_{sts't'}(4) \bigl
\{\p_{\beta}V^{-1}(x,\beta_{0})\bigl[v''_{1},c^{(\cdot s)}(x,
\beta _{0}),c^{(\cdot t)}(x,\beta_{0})\bigr] \bigr\}
\\
&&\hspace*{36pt}{} \times \bigl\{\p_{\beta}V^{-1}(x,\beta_{0})
\bigl[v''_{2},c^{(\cdot s')}(x,\beta
_{0}),c^{(\cdot t')}(x,\beta_{0})\bigr] \bigr\}
\pi_{0}(dx).
\end{eqnarray*}
Finally, put
\[
\Sig_{0}=\pmatrix{ \bigl(-
\mbbg_{\infty}^{\prime\al}\bigr)^{-1}(\theta_{0}) &
\bigl\{\bigl(\mbbg_{\infty}^{\prime\al}\bigr)^{-1}
\mbbv_{\al\beta}\bigl(\mbbg _{\infty
}^{\prime\beta}
\bigr)^{-1}\bigr\}(\theta_{0})
\vspace*{2pt}\cr
\operatorname{ Sym.} & \bigl\{\bigl(\mbbg_{\infty}^{\prime\beta}
\bigr)^{-1}\mbbv_{\beta\beta}\bigl(\mbbg _{\infty
}^{\prime\beta}
\bigr)^{-1}\bigr\}(\theta_{0}) }.
\]
Now we can state our main result, the proof of which is deferred to
Section~\ref{sec_proof1}.

%th2.7 #&#
\begin{thmm}\label{thmm_main}
Suppose Conditions~\ref{A_J},~\ref{A_coeff},~\ref{A_ergo},~\ref{A_iden}
and~\ref{A_nond}.
Then we have
\[
E_{0} \bigl[f \bigl(\sqrt{T_{n}}(\hat{
\theta}_{n}-\theta_{0}) \bigr) \bigr] \to\int f(u)\phi
(u;0,\Sig_{0} )\,du, \qquad n\to\infty
\]
for every continuous function $f\dvtx \mbbr^{p}\to\mbbr$ of at most
polynomial growth,
where $\phi(\cdot;0,\Sig_{0})$ denotes the centered Gaussian density
with covariance matrix $\Sig_{0}$.
\end{thmm}

The following two remarks are immediate:
\begin{itemize}
\item The estimators $\hat{\al}_{n}$ and $\hat{\beta}_{n}$ are
asymptotically independent if $\nu(3)=0$,
implying that $\hat{\al}_{n}$ and $\hat{\beta}_{n}$ may not be
asymptotically independent if $\nu$ is skewed.
If $c\equiv0$ so that $X$ is a diffusion, then $\nu(4)=0$, so that
$\mbbv_{\beta\beta}=0$ and
$\sqrt{T_{n}}(\hat{\beta}_{n}-\beta_{0})$ is asymptotically degenerate
at $0$.
This is in accordance with the case of diffusion, where the GQMLE of
$\beta$ is $\sqrt{n}$-consistent.
See Section~\ref{sec_discussion_eff} for a discussion on the
efficiency issue.

\item The revealed convergence rate $\sqrt{T_{n}}$ of the GQMLE $\hat
{\beta}_{n}$
alerts us to take precautions against the presence of jumps.
For instance, suppose that one has adopted the parametric diffusion
model [i.e., (\ref{SDE}) with $c\equiv0$]
although there actually does exist a nonnull jump part.
Then one takes $\sqrt{n}$ for the convergence rate of~$\hat{\beta}_{n}$,
although the true one is $\sqrt{T_{n}}$, which may lead to a seriously
inappropriate confidence zone.
This point can be sufficient grounds for importance of testing the
presence of jumps.
In case of one-dimensional $X$, Masuda~\cite{Mas_rims}, Section~4, constructed
an analogue to Jarque--Bera normality test and studied its asymptotic behavior.
See Masuda~\cite{SNRP} for a multivariate extension.
\end{itemize}

In order to construct confidence regions for $\theta_{0}$ as well as to
perform statistical tests,
we need a consistent estimator of the asymptotic covariance matrix
$\Sig_{0}$.
Although $\Sig_{0}$ contains unknown third and fourth mixed moments of
$\nu$,
it turns out to be possible to provide a consistent estimator of $\Sig_{0}$
without any specific knowledge of $\nu$ other than Assumption \ref
{A_J}. Let\looseness=1
\[
\hat{\Sig}_{n}=\pmatrix{ \bigl(-
\hat{\mbbg}_{n}^{\prime\al}\bigr)^{-1} & \bigl(\hat{
\mbbg}_{n}^{\prime\al}\bigr)^{-1}\hat{
\mbbv}_{\al\beta,n}\bigl(\hat {\mbbg }_{n}^{\prime\beta}
\bigr)^{-1}
\vspace*{2pt}\cr
\operatorname{ Sym.} & \bigl(\hat{\mbbg}_{n}^{\prime\beta}
\bigr)^{-1}\hat{\mbbv}_{\beta\beta,n}\bigl(\hat {\mbbg}_{n}^{\prime\beta}
\bigr)^{-1} },
\]\looseness=0
where, for each $v'_{1}, v'_{2}\in\mbbr^{p_{\al}}$ and $v_{1}'',
v''_{2}\in\mbbr^{p_{\beta}}$,
\begin{eqnarray*}
&&\hat{\mbbg}_{n}^{\prime\al}\bigl[v'_{1},v'_{2}
\bigr]:=-\frac{1}{n}\sum_{j=1}^{n}V_{j-1}^{-1}(
\hat{\beta}_{n}) \bigl[\p_{\al}a_{j-1}(\hat{
\al}_{n})\bigl[v'_{1}\bigr],
\p_{\al
}a_{j-1}(\hat{\al }_{n})\bigl[v'_{2}
\bigr] \bigr],
\\[2pt]
&&\hat{\mbbg}_{n}^{\prime\beta}\bigl[v''_{1},v''_{2}
\bigr]:=-\frac{1}{n}\sum_{j=1}^{n}
\operatorname{ trace} \bigl\{ \bigl(V_{j-1}^{-1}\p_{\beta}V_{j-1}
\bigr)\otimes \bigl(V_{j-1}^{-1}\,\p_{\beta}V_{j-1}
\bigr) ) (\hat{\beta}_{n})\bigl[v''_{1},
v''_{2}\bigr] \bigr\},
\\[2pt]
&&\hat{\mbbv}_{\al\beta,n}\bigl[v'_{1},v''_{1}
\bigr]\\[2pt]
&&\qquad:= -\sum_{j=1}^{n}\frac{1}{T_{n}}
\bigl(V_{j-1}^{-1}\otimes\p_{\beta
}V_{j-1}^{-1}
\bigr) (\hat{\beta}_{n})\\[2pt]
&&\hspace*{61pt}{}\times \bigl[ \bigl(\p_{\al}a_{j-1}(
\hat{\al}_{n})\bigl[v'_{1}\bigr],
\chi_{j}(\hat {\al }_{n}) \bigr), \bigl(v''_{1},
\chi_{j}(\hat{\al}_{n})^{\otimes2}\bigr) \bigr],
\\[5pt]
&&\hat{\mbbv}_{\beta\beta,n}\bigl[v''_{1},v''_{2}
\bigr]\\[2pt]
&&\qquad:= \sum_{j=1}^{n}\frac{1}{T_{n}}
\bigl(\p_{\beta}V_{j-1}^{-1}\otimes \,\p
_{\beta}V_{j-1}^{-1} \bigr) (\hat{\beta}_{n})
\bigl[\bigl(v''_{1},\chi_{j}(
\hat{\al}_{n})^{\otimes2}\bigr),\bigl(v''_{2},
\chi _{j}(\hat{\al}_{n})^{\otimes2}\bigr) \bigr].
\end{eqnarray*}
We will denote by $\cil$ the weak convergence under $P_{0}$.
%
%co2.8 #&#
\begin{cor}\label{cor_main}
Under the conditions of Theorem~\ref{thmm_main}, we have $\hat{\Sig
}_{n}\to^{p}\Sig_{0}$, and hence
%
%
%e2.12 #&#
\begin{equation}
\hat{\Sig}_{n}^{-1/2}\sqrt{T_{n}}(\hat{
\theta}_{n}-\theta _{0})\cil\mcn _{p}(0,I_{p})
\label{rev_Sclt}
\end{equation}
holds true.
\end{cor}
The proof of Corollary~\ref{cor_main} is given in Section~\ref{sec_proof2}.

The primary objective of this paper is to derive
the $L^{q}(P_{0})$-boundedness of $\sqrt{T_{n}}(\hat{\theta
}_{n}-\theta
_{0})$ for every $q>0$,
for which the moment conditions [Assumptions~\ref{A_J} plus \ref
{A_ergo}(b)] seem indispensable.
Nevertheless, as pointed out by the anonymous referee,
the existence of the moments of all orders is too much to ask in
Corollary~\ref{cor_main}.
Let us discuss a possibility of relaxing the moment condition in some detail;
to make the exposition more clear, we here do not seek the greatest generality.

Clearly, the really necessary order (of $J$, hence $X$ too) partly
depends on the growth of
the coefficients $(a,b,c)$ and its partial derivatives with respect to
$\theta$.
We will show that the consistency and asymptotic normality of $\hat
{\theta}_{n}$
follow on some weaker moment and stability assumptions
than the corresponding ones imposed in Theorem~\ref{thmm_main}.
We impose the following three conditions instead of
Assumptions~\ref{A_coeff},~\ref{A_J} and~\ref{A_ergo}:
%
%
%e2.13 #&#
\begin{equation}
\cases{\displaystyle \max_{k\in\{0,1,2\} \atop l\in\{0,1,\ldots,5\}}\sup
_{(x,\theta
)\in\mbbr
^{d}\times\Theta} \bigl\{ \bigl|\p_{x}^{k}\,
\p_{\al}^{l}a(x,\al) \bigr| + \bigl|\p_{x}^{k}\,
\p_{\beta}^{l}b(x,\beta) \bigr| \vspace*{2pt}\cr
\hspace*{151pt}{}+ \bigl|\p_{x}^{k}
\,\p_{\beta}^{l}c(x,\beta) \bigr| \bigr\}<\infty,
\vspace*{2pt}\cr
\displaystyle {\sup_{(x,\theta)\in\mbbr^{d}\times\Theta} \bigl\llvert V^{-1}(x,\beta)\bigr
\rrvert <\infty;} }
\label{A_coeff+}
\end{equation}
%
%
%e2.14 #&#
%e2.15 #&#
\begin{eqnarray}\label{moment_order}
E[J_{1}]&=&0, \qquad E\bigl[J_{1}^{\otimes2}
\bigr]=I_{r''} \quad \mbox{and}
\nonumber
\\[-8pt]
\\[-8pt]
\nonumber
E\bigl[|J_{1}|^{q}\bigr]&<&
\infty\qquad \mbox{for some }q>(p\vee4);
\end{eqnarray}
\begin{equation}
\begin{tabular}{p{280pt}@{}}
$X$ admits a unique invariant distribution $\pi_{0}$
such that (\ref{LLN}) holds true for every $f\in L^{1}(
\pi_{0})$. \label{LLN-}
\end{tabular}
\end{equation}
It is possible to deal with unbounded coefficients,
but then we inevitably need the uniform boundedness of moments as in
(\ref{g_moment_bound}),
where the minimal value of the index~$q$ must be determined according
to the growth orders of all the
coefficients as well as their partial derivatives, leading to a
somewhat messy description.

We then derive the asymptotic normality result as follows, proof of
which is given in Section~\ref{sec_proof3}.

%th2.9 #&#
\begin{thmm}\label{thmm_AN}
Suppose (\ref{A_coeff+}), (\ref{moment_order}), (\ref{LLN-})
and Assumptions~\ref{A_iden}\break and~\ref{A_nond}. Then we have
$\sqrt{T_{n}}(\hat{\theta}_{n}-\theta_{0})\cil\mcn_{p}(0,\Sig_{0})$.
\end{thmm}

In particular, we then do not need the exponential mixing property in
Assumption~\ref{A_ergo}, and
the ergodic theorem (\ref{LLN}) is enough.
This is of great advantage, as the exponential ergodicity is much
stronger than (\ref{LLN}) to hold;
see also Remark~\ref{rem_ergo_2}.
Finally, it also should be noted that it is possible to derive the
Studentized version (\ref{rev_Sclt})
under the assumptions in Theorem~\ref{thmm_AN} with ``$q>(p\vee4)$''
in (\ref{moment_order}) strengthened to ``$q>(p\vee8)$.''
Indeed, it is clear from the proof of Corollary~\ref{cor_main} why we
require that $q>(p\vee8)$,
and we omit the details.

We end this section with some remarks on the model setup.
\begin{itemize}
\item Although we are considering ``ergodic'' $X$, it is obvious that
we can target L\'evy processes as well,
according to the built-in independence of the increments $(\D
_{j}X)_{j\le n}$.

\item A general form of the martingale estimating functions is
\[
\theta\mapsto\sum_{j=1}^{n}W_{j-1}(
\theta) \bigl\{g(X_{t_{j-1}},X_{t_{j}};\theta)-E_{\theta}
\bigl[g(X_{t_{j-1}},X_{t_{j}};\theta)|\mcf_{t_{j-1}} \bigr]
\bigr\}
\]
for some $W\in\mbbr^{p}\otimes\mbbr^{m}$ and
$\mbbr^{m}$-valued function $g$ on $\mbbr^{d}\times\mbbr^{d}\times
\Theta$.
We would have a wide choice of $W$ and $g$.
When the conditional expectations involved do not admit closed forms,
then the leading-term approximation of them via the It\^o--Taylor
expansion can be used.
In view of this, as in Kessler~\cite{Kes97},
it would be formally possible to relax the condition $nh_{n}^{2}\to0$
in (\ref{RIED}) by gaining the order of
the It\^o--Taylor expansions of the conditional mean and conditional covariance,
\begin{eqnarray*}
E_{\theta}[X_{t_{j}}|\mcf_{t_{j-1}}]&=&X_{t_{j-1}}+a_{j-1}(
\al )h_{n}+\cdots,
\\
V_{\theta}[X_{t_{j}}|\mcf_{t_{j-1}}]&=&V_{j-1}(
\beta)h_{n}+\cdots,
\end{eqnarray*}
which we have implicitly used up to the $h_{n}$-order terms to build
$\mbbq_{n}$ of (\ref{GQL}).
However, we then need specific moment structures of $\nu$,
which appear in the higher orders of the above It\^o--Taylor expansion.
Moreover, we should note that the convergence rate $\sqrt{T_{n}}$ can
never be improved for both $\al$ and $\beta$,
even if $E_{\theta}[X_{t_{j}}|\mcf_{t_{j-1}}]$ and $V_{\theta
}[X_{t_{j}}|\mcf_{t_{j-1}}]$ have
closed forms, such as the case of linear drifts,
so that the rate of $h_{n}\to0$ may not matter as long as $T_{n}\to
\infty$.
See also Remark~\ref{rem_G_exte}.

\item As was mentioned in the \hyperref[sec1]{Introduction},
the sampling points $t_{1},\ldots,t_{n}$ may be irregularly spaced to
some extent.
Let $0\equiv t_{0}<t_{1}<\cdots<t_{n}=:T_{n}$, and put $\D
_{j}t:=t_{j}-t_{j-1}$.
We claim that it is possible to remove the equidistance condition,
while retaining that $h_{n}:=\max_{1\le j\le n}\D_{j}t\to0$.
We need the additional condition about asymptotic behavior of the spacing
%
%
%e2.16 #&#
\begin{equation}
\frac{1}{h_{n}}\min_{1\le j\le n}\D_{j}t\to1,
\label{irr_sampling}
\end{equation}
which obviously entails that $T_{n}\sim nh_{n}$ (the ratio of both
sides tends to~$1$).
Then the same statements as in Theorem~\ref{thmm_main}, Corollary \ref
{cor_main}
and Theorem~\ref{thmm_AN} remain valid under (\ref{irr_sampling}).
For this point, we only note that estimate (\ref{LLN_gap}) remains true
even under (\ref{irr_sampling}):
noting that
\begin{eqnarray*}
k_{n}&:=&\max_{j\le n}\biggl\llvert \biggl(
\frac{1}{n\D_{j}t}-\frac
{1}{T_{n}} \biggr)n\D_{j}t\biggr\rrvert
\\[-2pt]
&\le& \biggl(1-\frac{1}{h_{n}}\min_{j\le n}\D_{j}t
\biggr) + \biggl(\frac{nh_{n}}{T_{n}}-1 \biggr)=o(1),
\end{eqnarray*}
we have, for any $f$ such that both $f$ and $\p f$ are of at most
polynomial growth,
\begin{eqnarray*}
\del_{n}&:=&\Biggl\llvert \frac{1}{T_{n}}\int_{0}^{T_{n}}f(X_{t})\,dt-
\frac
{1}{n}\sum_{j=1}^{n}f(X_{t_{j-1}})
\Biggr\rrvert
\\[-2pt]
&= &\Biggl\llvert \sum_{j=1}^{n}
\frac{1}{T_{n}}\int_{t_{j-1}}^{t_{j}}f(X_{t})\,dt
-\sum_{j=1}^{n}\frac{1}{n\D_{j}t}\int
_{t_{j-1}}^{t_{j}}f(X_{t_{j-1}})\,dt\Biggr\rrvert
\\[-2pt]
&\le &k_{n}\frac{1}{n}\sum_{j=1}^{n}
\frac{1}{\D_{j}t}\int_{t_{j-1}}^{t_{j}}\bigl|f(X_{t})\bigr|\,dt
\\[-2pt]
&&{} +\frac{1}{n}\sum_{j=1}^{n}
\frac{1}{\D_{j}t}\int_{t_{j-1}}^{t_{j}}\bigl|f(X_{t})-f(X_{t_{j-1}})\bigr|\,dt
\\[-2pt]
&\lesssim& k_{n}\frac{1}{n}\sum_{j=1}^{n}
\frac{1}{\D_{j}t}\int_{t_{j-1}}^{t_{j}}\bigl(1+|X_{t}|\bigr)^{C}\,dt
\\[-2pt]
&& {} +\frac{1}{n}\sum_{j=1}^{n}
\frac{1}{\D_{j}t}\int_{t_{j-1}}^{t_{j}}\bigl(1+|X_{t}|\bigr)^{C}|X_{t}-X_{t_{j-1}}|\,dt
\end{eqnarray*}
for some $C>0$. Therefore, Schwarz's inequality together with Lemma~\ref{lem_inc_ME} leads to
the estimate $E_{0}[\del_{n}]\lesssim k_{n}+\sqrt{h_{n}}=o(1)$,
enabling us to use
$n^{-1}\times \sum_{j=1}^{n}f(X_{t_{j-1}})\to^{p}\int f(x)\pi_{0}(dx)$
as in the case of the equally-spaced sample. With this in mind, we can deduce
the same estimates and limit results in the proofs given in Sections~\ref{sec_proof2} to~\ref{sec_proof3}
in an entirely analogous way, the details being omitted.\vadjust{\goodbreak}
\end{itemize}

%%%

%s2.3 #&#
\subsection{Discussion}\label{sec_discussion}

%s2.3.1 #&#
\subsubsection{On the identifiability of the dispersion parameter}

Suppose that the coefficients $b(x,\beta)$ and $c(x,\beta)$ depend on
$\beta$ only through $\beta_{1}$ and $\beta_{2}$, respectively, where
$\beta=(\beta_{1},\beta_{2})$.
On the one hand, it should be theoretically possible to identify $\beta
_{1}$ and $\beta_{2}$ individually
by the (intractable) likelihood function;
for example, see A\"{\i}t-Sahalia and Jacod~\cite{AitJac08} for the
precise asymptotic behavior
of the Fisher information matrix for $\beta$ in case of univariate L\'
evy processes.
We also refer to A\"{\i}t-Sahalia and Jacod~\cite{AitJac07} for how
to construct
an asymptotically efficient estimator of $\beta_{1}$
through the use of a truncated power-variation statistics,
regarding $\beta_{2}$ as a nuisance parameter.
To perform individual estimation for more general diffusions with jumps,
it is unadvised to resort to the likelihood based estimation.
Instead, we may adopt a threshold-type estimator utilizing
only relatively small (resp., large) increments of $X$ for estimating
$\beta_{1}$ (resp., $\beta_{2}$),
which makes it possible to extract information of the diffusion and
jump parts separately,
in compensation for a nontrivial fine tuning of the threshold;
see Shimizu and Yoshida~\cite{ShiYos06} and Ogihara and Yoshida \cite
{OgiYos11} in case of compound-Poisson jumps
and Shimizu~\cite{Shi06} in the presence of infinitely many small jumps
of finite variation.

On the other hand, our identifiability condition on $\beta$ in
Assumption~\ref{A_iden}
can be unfortunately stringent in the simultaneous presence of
nondegenerate diffusion and jump components.
Let us look at the assumption in the multiplicative-parameter case
$b(x,\beta)=\beta_{1}b_{0}(x)$ and $c(x,\beta)=\beta_{2}c_{0}(x)$,
where $b_{0}$ and $c_{0}$ are known positive functions and
where we set $d= r'=r''=p_{\beta}=1$ for simplicity;
we implicitly suppose that the function equals $1$ if it is constant
because the constant then can be absorbed into $\beta$.
Further, we here suppose that $\overline{\Theta}_{\beta}\subset
(0,\infty
)\times(0,\infty)$, so that
$X$ admits both nonnull diffusion and jump parts.
Then direct computation gives $\mbbg_{\infty}^{\beta}(\beta
)=M(\beta
)[\beta_{0}-\beta]$, where
\[
M(\beta):=\lleft( \matrix{ 2\beta_{1}(\beta_{10}+
\beta_{1})I_{bb} & 2\beta_{1}(\beta
_{20}+\beta_{2})I_{bc} \vspace*{2pt}
\cr
2
\beta_{2}(\beta_{10}+\beta_{1})I_{bc}
& 2\beta_{2}(\beta _{20}+\beta_{2})I_{cc}
} \rright)
\]
with $I_{bb}:=\int b_{0}^{4}(x)V^{-2}(x,\beta)\pi_{0}(dx)$,
$I_{bc}:=\int b_{0}^{2}(x)c_{0}^{2}(x)V^{-2}(x,\beta)\pi_{0}(dx)$, and
$I_{cc}:=\int c_{0}^{4}(x)V^{-2}(x,\beta)\pi_{0}(dx)$.
We have $|M(\beta)|=C(\beta)|I_{bb}I_{cc}-I_{bc}|$ for some constant
$C(\beta)$ depending on $\beta$
such that $\inf_{\beta}C(\beta)>0$, so that the identifiability
condition on $\beta$ is satisfied
if $|I_{bb}I_{cc}-I_{bc}|>0$.
In view of Schwarz's inequality, we always have $I_{bb}I_{cc}-I_{bc}\ge
0$, the equality holding only when
there exists an $r\in\mbbr$ such that $b_{0}(x)=rc_{0}(x)$ for every
$x\in\mbbr$.
That is, the GQMLE fails to be consistent as soon as $b_{0}$ and $c_{0}$
are proportional to each other;
especially if both $b_{0}$ and $c_{0}$ are constant (hence $1$, as was
presupposed), then
the GQMLE indeed cannot identify $\beta_{1}$ and $\beta_{2}$ individually,
for there do exist infinitely many $\beta=(\beta_{1},\beta_{2})$
such that
\[
V(x,\beta)-V(x,\beta_{0})=\bigl(\beta_{1}^{2}+
\beta_{2}^{2}\bigr)-\bigl(\beta _{10}^{2}+
\beta_{20}^{2}\bigr)=0
\]
for every $x$.
This seems to be unavoidable as our contrast function $\mbbm_{n}$ is
constructed solely
based on fitting local conditional mean and covariance matrix.
Although our estimation procedure cannot in general separate
information of
diffusion and jump variances, it should be noted that,
when both $b_{0}$ and $c_{0}$ are constant,
we may instead consistently estimate the ``local variance'' $\beta
_{1}^{2}+\beta_{2}^{2}$.

Finally, we remark that the identifiability condition
``$|\mbbg_{\infty}^{\beta}(\beta)|^{2}\ge\chi_{\beta}|\beta
-\beta_{0}|^{2}$''
becomes much simpler when we know that $b(\cdot,\cdot)\equiv0$ from
the very beginning;
then, in view of expression (\ref{Gb_infty}) and Assumption~\ref{A_coeff}(b),
it would suffice to have $|\p_{\beta}c^{2}(x,\beta)|>0$ over a domain.

%

%Our parametrization (\ref{SDE}) is an extension of the one for
%diffusion case where
%$c(\cdot,\cdot)\equiv0$. As a matter of fact, it could be formally
%possible to deal with the more complicated model
%dX_{t}=a(X_{t},\al)dt+b(X_{t},\beta)dW_{t}+\int c(X_{t-},z,\gam)\tilde{
%where $\tilde{\mu}(dt,dz):=\tilde{\mu}(dt,dz)-\nu(dz)dt$
%stands for the compensated Poisson random measure with intensity (L
%But then, the GQL function is not explicit unless we have rather
%specific knowledge
%of the moments of $\nu$. We do not consider this case here.

%%%

%s2.3.2 #&#
\subsubsection{On the asymptotic efficiency}\label{sec_discussion_eff}

The efficiency issue for model (\ref{SDE}) based on high-frequency sampling
is a difficult problem and has been left unsolved over the years,
which hinders us to do quantitative study on how much information loss
occurs on our GQMLE;
as a matter of fact, we do not know any Haj\'ek bound on the asymptotic
covariances
especially when $J$ is of infinite activity.
This general issue is beyond the scope of this paper, but instead we
give some remarks in this direction.

\begin{itemize}
\item
Overall, the amount of efficiency loss in using our GQMLE may strongly
depend on the structure of the jump part and on its relation to the
possibly nondegenerate diffusion part; this would be a major drawback
of our GQMLE.
We do know the theoretical minimal asymptotic covariance matrix
when $X$ is a diffusion with compound-Poisson jumps with nondegenerate
diffusion part, where, in particular,
the optimal rate of convergence in estimating $\al$ is $\sqrt{T_{n}}$,
achieved by our GQMLE $\hat{\al}_{n}$; for details,
see Shimizu and Yoshida~\cite{ShiYos06} and Ogihara and Yoshida \cite
{OgiYos11}, as well as the references therein.
In order to observe the effect of the jump part in estimation of $\al$
in a concise way,
let us look at the univariate $X$ given by
\[
dX_{t}=a(X_{t},\al)\,dt+b(X_{t})\,dW_{t}+c(X_{t-})\,dJ_{t},
\]
where $\al\in\mbbr$, $\inf_{x}b(x)\wedge\inf_{x}c(x)>0$, and $J$
is a
centered compound-Poisson process.
The asymptotic variance of $\hat{\al}_{n}$ is then given by the
inverse of
\[
-\mbbg_{\infty}^{\prime\al}(\theta_{0}) =\int \bigl
\{b^{2}(x)+c^{2}(x) \bigr\}^{-1} \bigl\{
\p_{\al}a(x,\al _{0}) \bigr\}^{2}
\pi_{0}(dx),
\]
while the minimal asymptotic variance of the asymptotically efficient
estimator is the inverse of
$A_{0}^{\ast}:=\int b^{-2}(x)\{\p_{\al}a(x,\al_{0})\}^{2}\pi_{0}(dx)$.
Hence, it would be natural to measure amount of efficiency loss in
using $\hat{\al}_{n}$ by the quantity
\[
A^{\ast}_{0}-\bigl\{-\mbbg_{\infty}^{\prime\al}(
\theta_{0})\bigr\} =\int\frac{\{\p_{\al}a(x,\al_{0})\}^{2}}{b^{2}(x)} \biggl(\frac{c^{2}(x)}{b^{2}(x)+c^{2}(x)}
\biggr)\pi_{0}(dx).
\]
From this expression, we may expect that the efficiency loss may be
large (resp., not so significant)
when the jump part is much larger (resp., smaller) compared with the
diffusion part.
This point comes into focus by looking at the Ornstein--Uhlenbeck process
\[
dX_{t}=-\al_{0}X_{t}\,dt+\beta_{1}\,dW_{t}+
\beta_{2}\,dJ_{t},
\]
where $\al_{0},\beta_{1},\beta_{2}>0$.
In this case, by means of the special relation\break $m\al_{0}\kappa
(m)=\kappa
_{Z}(m)$ for $m\in\mbbn$,
where $\kappa(m)$ and $\kappa_{Z}(m)$, respectively, denote the $m$th
cumulants of
$\pi_{0}$ and $\mcl(\beta_{1}W_{1}+\beta_{2}J_{1})$ (cf.
Barndorff-Niesen and Shephard~\cite{BarShe01}),
we have
\[
A^{\ast}_{0}-\bigl\{-\mbbg_{\infty}^{\prime\al}(
\theta_{0})\bigr\} =\frac{\beta_{2}^{2}}{\beta_{1}^{2}(\beta_{1}^{2}+\beta
_{2}^{2})}\int x^{2}
\pi_{0}(dx) =\frac{1}{2\al_{0}} \biggl(\frac{\beta_{2}}{\beta_{1}}
\biggr)^{2},
\]
which becomes larger (resp., smaller) with increasing (resp., decreasing)
$\beta_{2}^{2}/\beta_{1}^{2}$, the ratio of the jump-part variance to
the diffusion-part one.

Furthermore, if $X$ is supposed to be of pure-jump driven type (i.e.,
$b\equiv0$) from the very beginning,
the optimal rate of convergence in estimating $\al$ may be faster than
$\sqrt{T_{n}}$.
For example, if $X$ is the Ornstein--Uhlenbeck-type process
$dX_{t}=-\al X_{t}\,dt+dJ_{t}$
and if $\mcl(h^{-1/\gamma}J_{h})$ for small $h$ behaves like
the non-Gaussian $\gamma$-stable distribution [$\gamma\in(0,2)$],
then the least absolute deviation (LAD)-type estimator has asymptotic
normality at rate $\sqrt{n}h_{n}^{1-1/\gamma}$,
which is faster than $\sqrt{T_{n}}=\sqrt{nh_{n}}$; see Masuda \cite
{Mas10} for details.
Unfortunately, it is not clear whether or not it is possible to
generalize the LAD-type estimation method to deal with $X$ of (\ref
{SDE}) with nonlinear coefficients.

\item Let us consider
%
%
%e2.17 #&#
\begin{equation}
dX_{t}=a(X_{t},\al)\,dt+c(X_{t-},
\beta)\,dJ_{t}, \label{pj_SDE}
\end{equation}
where $J$ is a centered pure-jump L\'evy process of infinite activity
[i.e., $\nu(\mbbr)=\infty$]
such that $E[J_{1}^{2}]=1$.
Sometimes, a pure-jump L\'evy process $J$ can be approximated by a
standard Wiener process
if the parameter contained in the L\'evy measure $\nu(dz)$ behaves suitably;
for instance, $\mcl(J_{1})\to\mcn_{1}(0,1)$ as $\del\to\infty$ if
$\mcl
(J_{1})$ obeys
the symmetric centered normal inverse-Gaussian distribution $\operatorname{NIG}(\del,0,\del,0)$.
Although the rate of convergence $\sqrt{T_{n}}$ of our GQMLE $\hat
{\beta
}_{n}$ can be never improved
as long as we have a nonnull jump part, it is expected, in general,
that if $\mcl(J_{1})$ in (\ref{pj_SDE}) gets ``closer'' to the normal
distribution
[i.e., if both $|\nu(3)|$ and $\nu(4)$ become small], our GQMLE will
exhibit better performance;
see Table~\ref{sim_table_1} in Section~\ref{sec_sim} for some
simulation results in this setting.
As a matter of fact, Theorem~\ref{thmm_main} verifies that
\[
\sup_{n\in\mbbn}V_{0} \bigl[\sqrt{T_{n}}(
\hat{\beta}_{n}-\beta _{0}) \bigr]\lesssim\nu(4).
\]
[Recall that $\mbbv_{\beta\beta}$ depends on $\nu(4)$ linearly.]
It is worth mentioning that, even though $\hat{\beta}_{n}$ is here
$\sqrt{T_{n}}$-consistent,
$\sqrt{n}(\hat{\beta}_{n}-\beta_{0})$ behaves like a tight sequence
if $\kappa_{n}:=\nu(4)$ gets smaller as $\kappa_{n}=O(h_{n})$.
\end{itemize}

%f1 #&#
\begin{figure}[b]

\includegraphics{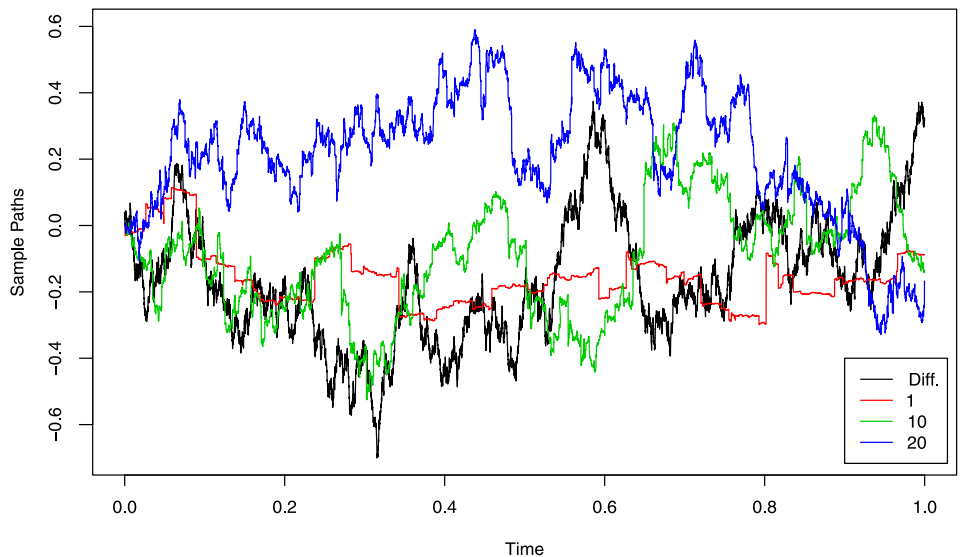}

\caption{Plots of sample paths of $X$ of (\protect\ref{sim_model}) for
$\del=1,10$ and $20$,
with a diffusion corresponding to $X$ with $J$ replaced by a standard
Wiener process.}
\label{GQMLE_fig1}
\end{figure}

%%%

%s2.4 #&#
\subsection{A numerical example}\label{sec_sim}

For simulation purposes, we consider the following concrete model:
%
%
%e2.18 #&#
\begin{equation}
dX_{t}=\frac{-\al X_{t}}{\sqrt{1+X_{t}^{2}}}\,dt+\sqrt{\beta}\,dJ_{t},\qquad
X_{0}=0, \label{sim_model}
\end{equation}
where the true value is $(\al_{0},\beta_{0})=(1,1)$, the driving
process is
the normal inverse Gaussian L\'evy process such that $\mcl
(J_{t})=\operatorname{NIG}(\del,0,\del t,0)$, where $\del=1,10$ or $20$.
It holds that $E[J_{t}]=0$, $E[J_{1}^{2}]=t$, and $\mcl(J_{t})\to\mcn
(0,t)$ in total variation as $\del\to\infty$,
and that $\nu(3)=0$ and $\nu(4)=3/\del^{2}$.
Model (\ref{sim_model}) is a normal-inverse Gaussian counterpart to the
hyperbolic diffusion,
for which $J$ is replaced by a standard Wiener process.
For this $X$, we can verify all the assumptions; see Proposition~\ref
{prop_ergo} for the verification of
the stability conditions.

We simulated $1000$ independent paths by Euler scheme with sufficiently
fine step size to obtain
$1000$ independent estimates $\hat{\theta}_{n}=(\hat{\al}_{n},\hat
{\al}_{n})$,
and then computed their empirical mean and standard deviations.

Figure~\ref{GQMLE_fig1} shows typical sample paths of $X$ for $\del
=1,10$, and $20$,
with a diffusion corresponding to $X$ with $J$ replaced with a standard
Wiener process, just for comparison.

Table~\ref{sim_table_1} reports the results;
just for comparison, we included the case of diffusion, where $J$ is a
standard Wiener process.
From the table, we can observe the following:
\begin{itemize}
\item the performance of $\hat{\al}_{n}$ are rather similar for all the
three cases;
\item the performance of $\hat{\beta}_{n}$ gets better for larger
$\del
$, which can be expected from
the fact that the asymptotic variance of $\hat{\beta}_{n}$ is a
constant multiple of $\nu(4)=3\del^{-2}$;
we have $\mbbv_{\beta\beta}\to0$ as $\del\to\infty$.
\end{itemize}

%t1 #&#
\begin{table}
\caption{Finite sample performance of $\hat{\theta}_{n}$ concerning the
model (\protect\ref{sim_model});
just for comparison, the case of diffusion is also included.
In each case, the sample mean is given with
the sample standard deviation in parenthesis}\label{sim_table_1}
\begin{tabular*}{\textwidth}{@{\extracolsep{\fill}}lccccccccc@{}}
\hline
 &  & \multicolumn{2}{c}{\textbf{Diffusion}} & \multicolumn
{2}{c}{$\bolds{\del=1}$} & \multicolumn{2}{c}{$\bolds{\del=10}$}
& \multicolumn{2}{c@{}}{$\bolds{\del=20}$} \\[-6pt]
 &  & \multicolumn{2}{c}{\hrulefill} & \multicolumn
{2}{c}{\hrulefill} & \multicolumn{2}{c}{\hrulefill}
& \multicolumn{2}{c@{}}{\hrulefill} \\
\multicolumn{1}{@{}l}{$\bolds{T_n}$}& \multicolumn{1}{c}{$\bolds{h_n}$}& \multicolumn{1}{c}{$\bolds{\al}$} &
\multicolumn{1}{c}{$\bolds{\beta}$} & \multicolumn{1}{c}{$\bolds{\al}$} & \multicolumn{1}{c}{$\bolds{\beta}$} & \multicolumn{1}{c}{$\bolds{\al}$} &
\multicolumn{1}{c}{$\bolds{\beta}$}
& \multicolumn{1}{c}{$\bolds{\al}$} & \multicolumn{1}{c@{}}{$\bolds{\beta}$} \\
\hline
\phantom{0}10 & 0.05 & 1.16 & 0.96 & 1.15 & 0.98 & 1.18 & 0.97 & 1.18 & 0.96 \\
& & (0.63) & (0.10) & (0.62) & (0.58) & (0.65) & (0.11) & (0.65) &
(0.10) \\
& 0.01 & 1.19 & 0.99 & 1.17 & 0.97 & 1.21 & 0.99 & 1.19 & 0.99 \\
& & (0.67) & (0.04) & (0.64) & (0.48) & (0.66) & (0.07) & (0.68) &
(0.05) \\[3pt]
100& 0.05  & 1.00 & 0.97 & 1.00 & 0.98 & 1.00 & 0.97 & 1.01 & 0.97 \\
& & (0.18) & (0.03) & (0.19) & (0.17) & (0.18) & (0.04) & (0.17) &
(0.03) \\
& 0.01 &  1.02 & 0.99 & 1.02 & 1.00 & 1.02 & 0.99 & 1.03 & 1.00 \\
& & (0.18) & (0.01) & (0.19) & (0.17) & (0.18) & (0.02) & (0.19) &
(0.02) \\
\hline
\end{tabular*}
%%
%%
%$T_n$ & $h_n$ && \multicolumn{2}{c}{Diffusion} & \multicolumn
%{2}{c}{$\del=1$} & \multicolumn{2}{c}{$\del=10$}
%& \multicolumn{2}{c}{$\del=20$} \\
%& && $\al$ & $\beta$ & $\al$ & $\beta$ & $\al$ & $\beta$
%& $\al$ & $\beta$ \\
%10 & 0.05 & & 1.158 & 0.964 & 1.149 & 0.982 & 1.179 & 0.965 & 1.181 &
%0.964 \\
%& & & (0.629) & (0.100) & (0.621) & (0.576) & (0.645) & (0.113) &
%(0.649) & (0.103) \\[1mm]
%& 0.01 & & 1.193 & 0.993 & 1.172 & 0.968 & 1.210 & 0.993 & 1.190 &
%0.994 \\
%& & & (0.673) & (0.044) & (0.635) & (0.476) & (0.658) & (0.069) &
%(0.682) & (0.054) \\
%100& 0.05 & & 0.999 & 0.970 & 0.997 & 0.976 & 0.996 & 0.971 & 1.006 &
%0.974 \\
%& & & (0.184) & (0.031) & (0.188) & (0.174) & (0.177) & (0.035) &
%(0.171) & (0.032) \\[1mm]
%& 0.01 & & 1.018 & 0.994 & 1.017 & 0.996 & 1.017 & 0.994 & 1.029 &
%0.995 \\
%& & & (0.184) & (0.014) & (0.189) & (0.174) & (0.181) & (0.022) &
%(0.191) & (0.017) \\[1mm]
%%
\end{table}

%%%%%
%%%%%

%s3 #&#
\section{Mighty convergence of a class of continuous random
fields}\label{sec_MC}

In this section, we prove a fundamental result concerning
the ``single-norming'' mighty convergence of a continuous statistical
random fields
associated with general vector-valued estimating functions; here, the
``single-norming'' means that
the rates of convergence are the same for all the arguments of the
corresponding estimator.
Theorem~\ref{thmm_PLDI_G} below will serve as a fundamental tool in the
proof of Theorem~\ref{thmm_main};
the content of this section can be read independently of the main body.

To proceed, we need some notation. Denote by $\{\mcx_{n},\mca
_{n},(P_{\theta})_{\theta\in\Theta}\}_{n\in\mbbn}$
underlying statistical experiments, where $\Theta\subset\mbbr^{p}$
is a
bounded convex domain.
Let $\theta_{0}\in\Theta$, and write $P_{0}=P_{\theta_{0}}$.
Let $\mbbg_{n}=(\mbbg_{j,n})_{j=1}^{p}\dvtx \mcx_{n}\times\Theta\to
\mbbr
^{p}$ be vector-valued random functions;
as usual, we will simply write $\mbbg_{n}(\theta)$, dropping the
argument of~$\mcx_{n}$.
Our target ``contrast'' function is
%
%
%e3.1 #&#
\begin{equation}
\mbbm_{n}(\theta):=-\frac{1}{T_{n}}\bigl|\mbbg_{n}(
\theta)\bigr|^{2}, \label{MC_H_def}
\end{equation}
where $(T_{n})$ is a nonrandom positive real sequence such that
$T_{n}\to\infty$.
The corresponding ``$M$-estimator'' is defined to be any measurable mapping
$\hat{\theta}_{n}\dvtx \mcx_{n}\to\overline{\Theta}$ such that
\[
\hat{\theta}_{n}\in\argmax_{\theta\in\overline{\Theta}}\mbbm _{n}(
\theta).
\]
Due to the compactness of $\overline{\Theta}$ and the continuity of
$\mbbm_{n}$ imposed later on,
we can always find such a $\hat{\theta}_{n}$.
The estimate $\hat{\theta}_{n}$ can be any root of $\mbbg_{n}(\theta
)=0$ as soon as it exists.

Set $U_{n}(\theta_{0}):=\{u\in\mbbr^{p}\dvtx  \theta
_{0}+T_{n}^{-1/2}u\in
\Theta\}$
and define random fields $\mbbz_{n}\dvtx\break U_{n}(\theta_{0})\to(0,\infty
)$ by
%
%
%e3.2 #&#
\begin{equation}
\mbbz_{n}(u)=\mbbz_{n}(u;\theta_{0}):=\exp
\bigl\{\mbbm_{n}\bigl(\theta _{0}+T_{n}^{-1/2}u
\bigr)-\mbbm_{n}(\theta_{0}) \bigr\}. \label{MC_Z_def}
\end{equation}
Obviously, it holds that
\[
\hat{u}_{n}:=\sqrt{T_{n}}(\hat{\theta}_{n}-
\theta_{0})\in\argmax _{\theta
\in\overline{\Theta}}\mbbz_{n}(\theta).
\]
We consider the following two conditions for the random fields $\mbbz_{n}$.
\begin{itemize}
\item{[\textit{Polynomial type Large Deviation Inequality (PLDI)}]}.
For every $M>0$, we have
%
%
%e3.3 #&#
\begin{equation}
\sup_{r>0} \Bigl\{r^{M}\sup_{n\in\mbbn}P_{0}
\Bigl[\sup_{|u|>r}\mbbz _{n}(u)\ge e^{-r}
\Bigr] \Bigr\}<\infty. \label{PLDI_def}
\end{equation}
\item(\textit{Weak convergence on compact sets}).
There exists a random field $\mbbz_{0}(\cdot)=\mbbz_{0}(\cdot;\theta
_{0})$ such that
$\mbbz_{n}\to^{\mcl}\mbbz_{0}$ in $\mcc(\overline{B(R)})$ for
each $R>0$,
where $\overline{B(R)}:=\{u\in\mbbr^{p}; |u|\le R\}$.
\end{itemize}
Under these conditions, the mode of convergence of $\mbbz_{n}(\cdot)$
is mighty enough to deduce that
the maximum-point sequence $(\hat{u}_{n})_{n}$ is
$L^{q}(P_{0})$-bounded for every \mbox{$q>0$},
which especially implies that $(\hat{u}_{n})_{n}$ is tight:
indeed, if (\ref{PLDI_def}) is in force,
\begin{eqnarray*}
\sup_{n\in\mbbn}P_{0}\bigl[|\hat{u}_{n}|>r\bigr] &\le&
\sup_{n\in\mbbn}P_{0} \Bigl[\sup_{|u|>r}
\mbbz_{n}(u)\ge \mathbb {Z}_{n}(0) \Bigr]
\\
&=&\sup_{n\in\mbbn}P_{0} \Bigl[\sup_{|u|>r}
\mbbz_{n}(u)\ge1 \Bigr]\lesssim\frac{1}{r^{M}}
\end{eqnarray*}
for every $r>0$, so that
\[
\sup_{n\in\mbbn}E_{0}\bigl[|\hat{u}_{n}|^{q}
\bigr]= \int_{0}^{\infty}\sup_{n\in\mbbn}P_{0}
\bigl[|\hat{u}_{n}|>s^{1/q}\bigr]\,ds \lesssim1+\int
_{1}^{\infty}s^{-M/q}\,ds<\infty.
\]
If $u\mapsto\mbbz_{0}(u)$ is a.s. maximized at a unique point $\hat
{u}_{\infty}$,
then it follows from the tightness of $(\hat{u}_{n})_{n\in\mbbn}$ that
$\hat{u}_{n}\to^{\mcl}\hat{u}_{\infty}$;
let us remind the reader that the weak convergence on any compact set alone
is not enough to deduce the weak convergence of $\hat{u}_{n}$,
since $U_{n}(\theta_{0})\uparrow\mbbr^{p}$ and we have no guarantee
that $(\hat{u}_{n})$ is tight.
Moreover, owing to the PLDI, the moment of $f(\hat{u}_{n})$ converges
to that of\vadjust{\goodbreak} $f(\hat{u}_{\infty})$
for every continuous function $f$ on $\mbbr^{p}$ of at most polynomial growth.
In our framework, $\log\mbbz_{0}$ admits a quadratic structure with
a normally distributed linear term and a nonrandom positive definite
quadratic term,
so that $\hat{u}_{\infty}$ is asymptotically normally distributed.

We now introduce regularity conditions.

%as3.1 #&#
\begin{ass}[(Smoothness)]\label{MC_A1}
The functions $\theta\mapsto\mbbg_{n}(\theta)$ are continuously
extended to the boundary of $\Theta$,
and belong to $\mcc^{3}(\Theta)$, $P_{0}$-a.s.
\end{ass}

%as3.2 #&#
\begin{ass}[(Bounded moments)]\label{MC_A2}
For every $K>0$,
\[
\sup_{n\in\mbbn}E_{0} \biggl[\biggl\llvert
\frac{1}{\sqrt{T_{n}}}\mbbg _{n}(\theta _{0})\biggr\rrvert
^{K} \biggr] +\max_{k\in\{0,1,2,3\}} \sup_{n\in\mbbn}E_{0}
\biggl[\sup_{\theta\in\Theta} \biggl\llvert \frac{1}{T_{n}}
\p_{\theta}^{k}\mbbg_{n}(\theta)\biggr\rrvert
^{K} \biggr] <\infty.
\]
\end{ass}

Let $M>0$ be a given constant.

%as3.3 #&#
\begin{ass}[(Limits)]\label{MC_A3}
(a) There exist a nonrandom function $\mbbg_{\infty}\dvtx\break \Theta
\to
\mbbr^{p}$ and
positive constants $\chi=\chi(\theta_{0})$ and $\ep$ such that:
$\mbbg_{\infty}(\theta_{0})=0$; $\sup_{\theta}|\mbbg_{\infty
}(\theta
)|<\infty$;
$|\mbbg_{\infty}(\theta)|^{2}\ge\chi|\theta-\theta_{0}|^{2}$ for every
$\theta\in\Theta$; and
\[
\sup_{n\in\mbbn}E_{0} \biggl[\sup_{\theta\in\Theta}
\biggl\llvert \sqrt{T_{n}} \biggl(\frac{1}{T_{n}}
\mbbg_{n}(\theta)-\mbbg_{\infty}(\theta ) \biggr)\biggr\rrvert
^{M+\ep} \biggr]<\infty.
\]
\begin{longlist}[(b)]
\item[(b)] There exists a nonrandom
$\mbbg_{\infty}'(\theta_{0})\in\mbbr^{p}\otimes\mbbr^{p}$ of
rank $p$
such that
\[
\sup_{n\in\mbbn}E_{0} \biggl[ \biggl\llvert
\sqrt{T_{n}} \biggl(\frac{1}{T_{n}}\p_{\theta}\mbbg
_{n}(\theta_{0}) -\mbbg_{\infty}'(
\theta_{0}) \biggr)\biggr\rrvert ^{M} \biggr]<\infty.
\]
\end{longlist}
\end{ass}

%as3.4 #&#
\begin{ass}[(Weak convergence)]\label{MC_A4}
$T_{n}^{-1/2}\mbbg_{n}(\theta_{0})\to^{\mcl}\mcn_{p}(0,\mbbv
(\theta_{0}))$
for some positive definite $\mbbv(\theta_{0})\in\mbbr^{p}\otimes
\mbbr^{p}$.
\end{ass}

Let $\Sig(\theta_{0}):=(\mbbg'_{\infty})^{-1}\mbbv(\mbbg'_{\infty})^{-1
\top}(\theta_{0})$.
The main claim of this section is the following.

%th3.5 #&#
\begin{thmm}\label{thmm_PLDI_G}
Let $M>0$.
\begin{longlist}[(a)]
\item[(a)] Suppose that Assumptions~\ref{MC_A1},~\ref{MC_A2} and
\ref
{MC_A3} hold.
Then the PLDI (\ref{PLDI_def}) holds true.
\item[(b)] If Assumption~\ref{MC_A4} is additionally met, then
\[
E_{0}\bigl[f(\hat{u}_{n})\bigr]\to\int f(u)\phi
\bigl(u;0,\Sig(\theta _{0}) \bigr)\,du
\]
for every continuous function $f\dvtx \mbbr^{p}\to\mbbr$ satisfying that
$\limsup_{|u|\to\infty}|u|^{-q}\times |f(u)|<\infty$ for some $q\in(0,M)$.
\end{longlist}
\end{thmm}

\begin{pf}
Applying Taylor's expansion to (\ref{MC_Z_def}), we get
%
%
%e3.4 #&#
\begin{equation}
\log\mbbz_{n}(u)=\D_{n}(\theta_{0})[u]-
\tfrac{1}{2}\Gam(\theta _{0})[u,u]+\xi_{n}(u),
\label{logZn_expa}\vadjust{\goodbreak}
\end{equation}
where $\D_{n}(\theta_{0}):=T_{n}^{-1/2}\,\p_{\theta}\mbbm_{n}(\theta_{0})$,
$\Gam_{n}(\theta_{0}):=-T_{n}^{-1}\,\p_{\theta}^{2}\mbbm_{n}(\theta_{0})$,
$\Gam(\theta_{0}):=\break 2\mbbg'_{\infty}(\theta_{0})^{\top}\mbbg
'_{\infty
}(\theta_{0})$ and
%
%
%e3.5 #&#
\begin{eqnarray}\label{xi_form}
\xi_{n}(u)&:=&\frac{1}{2}\bigl\{\Gam(\theta_{0})-
\Gam_{n}(\theta_{0})\bigr\}[u,u]
\nonumber
\\[-8pt]
\\[-8pt]
\nonumber
&&{} -\int_{0}^{1}(1-s)
\int\p_{\theta}\Gam_{n}\bigl(\theta_{0}+stT_{n}^{-1/2}u
\bigr) \bigl[sT_{n}^{-1/2}u,u^{\otimes2}\bigr]\,dt\,ds.
\end{eqnarray}

We will prove (a) by making use of Yoshida~\cite{Yos11}, Theorem 3(c).
The task is then to verify conditions [A1$''$], [A4$'$], [A6], [B1] and
[B2] of that paper.
For convenience and clarity, we will list them in a reduced form with
our notation.
First we look at [B1] and [B2]:
\begin{longlist}[{[B1]}]%{\leftmargin=3pt,\itemindent=0pt}
\item[{[B1]}] the matrix $\Gam(\theta_{0})$ is positive definite;
\item[{[B2]}] there exists a constant $\chi>0$ such that
$\mbby(\theta)\le-\chi^{2}|\theta-\theta_{0}|^{2}$ for each
$\theta\in
\Theta$.
\end{longlist}
Here $\mbby(\theta):=-|\mbbg_{\infty}(\theta)|^{2}$,
where $\mbbg_{\infty}(\theta)$ is the one appearing in Assumption~\ref{MC_A3}.
Obviously, Assumption~\ref{MC_A3} assures [B1] and [B2] (the identifiability);
in particular, we have the convergence
$T_{n}^{-1}\mbbm_{n}(\theta)\to^{p}-|\mbbg_{\infty}(\theta)|^{2}$ for
each $\theta\in\Theta$, so that
\[
\mbby_{n}(\theta):=\frac{1}{T_{n}} \bigl\{ \mbbm_{n}(
\theta)-\mbbm_{n}(\theta_{0}) \bigr\} =\frac{1}{T_{n}}
\log\mbbz_{n} \bigl(\sqrt{T_{n}}(\theta-\theta
_{0}) \bigr) \to^{p}\mbby(\theta).
\]
Next, given constants $M>0$ [the number in (\ref{PLDI_def})] and $\al
\in(0,1)$,
conditions [A6], [A1$''$] and [A4$'$] read as follows:
\begin{enumerate}[{[A4$'$]}]%{\leftmargin=3pt,\itemindent=0pt}
\item[{[A6]}]
\hspace*{3pt}(i)\hspace*{3pt} $\sup_{n}E_{0} [\llvert \D_{n}(\theta_{0})\rrvert ^{M_{1}} ]<\infty$ for $M_{1}:=M/(1-\rho_{1})$.
\begin{enumerate}[(ii)]
\item[(ii)] $\sup_{n}E_{0} [\sup_{\theta}
\llvert T_{n}^{1/2-\beta_{2}}(\mbby_{n}(\theta)-\mbby(\theta
))\rrvert ^{M_{2}} ]<\infty$,
for $M_{2}:=M/(1-2\beta_{2}-\rho_{2})$.
\end{enumerate}
\item[{[A1$''$]}]
\begin{enumerate}
\item[(i)] $\sup_{n}E_{0} [\sup_{\theta}
\llvert T_{n}^{-1}\p_{\theta}^{3}\mbbm_{n}(\theta)\rrvert ^{M_{3}}
]<\infty$
for $M_{3}:=M/\{\al/(1-\al)-\rho_{1}\}$.
\item[(ii)]
$\sup_{n}E_{0} [\llvert T_{n}^{\beta_{1}}(\Gam_{n}(\theta
_{0})-\Gam
(\theta_{0})\rrvert ^{M_{4}} ]<\infty$
for $M_{4}:=M/\{2\beta_{1}/(1-\al)-\rho_{1}\}$.
\end{enumerate}
\item[{[A4$'$]}] The parameters $\al$, $\beta_{1}$, $\beta_{2}$,
$\rho_{1}$ and $\rho_{2}$ fulfil the inequalities
\begin{eqnarray*}
&\displaystyle 0<\beta_{1}<1/2, \qquad 0<\rho_{1}<\min \biggl(1,
\frac
{\al
}{1-\al}, \frac{2\beta_{1}}{1-\al} \biggr),&
\\
&\displaystyle  2\al<\rho_{2},\qquad \beta_{2}\ge0,\qquad 1-2\beta_{2}-
\rho_{2}>0.&
\end{eqnarray*}
\end{enumerate}
These conditions involve several ``moment-index'' parameters to be controlled,
which do not seem straightforward to handle.
Nevertheless, under our assumptions we can provide a rather simplified version.
Instead of ``[A1$''$], [A4$'$] and [A6]'' we will verify the following
``[A1$''{}^\sharp$] and [A6${}^\sharp$]'':
\begin{enumerate}[{[A1$''{}^\sharp$]}]%{\leftmargin=3pt,\itemindent=0pt}
\item[{[A1$''{}^\sharp$]}]
\hspace*{3pt}(i)\hspace*{3pt} $\sup_{n}E_{0} [\sup_{\theta}
\llvert T_{n}^{-1}\p_{\theta}^{3}\mbbm_{n}(\theta)\rrvert ^{K}
]<\infty$ for every $K>0$.\vspace*{2pt}
\begin{enumerate}[(ii)]
\item[(ii)] $\sup_{n}E_{0} [
\llvert \sqrt{T_{n}}(\Gam_{n}(\theta_{0})-\Gam(\theta_{0})\rrvert ^{M-\ep
_{1}} ]<\infty$
for every $\ep_{1}>0$ small enough.
\end{enumerate}
\item[{[A6${}^\sharp$]}]
\hspace*{3pt}(i)\hspace*{3pt} $\sup_{n}E_{0} [\llvert \D_{n}(\theta_{0})\rrvert ^{K}
]<\infty$ for every $K>0$.\vspace*{2pt}
\begin{enumerate}
\item[(ii)]
$\sup_{n}E_{0} [\sup_{\theta}
\llvert \sqrt{T_{n}}(\mbby_{n}(\theta)-\mbby(\theta))\rrvert ^{M+\ep
/2} ]<\infty$,
for $\ep$ given in Assumption~\ref{MC_A3}.
\end{enumerate}
\end{enumerate}
Let us show that ``[A1$''{}^\sharp$] and [A6${}^\sharp$]'' imply
``[A1$''$], [A4$'$] and [A6].''
First, by [A1$''{}^\sharp$](i) and [A6${}^\sharp$](i), the numbers
$M_{1}$ and $M_{3}$ can be arbitrarily large,
so that we may in particular take $\al$ and $\rho_{1}$ arbitrarily
small (i.e., nearly zero).
Then we have [A1$''$](i) and [A6](i).
Next, we note that in [A1$''{}^\sharp$](ii) the exponent of ``$T_{n}$''
is $1/2$, hence
we may let $\beta_{2}$ be sufficiently close to $1/2$.
Then, taking $\al$ and $\rho_{1}$ small enough with $\rho_{1}<\al
/(1-\al)$,
we can obtain the first two inequalities in [A4$'$].
Next, in view of [A6${}^\sharp$](ii), we can take $\beta_{2}=0$
and $\rho_{2}$ small enough to make [A6](ii) and the last three ones in
[A4$'$] valid.
Finally, as for $M_{4}$, we note that a suitable control of $(\al,\rho
_{1},\beta_{1})$ leads to
\[
\frac{2\beta_{1}}{1-\al}-\rho_{1} =1+ \biggl(\frac{\al}{1-\al}-
\rho_{1} \biggr)+\frac{2\beta
_{1}-1}{1-\al}>1,
\]
so that [A1$''$](ii) follows. In sum, under ``[A1$''{}^\sharp$] and
[A6${}^\sharp$],''
we can pick $\rho_{1},\rho_{2},\al\approx0$ and $\beta_{2}=0$, and
then $\beta_{1}\approx1/2$,
in order to make all of ``[A1$''$], [A4$'$] and [A6]'' valid.
Thus we are left to proving [A1$''{}^\sharp$] and [A6${}^\sharp$]
above.\vspace*{2pt}

We begin with [A1$''{}^\sharp$]. Since\vspace*{2pt} $|T_{n}^{-1}\p_{\theta
}^{3}\mbbm
_{n}(\theta)|\lesssim
|T_{n}^{-1}\mbbg_{n}(\theta)\Vert T_{n}^{-1}\p_{\theta}^{3}\mbbg
_{n}(\theta)|
+|T_{n}^{-1}\p_{\theta}\mbbg_{n}(\theta)\Vert T_{n}^{-1}\p_{\theta
}^{2}\mbbg
_{n}(\theta)|$,
we have for every $K>0$,
\[
\sup_{n\in\mbbn}E_{0} \biggl[\sup_{\theta\in\Theta}
\biggl\llvert \frac{1}{T_{n}}\p_{\theta}^{3}
\mbbm_{n}(\theta)\biggr\rrvert ^{K} \biggr]<\infty.
\]
Noting that $\p_{\theta_{i}}\,\p_{\theta_{j}}\mbbm_{n}
=-2T_{n}^{-1}\{\p_{\theta_{i}}\,\p_{\theta_{j}}\mbbg_{n}[\mbbg
_{n}]+\p
_{\theta_{i}}\mbbg_{n}
[\p_{\theta_{j}}\mbbg_{n}]\}$, we also have
\begin{eqnarray*}
&& \sqrt{T_{n}}\bigl|\Gam_{n}(\theta_{0})-\Gam(
\theta_{0})\bigr|
\\
&&\qquad\lesssim \biggl\llvert \frac{1}{\sqrt{T_{n}}}\mbbg_{n}(
\theta_{0})\biggr\rrvert \biggl\llvert \frac{1}{T_{n}}
\p_{\theta}^{2}\mbbg_{n}(\theta_{0})
\biggr\rrvert
\\
&&\qquad\quad{}{} + \biggl(\bigl|\Gam(\theta_{0})\bigr|+\biggl\llvert \frac{1}{T_{n}}
\p_{\theta
}\mbbg _{n}(\theta_{0})\biggr\rrvert
\biggr) \biggl\llvert \sqrt{T_{n}} \biggl( \frac{1}{T_{n}}
\p_{\theta}\mbbg_{n}(\theta_{0})-
\mbbg'_{\infty
}(\theta_{0}) \biggr)\biggr\rrvert.
\end{eqnarray*}
Therefore, Assumptions~\ref{MC_A2} and~\ref{MC_A3} combined with H\"
older's inequality yield that
for $\ep_{1}\in(0,M)$,
\begin{eqnarray*}
& &\sup_{n\in\mbbn}E_{0} \bigl[\bigl\llvert
\sqrt{T_{n}}(\Gam_{n}(\theta_{0})-\Gam(
\theta_{0})\bigr\rrvert ^{M-\ep
_{1}} \bigr]
\\
&&\qquad\lesssim 1+ \biggl\{\sup_{n\in\mbbn}E_{0} \biggl[
\biggl\llvert \sqrt{T_{n}} \biggl( \frac{1}{T_{n}}\p_{\theta}
\mbbg_{n}(\theta_{0})-\mbbg'_{\infty
}(
\theta_{0}) \biggr)\biggr\rrvert ^{M} \biggr] \biggr
\}^{(M-\ep_{1})/M}<\infty.
\end{eqnarray*}
Thus [A1$''{}^\sharp$] follows.\vadjust{\goodbreak}

Next we prove [A6${}^\sharp$]. Statement (i) is obvious from Assumption
\ref{MC_A2},
%
%
%e3.6 #&#
\begin{equation}\quad
\sup_{n\in\mbbn}E_{0} \bigl[\bigl|\D_{n}(
\theta_{0})\bigr|^{K} \bigr] \lesssim \sup_{n\in\mbbn}E_{0}
\biggl[ \biggl\llvert \frac{1}{T_{n}}\p_{\theta}\mbbg_{n}(
\theta_{0})\biggr\rrvert ^{K} \biggl\llvert
\frac{1}{\sqrt{T_{n}}}\mbbg_{n}(\theta_{0})\biggr\rrvert
^{K} \biggr]<\infty. \label{thmm_PLDI_G_p1}
\end{equation}
Using the estimate
\begin{eqnarray*}
&& \bigl\llvert \sqrt{T_{n}}\bigl(\mbby_{n}(\theta)-\mbby(
\theta)\bigr)\bigr\rrvert
\\
&&\qquad\le\frac{1}{\sqrt{T_{n}}}\biggl\llvert \frac{1}{\sqrt{T_{n}}}\mbbg _{n}(
\theta _{0})\biggr\rrvert ^{2} \\
&&\qquad\quad{}+ \biggl(\bigl|
\mbbg_{\infty}(\theta)\bigr|+\biggl\llvert \frac{1}{T_{n}}\mbbg
_{n}(\theta )\biggr\rrvert \biggr) \biggl\llvert \sqrt{T_{n}}
\biggl(\frac{1}{T_{n}}\mbbg_{n}(\theta)-\mbbg _{\infty
}(
\theta) \biggr)\biggr\rrvert,
\end{eqnarray*}
it follows under Assumptions~\ref{MC_A2} and~\ref{MC_A3} that
\begin{eqnarray*}
&& \sup_{n\in\mbbn}E_{0} \Bigl[\sup_{\theta\in\Theta}
\bigl\llvert \sqrt{T_{n}}\bigl(\mbby_{n}(\theta)-\mbby(
\theta)\bigr)\bigr\rrvert ^{M+\ep
/2} \Bigr]
\\
&&\quad\lesssim1+\sup_{n\in\mbbn}E_{0} \biggl[\sup
_{\theta\in\Theta} \biggl\llvert \sqrt{T_{n}} \biggl(
\frac{1}{T_{n}}\mbbg_{n}(\theta)-\mbbg _{\infty
}(\theta)
\biggr)\biggr\rrvert ^{M+\ep} \biggr] ^{(M+\ep/2)/(M+\ep)}<\infty.
\end{eqnarray*}
Thus [A6${}^\sharp$] is ensured, and the proof of (a) is complete.

We now turn to the proof of (b). Fix any $R>0$.
Since we know that the sequence $(\hat{u}_{n})$ is
$L^{q}(P_{0})$-bounded for each $q\in(0,M)$
and that the set $\argmax_{u}\log\mbbz_{\infty}(u)$ a.s. consists of
the only point
\[
\hat{u}_{\infty}:=\Gam(\theta_{0})^{-1}
\D_{\infty}(\theta _{0})\sim\mcn _{p}\bigl(0,\Sig(
\theta_{0})\bigr),
\]
it suffices to show that $\log\mbbz_{n}\to^{\mcl}\log\mbbz
_{\infty}$
in $\mcc(\overline{B(R)})$, where
\begin{eqnarray*}
\log\mbbz_{\infty}(u)&:=&\D_{\infty}(\theta_{0})[u]-
\tfrac
{1}{2}\Gam (\theta_{0})[u,u],
\\
\D_{\infty}(\theta_{0}) &\sim& \mcn_{p} \bigl(0,4
\mbbg'_{\infty}(\theta_{0})^{\top} \mbbv(
\theta_{0})\mbbg'_{\infty}(\theta_{0})
\bigr)
\end{eqnarray*}
(e.g., Yoshida~\cite{Yos11}, Theorem 5).
We have $T_{n}^{-1}\p_{\theta}\mbbg_{n}(\theta_{0})\to^{p}\mbbg
_{\infty
}'(\theta_{0})$ from Assumption~\ref{MC_A3},
hence Slutsky's lemma and Assumption~\ref{MC_A4} imply that
\[
\D_{n}(\theta_{0})=-\frac{2}{T_{n}}\p_{\theta}
\mbbg_{n}(\theta_{0}) \biggl[\frac{1}{\sqrt{T_{n}}}
\mbbg_{n}(\theta_{0}) \biggr]\to ^{\mcl}\D
_{\infty}(\theta_{0}).
\]
Also, we have
%
%
%e3.7 #&#
\begin{equation}
\bigl|\xi_{n}(u)\bigr|\lesssim|u|^{2}\bigl|\Gam_{n}(
\theta_{0})-\Gam(\theta_{0})\bigr| +\frac{|u|^{3}}{\sqrt{T_{n}}}\sup
_{\theta\in\Theta} \biggl\llvert \frac{1}{T_{n}}\p_{\theta}^{3}
\mbbm_{n}(\theta)\biggr\rrvert =o_{p}(1)
\label{thmm_PLDI_G_p2}
\end{equation}
for every $u\in\overline{B(R)}$. Thus, recalling expression (\ref
{logZn_expa}),
we get $\log\mbbz_{n}(u)\to^{\mcl}\log\mbbz_{0}(u)$ for every
$u\in
\overline{B(R)}$, and moreover,
due to the linearity in $u$ of the weak convergence term $\D
_{n}(\theta
_{0})[u]$,
the Cram\'er--Wold device ensures the finite-dimensional convergence.
Therefore, it remains to check the tightness of $\{\log\mbbz_{n}(u)\}
_{u\in\overline{B(R)}}$.
In view of the classical Kolmogorov tightness criterion for continuous
random fields
(e.g., Kunita~\cite{Kun90}, Theorem 1.4.7),
it suffices to show that there exists a constant $\gamma>p(=\operatorname{
dim} \Theta)$ such that
%
%
%e3.8 #&#
\begin{equation}
\sup_{|u|\le R}\sup_{n\in\mbbn}E_{0}
\bigl[\bigl|\log\mbbz _{n}(u)\bigr|^{\gamma
} \bigr] +\sup
_{n\in\mbbn}E_{0} \Bigl[\sup_{|u|\le R}\bigl|
\p_{u}\log\mbbz _{n}(u)\bigr|^{\gamma} \Bigr]<\infty.
\label{tightness_ineq}
\end{equation}
In view of the estimates in (\ref{thmm_PLDI_G_p1}) and (\ref
{thmm_PLDI_G_p2}) as well as
the expressions (\ref{logZn_expa}) and (\ref{xi_form}),
\begin{eqnarray*}
& &\sup_{u\in\overline{B(R)}}\sup_{n\in\mbbn}E_{0}
\bigl[\bigl|\log \mbbz _{n}(u)\bigr|^{\gamma} \bigr]
\\
&&\qquad\lesssim\sup_{n\in\mbbn}E_{0} \bigl[\bigl|
\D_{n}(\theta_{0})\bigr|^{\gamma
} \bigr]+1 +\sup
_{u\in\overline{B(R)}}\sup_{n\in\mbbn}E_{0} \bigl[\bigl|\xi
_{n}(u)\bigr|^{\gamma
} \bigr]
\\
&&\qquad\lesssim1+E_{0} \bigl[\bigl|\Gam_{n}(\theta_{0})-
\Gam(\theta _{0})\bigr|^{\gamma
} \bigr] +\sup_{n\in\mbbn}E_{0}
\biggl[\sup_{\theta\in\Theta} \biggl\llvert \frac{1}{T_{n}}
\p_{\theta}^{3}\mbbm_{n}(\theta)\biggr\rrvert
^{\gamma
} \biggr]<\infty.
\end{eqnarray*}
Furthermore, since
\begin{eqnarray*}
\p_{u}\log\mbbz_{n}(u) &=&\p_{u} \biggl\{
\mbbm_{n} \biggl(\theta_{0}+\frac{1}{\sqrt {T_{n}}}u \biggr)-
\mbbm_{n}(\theta_{0}) \biggr\}
\\
&=&\frac{1}{\sqrt{T_{n}}}\p_{\theta}\mbbm_{n} \biggl(\theta
_{0}+\frac
{1}{\sqrt{T_{n}}}u \biggr)
\\
&=&\frac{1}{\sqrt{T_{n}}} \biggl\{\p_{\theta}\mbbm_{n}(
\theta_{0}) +\frac{1}{\sqrt{T_{n}}}\int_{0}^{1}
\p_{\theta}^{2}\mbbm_{n} \biggl(\theta_{0}+
\frac{s}{\sqrt{T_{n}}}u \biggr)[u]\,ds \biggr\},
\end{eqnarray*}
the finiteness of $\sup_{n}E_{0} [\sup_{|u|\le R}|\p_{u}\log
\mbbz
_{n}(u)|^{\gamma} ]$
follows on applying Assumption~\ref{MC_A2} to the estimate
\begin{eqnarray*}
&& \sup_{|u|\le R}\bigl\llvert \p_{u}\log
\mbbz_{n}(u)\bigr\rrvert
\\
&&\qquad\lesssim \biggl\llvert \frac{1}{\sqrt{T_{n}}}\mbbg_{n}(
\theta_{0})\biggr\rrvert \biggl\llvert \frac
{1}{T_{n}}
\p_{\theta}\mbbg_{n}(\theta_{0})\biggr\rrvert +\sup
_{\theta\in\Theta}\biggl\llvert \frac{1}{T_{n}}\p_{\theta
}^{2}
\mbbm _{n}(\theta)\biggr\rrvert
\\
&&\qquad\lesssim \biggl\llvert \frac{1}{\sqrt{T_{n}}}\mbbg_{n}(
\theta_{0})\biggr\rrvert \biggl\llvert \frac{1}{T_{n}}
\p_{\theta}\mbbg_{n}(\theta_{0})\biggr\rrvert
\\
&&\qquad\quad {} +\sup_{\theta\in\Theta} \biggl\{ \biggl\llvert \frac{1}{T_{n}}
\mbbg_{n}(\theta)\biggr\rrvert \biggl\llvert \frac
{1}{T_{n}}\p
_{\theta}^{2}\mbbg_{n}(\theta)\biggr\rrvert +\biggl
\llvert \frac{1}{T_{n}}\p_{\theta}\mbbg_{n}(\theta)\biggr
\rrvert ^{2} \biggr\}.
\end{eqnarray*}
Thus we have obtained (\ref{tightness_ineq}), thereby achieving the
proof of (b).
\end{pf}

%re3.6 #&#
\begin{rem}
We have confined ourselves to the ``single-norming (i.e.,
scalar-$T_{n}$)'' case
for the squared quasi-score function.
Nevertheless, as in the original formulation of Yoshida~\cite{Yos11}, Theorem\vadjust{\goodbreak}
1,
it would be also possible to deal with ``multi-norming'' cases
where elements of $\hat{\theta}_{n}$ possibly converge at different rates,
that is, cases of a matrix norming instead of the scalar norming $\sqrt {T_{n}}$.
This would require somewhat more complicated arguments, but we do not
need such an extension in this paper.
\end{rem}

%%%%%
%%%%%

%s4 #&#
\section{\texorpdfstring{Proofs of Theorem \protect\ref{thmm_main} and Corollary \protect\ref{cor_main}}
{Proofs of Theorem 2.7 and Corollary 2.8}}\label{sec_proofs}

%s4.1 #&#
\subsection{\texorpdfstring{Proof of Theorem \protect\ref{thmm_main}}{Proof of Theorem 2.7}}\label{sec_proof1}
The proof of Theorem~\ref{thmm_main} is achieved by applying Theorem
\ref{thmm_PLDI_G}.
When we have a reasonable estimating function $\theta\mapsto\mbbg
_{n}(\theta)$ with which
an estimator of $\theta$ is defined by a random root of the estimating
equation $\mbbg_{n}(\theta)=0$,
it may be unclear what is the ``single'' associated contrast function
to be maximized or minimized;
for example, it would be often the case when $\mbbg_{n}$ is constructed
via a kind of (conditional-) moment fittings.
The setup (\ref{M_def}) below provides a way of converting the
situation from $Z$-estimation to $M$-estimation.
%

%s4.1.1 #&#
\subsubsection{Introductory remarks}

At first glance, it seems that, in order to investigate the asymptotic
behavior of $\hat{\theta}_{n}$,
we may proceed as in the case of diffusions, expanding the GQL $\mbbq
_{n}$ of (\ref{GQL}) and then
investigating asymptotic behaviors of the derivatives $\p_{\theta
}^{k}\mbbq_{n}$;
see Yoshida~\cite{Yos11}, Section~6, for details.
Following this route, however, leads to an inconvenience, essentially due
to the fact that
$(h_{n}^{-1/2}\D_{j}X)_{j\le n}$ is not $L^{q}(P_{0})$-bounded for
$q>2$. To see this more precisely,
let us take a brief look at the simple one-dimensional L\'evy process
$X_{t}=\al t+\sqrt{\beta}J_{t}$,
with $\theta=(\al,\beta)\in\mbbr\times(0,\infty)$ and $\mcl(J_{1})$
admitting finite moments.
In this case, $\mbbq_{n}(\theta)=-\sum_{j}\{(\log\beta)+(\beta
h_{n})^{-1}(\D_{j}X-\al h_{n})^{2}\}$,
\begin{eqnarray*}
 \p_{\al}\mbbq_{n}(\theta)&=&\sum
_{j=1}^{n}\frac{2}{\beta}(\D _{j}X-
\al h_{n}),
\\
 \p_{\beta}\mbbq_{n}(\theta)&=&\sum
_{j=1}^{n}\frac{1}{\beta^{2}h_{n}} \bigl\{(
\D_{j}X-\al h_{n})^{2}-\beta h_{n}
\bigr\},
\\
\p_{\al}^{2}\mbbq_{n}(\theta)&=&
\frac{-2T_{n}}{\beta},\qquad \p_{\al}\,\p_{\beta}\mbbq_{n}(
\theta)=-\sum_{j=1}^{n}\frac
{2}{\beta
^{2}}(
\D_{j}X-\al h_{n}),
\\
 \p_{\beta}^{2}\mbbq_{n}(\theta)&=&-\sum
_{j=1}^{n}\frac{2}{\beta^{3}h_{n}} \biggl\{(
\D_{j}X-\al h_{n})^{2}-\frac{\beta h_{n}}{2} \biggr
\}.
\end{eqnarray*}
We can deduce the convergences
\begin{eqnarray*}
\frac{1}{T_{n}}\p_{\al}^{2}\mbbq_{n}(
\theta_{0})&\to^{p}&-2\beta_{0}^{-1},\qquad
\frac{1}{\sqrt{n}\sqrt{T_{n}}}\p_{\al}\,\p_{\beta}\mbbq _{n}(
\theta_{0})\to^{p}0,\\
 \frac{1}{n}\p_{\beta}^{2}
\mbbq_{n}(\theta_{0})&\to^{p}&-
\beta_{0}^{-2},
\end{eqnarray*}
so that the normalized quasi observed-information matrix\break
$-D_{n}^{-1}\p_{\theta}^{2}\mbbq_{n}(\theta_{0})D_{n}^{-1}\hspace*{-1.5pt}\to
^{p}\operatorname{
diag}(2\beta_{0}^{-1},\beta_{0}^{-2})$,
where $D_{n}:=\operatorname{ diag}(\sqrt{T_{n}},\sqrt{n})$.
In view of the classical Cram\'er-type method for $M$-estimation,
we should then have a central limit theorem for the normalized quasi-score
$\{T_{n}^{-1/2}\p_{\al}\mbbq_{n}(\theta_{0}),\break n^{-1/2}\p_{\beta
}\mbbq
_{n}(\theta_{0})\}$
for an asymptotic normality at rate $D_{n}$ to be valid for the
$M$-estimator associated with $\mbbq_{n}$.
However, different from the drifted Wiener process,
the sequence $\{n^{-1/2}\p_{\beta}\mbbq_{n}(\theta_{0})\}$ does
\textit
{not} converge
because $(h_{n}^{-1/2}\D_{j}X)_{j\le n}$ cannot be $L^{q}$-bounded for
large $q>2$
as can be seen from the moment structure of L\'evy processes;
see Luschgy and Pag\`es~\cite{LusPag08} for general moment estimates in
small time with several concrete examples.
Although we only mentioned the L\'evy process with diagonal norming,
the situation remains the same even when $X$ is actually an ergodic
solution to (\ref{SDE}).

The observation made in the last paragraph says that
the situation is different from the case of diffusions,
when developing asymptotic theory concerning the Gaussian
quasi-likelihood for model (\ref{SDE})
under high-frequency sampling framework; it is also different from the
case of time series models,
where the usual $\sqrt{n}$-consistency holds in most cases (see the
references cited in the \hyperref[sec1]{Introduction}).
Earlier attempts to tackle this point have been made by Mancini \cite
{Man04}, Shimizu and Yoshida~\cite{ShiYos06},
Ogihara and Yoshida~\cite{OgiYos11}, where they incorporated
jump-detection filters in defining a contrast function.
The filter approach has its own advantage such as
$\sqrt{n}$-rate estimation of the diffusion parameter, even in the
presence of jumps; however,
we should have in mind that its implementation involves fine-tuning parameters,
thereby possibly preventing us from a straightforward use of the approach.

In order to prove Theorem~\ref{thmm_main}, we will look at not $\theta
\mapsto\mbbq_{n}(\theta)$, but
\[
\theta\mapsto\mbbg_{n}(\theta)= \bigl\{\mbbg^{\al}_{n}(
\theta ),\mbbg ^{\beta}_{n}(\theta) \bigr\},
\]
where $\mbbg^{\al}_{n}\dvtx \Theta\to\mbbr^{p_{\al}}$ and $\mbbg
^{\beta
}_{n}\dvtx \Theta\to\mbbr^{p_{\beta}}$
are defined by
%
%
%e4.1 #&#
%e4.2 #&#
\begin{eqnarray}
\mbbg^{\al}_{n}(\theta)&=&\sum_{j=1}^{n}
\p_{\al}a_{j-1}(\al) \bigl[V_{j-1}^{-1}(
\beta)\bigl[\chi_{j}(\al)\bigr] \bigr], \label{G_def_al}
\\
\mbbg^{\beta}_{n}(\theta)&=&\sum_{j=1}^{n}
\biggl( \bigl\{-\p_{\beta}V_{j-1}^{-1}(\beta)\bigr\}
\bigl[\chi_{j}(\al)^{\otimes2}\bigr] -h_{n}
\frac{\p_{\beta}|V_{j-1}(\beta)|}{|V_{j-1}(\beta)|} \biggr). \label{G_def_be}
\end{eqnarray}
Our contrast function $\mbbm_{n}(\theta)$ is then defined to be the
``squared quasi-score'' as in~(\ref{MC_H_def}),
%
%
%e4.3 #&#
\begin{equation}
\mbbm_{n}(\theta)=-\frac{1}{T_{n}}\bigl|\mbbg_{n}(
\theta)\bigr|^{2}. \label{M_def}
\end{equation}
Trivially, $\mbbg_{n}\dvtx \Theta\to\mbbr^{p}$ fulfil that
$\mbbg_{n}(\theta)= \{(1/2)\p_{\al}\mbbq_{n}(\theta),h_{n}\,\p
_{\beta
}\mbbq_{n}(\theta) \}$.
The difference is that we put the factor ``$h_{n}$'' in front of $\p
_{\beta}\mbbq_{n}(\theta)$;
our estimating procedure is formally not the usual $M$-estimation based on
the Taylor expansion of $\theta\mapsto\mbbq_{n}(\theta)$ around
$\theta_{0}$,
but rather a kind of minimum distance estimation concerning the
Gaussian quasi-score function.
The optimization with respect to $\theta$ is asymptotically the same
for both of $\mbbq_{n}$ and $\mbbm_{n}$:
if there is no root $\theta\in\Theta$ for $\mbbg_{n}(\theta)=0$,
then we may assign any value (e.g., any element of $\Theta$) to $\hat
{\theta}_{n}$,
upholding the claim of Theorem~\ref{thmm_main}.

%re4.1 #&#
\begin{rem}\label{rem_G_exte}
More general cases than (\ref{G_def_al}) and (\ref{G_def_be}) can be
treated, such as
\begin{eqnarray*}
\mbbg^{\al}_{n}(\theta)&=&\sum_{j=1}^{n}
\bar{W}^{\al}_{j-1}(\theta )\bigl\{ X_{t_{j}}-m_{j-1}(
\theta)\bigr\},
\\
\mbbg^{\beta}_{n}(\theta)&=&\sum_{j=1}^{n}
\bigl(\bar{W}^{\beta,1}_{j-1}(\theta)\bigl[\bigl
\{X_{t_{j}}-m_{j-1}(\theta)\bigr\} ^{\otimes2}\bigr]
-h_{n}\bar{W}^{\beta,2}_{j-1}(\theta) \bigr)
\end{eqnarray*}
for some measurable $m\dvtx \mbbr^{d}\times\Theta\to\mbbr^{d}$,
$\bar{W}^{\al}\dvtx \mbbr^{d}\times\Theta\to\mbbr^{p_{\al}}\otimes
\mbbr^{d}$,
$\bar{W}^{\beta,1}\dvtx \mbbr^{d}\times\Theta\to\mbbr^{p_{\beta
}}\otimes
(\mbbr^{d}\otimes\mbbr^{d})$
and $\bar{W}^{\beta,2}\dvtx \mbbr^{d}\times\Theta\to\mbbr^{p_{\beta}}$.
This may be called a GQMLE as well, for we are still solely fitting the
local mean vectors and covariance matrices.
This setting allows us to deal with, for example, the parametric model
\[
dX_{t}=a(X_{t},\theta)\,dt+b(X_{t},
\theta)\,dW_{t}+c(X_{t-},\theta)\,dJ_{t}
\]
with possibly degenerate $b$ and $c$,
the resulting GQMLE $\hat{\theta}_{n}$ still being asymptotically
normal at rate $\sqrt{T_{n}}$
under suitable conditions. To avoid unnecessarily messy notation and
regularity conditions
without losing essence, we have decided to treat~(\ref{SDE}) in this paper.
\end{rem}

For later use, we here introduce some convention and recall a couple of
basic facts
that we will make use often without notice:
\begin{itemize}
\item We will often suppress ``$(\theta_{0})$'' from the notation:
$\chi_{j}:=\chi_{j}(\al_{0})$, $a_{j-1}:=a_{j-1}(\al_{0})$, $\mbbg
_{n}^{\al}=\mbbg_{n}^{\al}(\theta_{0})$, and so forth.
\item$\int_{j}$ denotes a shorthand for $\int_{t_{j-1}}^{t_{j}}$.
\item$M'_{j-1}(\theta):=\p_{\al}a_{j-1}(\al)^{\top
}V_{j-1}^{-1}(\beta
)\in\mbbr^{p_{\al}}\otimes\mbbr^{d}$.
\item$M''_{j-1}(\beta):=-\p_{\beta}V_{j-1}^{-1}(\beta)
=\{V_{j-1}^{-1}(\p_{\beta}V_{j-1})V_{j-1}^{-1}\}(\beta)
\in\mbbr^{p_{\beta}}\otimes\mbbr^{d}\otimes\mbbr^{d}$.
\item$d_{j-1}(\beta):=|V_{j-1}(\beta)|^{-1}\p_{\beta
}|V_{j-1}(\beta
)|\in\mbbr^{p_{\beta}}$.
\item Given real sequence $a_{n}$ and random variables $Y_{n}$ possibly
depending on $\theta$, we write
$Y_{n}=O^{\ast}_{p}(a_{n})$ if $\sup_{n,\theta}E_{0}
[|a_{n}^{-1}Y_{n}|^{K} ]<\infty$ for every $K>0$.

\item$E_{0}^{j-1}[\cdot]:=E_{0}[\cdot|\mcf_{t_{j-1}}]$.

\item$R$ denotes a generic function on $\mbbr^{d}$, possibly depending
on $n$ and $\theta$, for which there exists a constant
$C\ge0$ such that $\sup_{n,\theta}|R(x)|\le C(1+|x|)^{C}$ for every
$x\in\mbbr^{d}$.
\item Burkholder's inequality: for a martingale difference array
$(\zeta
_{nj})_{j\le n}$ and every $q\ge2$,
\[
E_{0} \biggl[\max_{k\le n} \biggl|\sum
_{j\le k}\frac{1}{\sqrt {n}}\zeta _{nj}
\biggr|^{q} \biggr] \lesssim E_{0} \biggl[ \biggl(
\frac{1}{n}\sum_{j\le n}\zeta _{nj}^{2}
\biggr)^{q/2} \biggr] \lesssim\frac{1}{n}\sum
_{i\le n}E \bigl[|\zeta_{nj}|^{q} \bigr].
\]
Moreover, if $b_{\cdot}$ and $c_{\cdot}$ are sufficiently integrable
adapted processes, then
\begin{eqnarray*}
E_{0} \biggl[ \biggl|\int_{0}^{T}b_{s-}\,dW_{s}
\biggr|^{q} \biggr] &\lesssim& T^{q/2-1}\int_{0}^{T}E_{0}
\bigl[|b_{s}|^{q}\bigr]\,ds,
\\
E_{0} \biggl[ \biggl|\int_{0}^{T}c_{s-}\,dJ_{s}
\biggr|^{q} \biggr] &\lesssim&(1\vee T)^{q/2-1}\int
_{0}^{T}E_{0}\bigl[|c_{s}|^{q}
\bigr]\,ds
\end{eqnarray*}
for every $T>0$ and $q\ge2$ such that $E[|J_{1}|^{q}]<\infty$.

\item Sobolev's inequality (e.g., Friedman~\cite{Fri06}, Section~10.2),
\[
E_{0} \Bigl[\sup_{\theta\in\Theta}\bigl|u(\theta)\bigr|^{q}
\Bigr] \lesssim\sup_{\theta\in\Theta} \bigl\{E_{0}\bigl[\bigl|u(
\theta )\bigr|^{q}\bigr]+E_{0}\bigl[\bigl|\p _{\theta}u(
\theta)\bigr|^{q}\bigr] \bigr\}
\]
for $q>p$ and a random field $u\in\mcc^{1}(\Theta)$;
recall that $p$ denotes the dimension of $\theta$
and that we are presupposing the boundedness and convexity of $\Theta$.
We will make use of this type of inequality to derive some
uniform-in-$\theta$ moment estimates for martingale terms.
\end{itemize}

We now turn to the proof of Theorem~\ref{thmm_main} by verifying the
conditions of Theorem~\ref{thmm_PLDI_G}.

%

%s4.1.2 #&#
\subsubsection{Verification of the conditions on $\mbbg_{n}$}

We rewrite $\mbbg_{n}$ as follows:
%
%
%e4.4 #&#
%e4.5 #&#
\begin{eqnarray}
\label{G_al}
\mbbg^{\al}_{n}(\theta)&=&\sum_{j=1}^{n}M'_{j-1}(
\theta)[\chi_{j}] -h_{n}\sum_{j=1}^{n}M'_{j-1}(
\theta) \bigl[a_{j-1}(\al )-a_{j-1} \bigr],
\\
\label{G_beta}
\mbbg^{\beta}_{n}(\theta)&=&\sum_{j=1}^{n}
\bigl\{M''_{j-1}(\beta )\bigl[\chi
_{j}^{\otimes2}\bigr]-h_{n}\,d_{j-1}(\beta)
\bigr\}
\nonumber
\\
&&{} +2h_{n}\sum_{j=1}^{n}M''_{j-1}(
\beta)\bigl[\chi _{j},a_{j-1}-a_{j-1}(\al)\bigr]
\\
&&{} +h_{n}^{2}\sum_{j=1}^{n}M''_{j-1}(
\beta) \bigl[\bigl\{a_{j-1}-a_{j-1}(\al )\bigr\} ^{\otimes2}
\bigr].\nonumber
\end{eqnarray}
We have
%
%
%e4.6 #&#
\begin{equation}
\chi_{j}=\zeta_{j}+r_{j},
\label{chi_def}
\end{equation}
where
%
%
%e4.7 #&#
%e4.8 #&#
\begin{eqnarray}
\zeta_{j}&:=&\int_{j}\tilde{a}_{j-1}(s)\,ds
+\int_{j}b(X_{s},\beta_{0})\,dW_{s}+
\int_{j}c(X_{s-},\beta_{0})\,dJ_{s},
\label{zeta_1}
\\
r_{j}&:=&\int_{j} \bigl\{E_{0}^{j-1}
\bigl[a(X_{s},\al_{0})\bigr]-a_{j-1} \bigr\}\,ds,
\label{zeta_2}
\end{eqnarray}
with $\tilde{a}_{j-1}(s):=a(X_{s},\al_{0})-E_{0}^{j-1}[a(X_{s},\al_{0})]$.
Obviously, $(\zeta_{j})_{j\le n}$ forms a martingale difference array
with respect to
the discrete-time filtration $(\mcf_{t_{j}})_{j\le n}$.

It\^o's formula and the present integrability condition lead to
%
%
%e4.9 #&#
\begin{equation}
E_{0}^{j-1}\bigl[a(X_{s},\al_{0})
\bigr]-a_{j-1}=\int_{j}E_{0}^{j-1}
\bigl[\mca a(X_{u},\al_{0})\bigr]\,du=h_{n}R_{j-1},
\label{eg_add1}
\end{equation}
where $\mca$ denotes the (extended) generator associated with $X$ under
$P_{0}$, that is,
for $f\in\mcc^{2}(\mbbr^{d})$
\begin{eqnarray*}
\mca f(x)&=&\p f(x) \bigl[a(x,\al_{0})\bigr]+\frac{1}{2}
\p^{2}f(x)\bigl[b(x,\beta _{0})^{\otimes2}\bigr]
\\
&&{} +\int \bigl\{f\bigl(x+c(x,\beta_{0})z\bigr)-f(x)-\p f(x)
\bigl[c(x,\beta _{0})z\bigr] \bigr\} \nu(dz).
\end{eqnarray*}
Putting (\ref{zeta_2}) and (\ref{eg_add1}) together gives
$r_{j}=h_{n}^{2}R_{j-1}$, therefore
%
%
%e4.10 #&#
\begin{equation}
\chi_{j}=\zeta_{j}+h_{n}^{2}R_{j-1}.
\label{chi_deco}
\end{equation}

Assumption~\ref{MC_A1} obviously holds under the present
differentiability conditions.
We begin with verifying Assumption~\ref{MC_A2}.

%le4.2 #&#
\begin{lem}\label{lem_check_me}
For every $K>0$, we have
\[
\sup_{n\in\mbbn}E_{0} \biggl[\biggl\llvert
\frac{1}{\sqrt{T_{n}}}\mbbg _{n}(\theta _{0})\biggr\rrvert
^{K} \biggr] +\sup_{n\in\mbbn}E_{0} \biggl[
\sup_{\theta\in\Theta} \biggl\llvert \frac{1}{T_{n}}\mbbg_{n}(
\theta)\biggr\rrvert ^{K} \biggr]<\infty.
\]
\end{lem}

\begin{pf}
By substituting (\ref{chi_deco}) in (\ref{G_al}) and (\ref{G_beta})
and then rearranging the resulting terms, we have
%
%
%e4.11 #&#
%e4.12 #&#
\begin{eqnarray} \label{gme_a1}
\mbbg_{n}^{\al}(\theta)&=&\sum_{j=1}^{n}M'_{j-1}(
\theta)\zeta_{j} +h_{n}\sum_{j=1}^{n}M'_{j-1}(
\theta)\bigl\{a_{j-1}-a_{j-1}(\al)\bigr\}
\nonumber
\\[-8pt]
\\[-8pt]
\nonumber
&&{} +h_{n}^{2}\sum_{j=1}^{n}M'_{j-1}(
\theta)R_{j-1},
\\
\label{gme_b1}
\mbbg^{\beta}_{n}(\theta)&=&\sum_{j=1}^{n}
\bigl\{M''_{j-1}(\beta )\bigl[\zeta
_{j}^{\otimes2}\bigr] -h_{n}d_{j-1}(\beta)
\bigr\}
\nonumber
\\[-8pt]
\\[-8pt]
\nonumber
&&{} +2h_{n}\sum_{j=1}^{n}M''_{j-1}(
\beta)\bigl[\zeta _{j},a_{j-1}-a_{j-1}(\al )
\bigr]+h_{n}^{2}\sum_{j=1}^{n}R_{j-1}.
\end{eqnarray}
To achieve the proof, we will separately look at
$T_{n}^{-1/2}\mbbg^{\al}_{n}$, $T_{n}^{-1/2}\mbbg^{\beta}_{n}$,
$T_{n}^{-1}\mbbg^{\al}_{n}(\theta)$ and $T_{n}^{-1}\mbbg^{\beta
}_{n}(\theta)$.
Fix any integer $K>(2\vee p)$ in the sequel.

First we prove $T_{n}^{-1/2}\mbbg^{\al}_{n}=O^{\ast}_{p}(1)$.
Observe that
\begin{eqnarray*}
\frac{1}{\sqrt{T_{n}}}\mbbg^{\al}_{n} &=&\sum
_{j=1}^{n}\frac{1}{\sqrt{T_{n}}}M'_{j-1}
\zeta_{j}+ \sqrt{T_{n}h_{n}^{2}}
\frac{1}{n}\sum_{j=1}^{n}M'_{j-1}R_{j-1}
\\
&=&\sum_{j=1}^{n}\frac{1}{\sqrt{T_{n}}}M'_{j-1}
\zeta_{j}+O^{\ast
}_{p} \Bigl(\sqrt
{T_{n}h_{n}^{2}} \Bigr).
\end{eqnarray*}
By (\ref{zeta_1}),
%
%
%e4.13 #&#
\begin{eqnarray}\label{M'z_form}
\sum_{j=1}^{n}\frac{1}{\sqrt{T_{n}}}M'_{j-1}
\zeta_{j} &=&\sum_{j=1}^{n}
\frac{1}{\sqrt{n}} \biggl(M'_{j-1}\frac{1}{\sqrt {h_{n}}}\int
_{j}b(X_{s},\beta_{0})\,dW_{s}
\biggr) \nonumber\\
&&{} +\sqrt{h_{n}}\sum_{j=1}^{n}
\frac{1}{\sqrt{n}} \biggl(M'_{j-1}\frac
{1}{h_{n}}\int
_{j}\tilde{a}_{j-1}(s)\,ds \biggr)
\\
&&{} +\sum_{j=1}^{n}\frac{1}{\sqrt{T_{n}}}M'_{j-1}
\int_{j}c(X_{s-},\beta _{0})\,dJ_{s}.\nonumber
\end{eqnarray}
Burkholder's inequality implies that
the first and second term on the right-hand side are $O^{\ast}_{p}(1)$
and $O^{\ast}_{p}(\sqrt{h_{n}})$, respectively.
As for the last term, by writing $\mathbf{1}_{j}\dvtx (0,\infty)\to\{
0,1\}$ for
the identity function of the interval $(t_{j-1},t_{j}]$,
%
%
%e4.14 #&#
\begin{eqnarray}\label{mart_ineq_add1}
&& E_{0} \Biggl[\Biggl |\sum_{j=1}^{n}
\frac{1}{\sqrt {T_{n}}}M'_{j-1}\int_{j}c(X_{s-},
\beta_{0})\,dJ_{s} \Biggr|^{K} \Biggr] \nonumber
\\
&&\qquad\lesssim T_{n}^{-K/2}E_{0} \Biggl[\Biggl |\int
_{0}^{T_{n}}\sum_{j=1}^{n}
\mathbf {1}_{j}(s)M'_{j-1}c(X_{s-},
\beta_{0})\,dJ_{s} \Biggr|^{K} \Biggr]
\nonumber
\\
&&\qquad\lesssim T_{n}^{-K/2}T_{n}^{K/2-1}\int
_{0}^{T_{n}} E_{0} \Biggl[ \Biggl(\sum
_{j=1}^{n}\mathbf {1}_{j}(s)\bigl|M'_{j-1}c(X_{s-},
\beta _{0})\bigr| \Biggr)^{K} \Biggr]\,ds
\\
&&\qquad=\frac{1}{T_{n}}\int_{0}^{T_{n}}\sum
_{j=1}^{n}\mathbf {1}_{j}(s)E_{0}
\bigl[\bigl|M'_{j-1}c(X_{s-},\beta_{0})\bigr|^{K}
\bigr]\,ds
\nonumber
\\
&&\qquad\lesssim\frac{1}{T_{n}}\sum_{j=1}^{n}
\int_{j}\,ds=1,\nonumber
\end{eqnarray}
and hence we are done.

We now prove $T_{n}^{-1/2}\mbbg^{\beta}_{n}=O^{\ast}_{p}(1)$.
In the sequel, we may and do suppose that $d=p_{\beta}=r'=r''=1$: this
reduction is possible because of the
polarization identity
\[
\bigl[S',S''\bigr]=\frac{1}{4}
\bigl(\bigl[S'+S''\bigr]-
\bigl[S'-S''\bigr] \bigr),
\]
which is valid for any two semimartingales $S'$ and $S''$.
By (\ref{chi_deco}) and (\ref{G_beta}),
\[
\frac{1}{\sqrt{T_{n}}}\mbbg^{\beta}_{n}=\sum
_{j=1}^{n}\frac
{1}{\sqrt{T_{n}}} \bigl(M''_{j-1}
\zeta_{j}^{2}-h_{n}d_{j-1} \bigr)
+O^{\ast}_{p} \Bigl(\sqrt{T_{n}h_{n}^{2}}
\Bigr),
\]
so that it remains to verify
%
%
%e4.15 #&#
\begin{equation}
\sum_{j=1}^{n}\frac{1}{\sqrt{T_{n}}}M''_{j-1}
\bigl(\zeta _{j}^{2}-h_{n}V_{j-1}
\bigr)=O^{\ast}_{p}(1). \label{G_beta_est*}
\end{equation}
Define $\zeta_{j}(t)$ for $t\in(t_{j-1},t_{j}]$ by
\begin{eqnarray*}
\zeta_{j}(t)&=&\int_{t_{j-1}}^{t}
\tilde{a}_{j-1}(s)\,ds +\int_{t_{j-1}}^{t}b(X_{s},
\beta_{0})\,dW_{s}\\
&&{}+\int_{t_{j-1}}^{t}c(X_{s-},
\beta_{0})\,dJ_{s}.
\end{eqnarray*}
Let $N(ds,dz)$ denote the Poisson random measure associated with $J$, and $\tilde{N}$ its
compensated version
[i.e., $J_{t}=\int_{0}^{t}\int z\tilde{N}(ds,dz)$].
The quadratic variation at time $t$ is then given as follows
(cf. Jacod and Shiryaev~\cite{JacShi03}, I.4.49(d), I.4.55(b)):
\begin{eqnarray*}
\bigl[\zeta_{j}(\cdot)\bigr]_{t}&=&\int
_{t_{j-1}}^{t}b^{2}(X_{s-},
\beta_{0})\,ds+ \int_{t_{j-1}}^{t}\int
c^{2}(X_{s-},\beta_{0})z^{2}N(ds,dz)
\\
&=&(t-t_{j-1})V_{j-1}+\int_{t_{j-1}}^{t}
\int c^{2}(X_{s-},\beta _{0})\tilde{N}(ds,dz)+
\int_{t_{j-1}}^{t}g_{j-1}(s)\,ds,
\end{eqnarray*}
where we used the assumption $\int z^{2}\nu(dz)=1$ (with the temporary
assumption $r''=1$)
and $g_{j-1}(s):=b^{2}(X_{s},\beta_{0})-b_{j-1}^{2}+c^{2}(X_{s-},\beta
_{0})-c_{j-1}^{2}$.
Applying the integration-by-parts formula, we get
\begin{eqnarray*}
\zeta_{j}^{2}-h_{n}V_{j-1} &=& \biggl
\{2\int_{j}\zeta_{j}(s-)\,d\zeta_{j}(s)
+\int_{j}\int c^{2}(X_{s-},
\beta_{0})z^{2}\tilde{N}(ds,dz)
\\
&&\hspace*{72pt}{} +\int_{j} \bigl(g_{j-1}(s)-E_{0}^{j-1}
\bigl[g_{j-1}(s)\bigr] \bigr)\,ds \biggr\}\\
&&{} +\int_{j}E_{0}^{j-1}
\bigl[g_{j-1}(s)\bigr]\,ds
\\
&=:&\zeta^{(0)}_{j}+\zeta^{(1)}_{j}\qquad
\mbox{say.}
\end{eqnarray*}
We can deduce that $\sum_{j=1}^{n}T_{n}^{-1/2}M''_{j-1}\zeta
^{(0)}_{j}=O^{\ast}_{p}(1)$,
as is the case in the proof of $\sum_{j=1}^{n}T_{n}^{-1/2}M'_{j-1}\zeta
_{j}=O^{\ast}_{p}(1)$
via the expression (\ref{M'z_form}).
Moreover, we can apply It\^o's formula to get $\zeta
^{(1)}_{j}=h_{n}^{2}R_{j-1}$
under the $\mcc^{2}$ property of $x\mapsto(b(x,\beta_{0}), c(x,\beta_{0}))$,
from which it follows that $\sup_{n}E_{0}[|\sum_{j=1}^{n}T_{n}^{-1/2}M''_{j-1}\times \zeta^{(1)}_{j}|^{K}]
\lesssim\sup_{n}(T_{n}h_{n}^{2})^{K/2}<\infty$. We thus get (\ref
{G_beta_est*}).

Let us turn to prove $\sup_{\theta}|T_{n}^{-1}\mbbg^{\al
}_{n}(\theta
)|=O^{\ast}_{p}(1)$.
In the same way as in the proof of $T_{n}^{-1/2}\mbbg^{\al
}_{n}=O^{\ast
}_{p}(1)$,
we can prove $\sum_{j=1}^{n}T_{n}^{-1/2}M'_{j-1}(\theta)\zeta
_{j}=O^{\ast}_{p}(T_{n}^{-1/2})$
for each $\theta\in\Theta$,
since the explicit dependence on $\theta$ is only through the
predictable parts $M'_{j-1}(\theta)$;\vspace*{1pt}
similar arguments will apply in some places below.
Therefore, it follows from (\ref{gme_a1}) that, for each $\theta\in
\Theta$,
\begin{eqnarray}\label{add_eq_Ga}
\frac{1}{T_{n}}\mbbg_{n}^{\al}(\theta) &=&
\frac{1}{\sqrt{T_{n}}} \Biggl(\sum_{j=1}^{n}
\frac{1}{\sqrt {T_{n}}}M'_{j-1}(\theta)\zeta_{j}
\Biggr) +h_{n} \Biggl(\frac{1}{n}\sum
_{j=1}^{n}M'_{j-1}(
\theta)R_{j-1} \Biggr) \nn
\nonumber\\
&&{} +\frac{1}{n}\sum_{j=1}^{n}M'_{j-1}(
\theta)\bigl\{a_{j-1}-a_{j-1}(\al )\bigr\}
\nonumber
\\[-9pt]
\\[-9pt]
\nonumber
&=&O_{p}^{\ast} \biggl(\frac{1}{\sqrt{T_{n}}}\vee h_{n}
\biggr) +\frac{1}{n}\sum_{j=1}^{n}M'_{j-1}(
\theta)\bigl\{a_{j-1}-a_{j-1}(\al)\bigr\}
\nonumber
\\
&=&O_{p}^{\ast} \biggl(\frac{1}{\sqrt{T_{n}}} \biggr) +
\frac{1}{n}\sum_{j=1}^{n}M'_{j-1}(
\theta)\bigl\{a_{j-1}-a_{j-1}(\al)\bigr\},\nonumber
\end{eqnarray}
so that $T_{n}^{-1}\mbbg_{n}^{\al}(\theta)=O^{\ast}_{p}(1)$.
In a quite similar manner, we obtain [see (\ref{pdGa_al}) and~(\ref
{pdGa_beta}) below]
%
%
%e4.16 #&#
\begin{eqnarray}\label{add_eq_Gap}
\qquad&&\frac{1}{T_{n}}\p_{\theta}\mbbg_{n}^{\al}(
\theta)
\nonumber
\\[-8pt]
\\[-8pt]
\nonumber
&&\qquad=O_{p}^{\ast
} \biggl(\frac{1}{\sqrt{T_{n}}} \biggr) +
\frac{1}{n}\sum_{j=1}^{n}
\p_{\theta} \bigl[M'_{j-1}(\theta)\bigl\{
a_{j-1}-a_{j-1}(\al)\bigr\} \bigr] =O^{\ast}_{p}(1).
\end{eqnarray}
Therefore, we arrive at $\sup_{\theta}|T_{n}^{-1}\mbbg^{\al
}_{n}(\theta
)|=O^{\ast}_{p}(1)$
by means of the Sobolev inequality.

It remains to prove $\sup_{\theta}|T_{n}^{-1}\mbbg^{\beta
}_{n}(\theta
)|=O^{\ast}_{p}(1)$;
we remind the reader that we are supposing that $d=p_{\beta}=r'=r''=1$.
As in the proof of (\ref{G_beta_est*}),
we can prove
\[
\sum_{j=1}^{n}\frac{1}{\sqrt{T_{n}}}
\p_{\theta}^{k}M''_{j-1}(
\beta ) \bigl(\zeta_{j}^{2}-h_{n}V_{j-1}
\bigr) =O^{\ast}_{p}(1)
\]
for each $k=0,1$ and $\beta$, so that the Sobolev inequality gives
$\sum_{j=1}^{n}T_{n}^{-1/2}\times  M''_{j-1}(\beta)(\zeta
_{j}^{2}-h_{n}V_{j-1})=O^{\ast}_{p}(1)$.
Therefore, it follows from (\ref{gme_b1}) and simple manipulation that
%
%e4.17 #&#
\begin{eqnarray}\label{add_eq_Gb}
\frac{1}{T_{n}}\mbbg^{\beta}_{n}(\theta) &=&
\frac{1}{\sqrt{T_{n}}} \Biggl(\sum_{j=1}^{n}
\frac{1}{\sqrt {T_{n}}}M''_{j-1}(\beta) \bigl(
\zeta_{j}^{2}-h_{n}V_{j-1}\bigr)
\Biggr) \nonumber
\\
&&{} +\frac{2\sqrt{T_{n}}}{n}\sum_{j=1}^{n}
\frac{1}{\sqrt {T_{n}}}M''_{j-1}(\beta) \bigl
\{a_{j-1}-a_{j-1}(\al) \bigr\}\zeta_{j}
\nonumber
\\
&&{} +\frac{h_{n}}{n}\sum_{j=1}^{n}R_{j-1}
+\frac{1}{n}\sum_{j=1}^{n}M''_{j-1}(
\beta) \bigl\{ V_{j-1}-V_{j-1}(\beta ) \bigr\}
\\
&=&O^{\ast}_{p} \biggl(\frac{1}{\sqrt{T_{n}}}\vee
\frac{\sqrt {T_{n}}}{n}\vee h_{n} \biggr) +\frac{1}{n}\sum
_{j=1}^{n}M''_{j-1}(
\beta) \bigl\{ V_{j-1}-V_{j-1}(\beta ) \bigr\}
\nonumber
\\
&=&O^{\ast}_{p} \biggl(\frac{1}{\sqrt{T_{n}}} \biggr) +
\frac{1}{n}\sum_{j=1}^{n}M''_{j-1}(\nonumber
\beta) \bigl\{ V_{j-1}-V_{j-1}(\beta ) \bigr\}.
\end{eqnarray}
Thus $T_{n}^{-1}\mbbg_{n}^{\beta}(\theta)=O^{\ast}_{p}(1)$.
Quite similarly, we get $T_{n}^{-1}\p_{\theta}\mbbg^{\beta
}_{n}(\theta
)=O^{\ast}_{p}(1)$,
%
%
%e4.18 #&#
\begin{eqnarray}\label{add_eq_Gbp}
\qquad&&\frac{1}{T_{n}}\p_{\theta}\mbbg_{n}^{\beta}(
\theta)
\nonumber
\\[-8pt]
\\[-8pt]
\nonumber
&&\qquad=O_{p}^{\ast
} \biggl(\frac{1}{\sqrt{T_{n}}} \biggr) +
\frac{1}{n}\sum_{j=1}^{n}
\p_{\theta} \bigl[M''_{j-1}(\beta) \bigl
\{ V_{j-1}-V_{j-1}(\beta) \bigr\} \bigr] =O^{\ast}_{p}(1),
\end{eqnarray}
completing the proof.%; see (\ref{pdGb_al}) and (\ref{pdGb_beta}) below.
\end{pf}

Next we turn to verifying the uniform moment estimates in Assumptions~\ref{MC_A3}.
To this end, we prove a preliminary lemma.

%le4.3 #&#
\begin{lem}\label{lem_uem}
Suppose the following conditions:
\begin{itemize}
\item the measurable function $f\dvtx \mbbr^{d}\times\Theta\to\mbbr$
fulfils that
$\theta\mapsto f(x,\theta)$ is differentiable for each $x$ and that
\[
g(x):=\sup_{\theta\in\Theta}\bigl\{\bigl|f(x,\theta)\bigr|\vee\bigl|
\p_{\theta
}f(x,\theta )\bigr|\bigr\}
\]
%
%$g(x):=\sup_{\theta}\{|f(x,\theta)|\vee|\p_{\theta}f(x,\theta)|\}$
is of at most polynomial growth;
\item there exist a probability measure $\pi_{0}$ and a constant $a>0$
such that
$\|P_{t}(x,\cdot)-\pi_{0}(\cdot)\|_{g}\lesssim e^{-at}g(x)$;
\item$\sup_{t}E_{0}[|X_{t}|^{q}]<\infty$ for every $q>0$.
\end{itemize}
Then, for every $K>0$ we have
\[
\sup_{n\in\mbbn}E_{0} \Biggl[\sup_{\theta\in\Theta}
\Biggl|\sqrt{T_{n}} \Biggl( \frac{1}{n}\sum
_{j=1}^{n}f_{j-1}(\theta)-\int f(x,\theta)
\pi _{0}(dx) \Biggr) \Biggr|^{K} \Biggr]<\infty.
\]
\end{lem}

\begin{pf}
Put $n^{-1}\sum_{j=1}^{n}f_{j-1}(\theta)-\int f(x,\theta)\pi
_{0}(dx)=\Lam_{n}'(f;\theta)+\Lam_{n}''(f;\theta)$,
where $\Lam_{n}'(f;\theta):=n^{-1}\sum_{j=1}^{n}\{f_{j-1}(\theta
)-E_{0}[f_{j-1}(\theta)]\}$
and $\Lam_{n}''(f;\theta):=n^{-1}\times \sum_{j=1}^{n}\{
E_{0}[f_{j-1}(\theta
)]-\int f(x,\theta)\pi_{0}(dx)\}$.
Under the present assumptions, we can apply Yoshida~\cite{Yos11}, Lemma
4, to get
$E_{0}[|\p_{\theta}^{k}\Lam_{n}'(f;\theta)|^{K}]\lesssim
T_{n}^{-K/2}+T_{n}^{1-K}\lesssim T_{n}^{-K/2}$
for $k\in\{0,1\}$ and $K\ge2$,
yielding that $\max_{k=0,1}\sup_{\theta}\sup_{n}E_{0}[|\sqrt {T_{n}}\p
_{\theta}^{k}\times \Lam_{n}'(f;\theta)|^{K}]<\infty$.
The Sobolev inequality then gives
\[
\sup_{n\in\mbbn}E_{0} \Bigl[\sup_{\theta\in\Theta}
\bigl\llvert \sqrt {T_{n}}\Lam _{n}'(f;\theta)
\bigr\rrvert ^{K} \Bigr]<\infty.
\]
As for $\Lam''_{n}(f;\theta)$, we have for $k\in\{0,1\}$,
\begin{eqnarray*}
& &\bigl\llvert \sqrt{T_{n}}\p_{\theta}^{k}
\Lam_{n}''(f;\theta)\bigr\rrvert
\\
&&\qquad= \Biggl|\frac{\sqrt{T_{n}}}{n}\sum_{j=1}^{n}
\biggl(\iint\p_{\theta}^{k}f(y,\theta)P_{t_{j-1}}(x,dy)
\eta(dx)\\
&&\hspace*{90pt}{} -\iint\p_{\theta}^{k}f(y,\theta)\pi_{0}(dy)
\eta(dx) \biggr) \Biggr|
\\
&&\qquad= \Biggl|\frac{\sqrt{T_{n}}}{n}\sum_{j=1}^{n}\int
\biggl(\int\p_{\theta}^{k}f(y,\theta) \bigl\{P_{t_{j-1}}(x,dy)-
\pi _{0}(dx) \bigr\} \biggr)\eta(dx) \Biggr|
\\
&&\qquad\le\frac{\sqrt{T_{n}}}{n}\sum_{j=1}^{n}\int
\bigl\llVert P_{t_{j-1}}(x,\cdot)-\pi_{0}(\cdot)\bigr\rrVert
_{g}\eta(dx)
\\
&&\qquad\lesssim\frac{\sqrt{T_{n}}}{n} \sum_{j=1}^{n}
\exp(-at_{j-1})\lesssim\frac{1}{\sqrt{T_{n}}}.
\end{eqnarray*}
This completes the proof.
\end{pf}

%co4.4 #&#
\begin{cor}
Assumption~\ref{MC_A3}(a) holds true.
\end{cor}

\begin{pf}
Again we may and do suppose that $d=p_{\beta}=r'=r''=1$.
Recalling (\ref{add_eq_Ga}), (\ref{add_eq_Gap}), (\ref{add_eq_Gb}) and
(\ref{add_eq_Gbp}),
we apply Lemma~\ref{lem_uem} with
$f(x,\theta)=M'(x,\theta)\{a(x,\al_{0})-a(x,\al)\}$ and $f(x,\theta
)=M''(x,\beta)\{V(x,\beta_{0})-V(x,\beta)\}$
to conclude
\[
\sup_{n\in\mbbn}E_{0} \biggl[\sup_{\theta\in\Theta}
\biggl\llvert \sqrt{T_{n}} \biggl(\frac{1}{T_{n}}
\mbbg_{n}(\theta)-\mbbg_{\infty}(\theta ) \biggr)\biggr\rrvert
^{K} \biggr]<\infty
\]
for every $K>0$, where $\mbbg_{\infty}(\theta):=(\mbbg_{\infty
}^{\al
}(\theta),\mbbg_{\infty}^{\beta}(\theta))$
are given by (\ref{Ga_infty}) and (\ref{Gb_infty}), the integrals in
which are finite by the assumptions.
Trivially $\mbbg_{\infty}(\theta_{0})=0$, and Assumption~\ref{MC_A3}(a)
is verified with
$\chi=\chi_{\al}\wedge\chi_{\beta}$.
\end{pf}

Let us mention the fundamental fact concerning conditional size of
$X$'s increments.
For the convenience of reference we include a sketch of the proof.

%le4.5 #&#
\begin{lem}\label{lem_inc_ME}
Let $g(x):=|a(x,\al_{0})|\vee|b(x,\beta_{0})|\vee|c(x,\beta
_{0})|$, and
fix any $q\ge2$ such that $E[|J_{t}|^{q}]<\infty$. Then
\[
E_{0}^{j-1} \Bigl[\sup_{s\in
[t_{j-1},t_{j}]}|X_{s}-X_{t_{j-1}}|^{q}
\Bigr]\lesssim\cases{ %
h_{n}^{q/2}g^{q}(X_{t_{j-1}}),
& \quad$\mbox{if $c\equiv0$},$
\vspace*{2pt}\cr
h_{n}g^{q}(X_{t_{j-1}}), & \quad $\mbox{otherwise}.$}
\]
In particular, the left-hand side is essentially bounded if so is $g$.
\end{lem}

\begin{pf}
Let $c\not\equiv0$. Given a constant $M>0$,
we let $\tau_{j-1,M}:=\inf\{s\ge t_{j-1}\dvtx |X_{s}|\ge M\}$ and
$\xi_{j-1,M}(s):=E_{0}^{j-1} [\sup\{|X_{u}-X_{t_{j-1}}|^{q}\dvtx u\in[t_{j-1},\break
 s\wedge\tau_{j-1,M}]\} ]$.
We can make use of the Lipschitz property of
the coefficients and Masuda~\cite{Mas05}, Lemma E.1, to derive
$\xi_{j-1,M}(t_{j})\lesssim\int_{t_{j-1}}^{t_{j}}\xi
_{j-1,M}(s)\,ds+h_{n}g^{q}(X_{t_{j-1}})$,
the upper bound being $P_{0}$-a.s. finite according to the definition
of $\tau_{j-1,M}$.
Hence the claim follows on applying Gronwall's inequality and then
letting $M\uparrow\infty$.
The case of $c\equiv0$ is similar.
\end{pf}

We now prove the central limit theorem required in Assumption~\ref{MC_A4}.

%le4.6 #&#
\begin{lem}\label{lem_aux_clt}
We have
%
%
%e4.19 #&#
\begin{equation}
\frac{1}{\sqrt{T_{n}}}\mbbg_{n}(\theta_{0})\to^{\mcl}
\mcn _{p} \bigl(0,\mbbv(\theta_{0}) \bigr),
\label{rev_clt}
\end{equation}
where $\mbbv(\theta_{0})$ is given by (\ref{mbbv_def}).
\end{lem}

\begin{pf}
We begin with extracting the leading martingale terms
of the sequences $T_{n}^{-1/2}\mbbg^{\al}_{n}$ and $T_{n}^{-1/2}\mbbg
^{\beta}_{n}$;
recall the expressions (\ref{gme_a1}) and (\ref{gme_b1}).
Let us rewrite (\ref{zeta_1}) as
%
%
%e4.20 #&#
\begin{equation}
\zeta_{j}=m_{j}+r'_{j},
\label{zeta_def}
\end{equation}
where
\begin{eqnarray*}
m_{j}&:=&b_{j-1}\D_{j}W+c_{j-1}
\D_{j}J,
\\
r'_{j}&:=&\int_{j}
\tilde{a}_{j-1}(s)\,ds +\int_{j}\bigl(b(X_{s},
\beta_{0})-b_{j-1}\bigr)\,dW_{s}+\int
_{j}\bigl(c(X_{s-},\beta _{0})-c_{j-1}
\bigr)\,dJ_{s}.
\end{eqnarray*}
We claim that it suffices to prove that
%
%
%e4.21 #&#
\begin{equation}
\sum_{j=1}^{n}\frac{1}{\sqrt{T_{n}}}\pmatrix{
\tilde{\gamma}^{\al}_{j}
\cr
\tilde{\gamma}^{\beta}_{j} }\to^{\mcl}\mcn_{p} \bigl(0,\mbbv(\theta_{0})
\bigr), \label{lem_aux_clt_1}\vadjust{\goodbreak}
\end{equation}
where $\tilde{\gamma}^{\al}_{j}:=M'_{j-1}m_{j}$ and
$\tilde{\gamma}^{\beta}_{j}:=M''_{j-1}[m_{j}^{\otimes2}]-h_{n}d_{j-1}$,
both of which form martingale difference arrays with respect to $(\mcf
_{t_{j}})_{j\le n}$;
we can verify that $E_{0}^{j-1}[\tilde{\gamma}^{\beta}_{j}[u]]=0$ for
each $u\in\mbbr^{p_{\beta}}$,
making use of the identity
$\operatorname{ trace}\{A(x)^{-1}\p_{x}\times A(x)\}=\p_{x}|A(x)|/|A(x)|$ for a
differentiable square-matrix function $A$.
In fact, recalling what we have seen in the proof of Lemma \ref
{lem_check_me}, we observe the following:
\begin{itemize}
\item We have
\begin{eqnarray*}
\frac{1}{\sqrt{T_{n}}}\mbbg^{\al}_{n} &=&\sum
_{j=1}^{n}\frac{1}{\sqrt{T_{n}}}M'_{j-1}
\biggl(\int_{j}b(X_{s},\beta_{0})\,dW_{s}
+\int_{j}c(X_{s-},\beta_{0})\,dJ_{s}
\biggr)+o_{p}(1)
\\[-1pt]
&=&\sum_{j=1}^{n}\frac{1}{\sqrt{T_{n}}}\tilde{
\gamma}^{\al}_{j} +\sum_{j=1}^{n}
\frac{1}{\sqrt{T_{n}}}M'_{j-1}\int_{j}
\bigl(b(X_{s},\beta _{0})-b_{j-1}
\bigr)\,dW_{s}
\\[-1pt]
&&{} +\sum_{j=1}^{n}\frac{1}{\sqrt{T_{n}}}M'_{j-1}
\int_{j}\bigl(c(X_{s-},\beta _{0})-c_{j-1}
\bigr)\,dJ_{s}+o_{p}(1).
\end{eqnarray*}
By means of Burkholder's inequality and Lemma~\ref{lem_inc_ME} combined
with the conditioning argument,
\begin{eqnarray*}
& &E_{0} \Biggl[ \Biggl|\sum_{j=1}^{n}
\frac{1}{\sqrt{T_{n}}} M'_{j-1}\int_{j}
\bigl(b(X_{s},\beta_{0})-b_{j-1}
\bigr)\,dW_{s} \Biggr|^{2} \Biggr]
\\[-1pt]
&&\qquad\lesssim E_{0} \Biggl[\sum_{j=1}^{n}
\frac{1}{T_{n}}\bigl|M'_{j-1}\bigr|^{2}|R_{j-1}|
\int_{j}h_{n}\,ds \Biggr]\lesssim h_{n}.
\end{eqnarray*}
Following the same line as in (\ref{mart_ineq_add1}), we also get
\[
E_{0} \Biggl[ \Biggl|\sum_{j=1}^{n}T_{n}^{-1/2}M'_{j-1}
\int_{j}\bigl(c(X_{s},\beta_{0})-c_{j-1}
\bigr)\,dJ_{s} \Biggr|^{2} \Biggr]\lesssim h_{n}.
\]
Therefore, it follows that
%
%
%e4.22 #&#
\begin{equation}
\frac{1}{\sqrt{T_{n}}}\mbbg^{\al}_{n}=\sum
_{j=1}^{n}\frac
{1}{\sqrt {T_{n}}}\tilde{
\gamma}^{\al}_{j}+o_{p}(1). \label{Ga+1}
\end{equation}

\item Put $B'_{n}=2\sum_{j=1}^{n}T_{n}^{-1/2}M''_{j-1}[m_{j},r'_{j}]$ and
$B''_{n}=\sum_{j=1}^{n}T_{n}^{-1/2}M''_{j-1}[r'_{j},r'_{j}]$, then
\begin{eqnarray*}
\frac{1}{\sqrt{T_{n}}}\mbbg^{\beta}_{n} &=&\sum
_{j=1}^{n}\frac{1}{\sqrt{T_{n}}}\bigl(M''_{j-1}
\bigl[\zeta_{j}^{\otimes
2}\bigr]-h_{n}d_{j-1}
\bigr)+o_{p}(1)
\\
&=&\sum_{j=1}^{n}\frac{1}{\sqrt{T_{n}}}\tilde{
\gamma}^{\beta
}_{j}+B'_{n}+B''_{n}+o_{p}(1).
\end{eqnarray*}
Since $\sup_{j\le n}E_{0}[|r'_{j}|^{q}]\lesssim h_{n}^{2}$ for every
$q\ge2$
and $E_{0}^{j-1}[|m_{j}|^{2}]\lesssim|R_{j-1}|^{2}h_{n}$, the
Cauchy--Schwarz inequality leads to
\begin{eqnarray*}
E_{0}\bigl[\bigl|B'_{n}\bigr|\bigr]&\lesssim&
\frac{1}{n}\sum_{j=1}^{n}\sqrt{
\frac{n}{h_{n}}} E_{0} \bigl[|R_{j-1}|^{2}E_{0}^{j-1}
\bigl[|m_{j}|^{2}\bigr] \bigr]^{1/2}
E_{0} \bigl[\bigl|r'_{j}\bigr|^{2}
\bigr]^{1/2}
\\
&\lesssim&\sqrt{nh_{n}^{2}}\to0.
\end{eqnarray*}
Moreover, for any $\ep\in(0,1/3)$, H\"older's inequality gives
\begin{eqnarray*}
E_{0}\bigl[\bigl|B''_{n}\bigr|\bigr]&
\lesssim&\frac{1}{n}\sum_{j=1}^{n}
\sqrt{\frac{n}{h_{n}}} E_{0} \bigl[|R_{j-1}|\bigl|r'_{j}\bigr|^{2}
\bigr]
\\
&\lesssim&\frac{1}{n}\sum_{j=1}^{n}
\sqrt{\frac{n}{h_{n}}} E_{0} \bigl[|R_{j-1}|^{(1+\ep)/\ep}
\bigr]^{\ep/(1+\ep)} E_{0} \bigl[\bigl|r'_{j}\bigr|^{2(1+\ep)}
\bigr]^{1/(1+\ep)}
\\
&\lesssim&\frac{1}{n}\sum_{j=1}^{n}
\sqrt{\frac{n}{h_{n}}} E_{0} \bigl[\bigl|r'_{j}\bigr|^{2(1+\ep)}
\bigr]^{1/(1+\ep)} \lesssim\sqrt{nh_{n}^{4/(1+\ep)-1}}
\\
&\lesssim&\sqrt{nh_{n}^{2}}\to0.
\end{eqnarray*}
Hence we have derived
%
%
%e4.23 #&#
\begin{equation}
\frac{1}{\sqrt{T_{n}}}\mbbg^{\beta}_{n}=\sum
_{j=1}^{n}\frac
{1}{\sqrt {T_{n}}}\tilde{
\gamma}^{\beta}_{j}+o_{p}(1). \label{Gb+1}
\end{equation}
\end{itemize}
Having (\ref{Ga+1}) and (\ref{Gb+1}) in hand, it remains to verify
(\ref
{lem_aux_clt_1}).
We are going to apply the classical martingale central limit theorem
(e.g., Dvoretzky~\cite{Dvo77}).

Put $\tilde{\gamma}_{j}=(\tilde{\gamma}^{\al}_{j},\tilde{\gamma
}^{\beta}_{j})$.
It is easy to verify the Lyapunov condition: in fact,
we have $E_{0}^{j-1}[|\tilde{\gamma}_{j}|^{K}]\lesssim h_{n}|R_{j-1}|$
for any $K>2$,
so that $\sum_{j=1}^{n}E_{0}[|T_{n}^{-1/2}\tilde{\gamma
}_{j}|^{K}]\lesssim T_{n}^{1-K/2}\to0$.
It remains to compute the convergence of the quadratic characteristics:
$\sum_{j=1}^{n}E_{0}^{j-1}[\tilde{\gamma}_{j}^{\otimes2}]\to
^{p}\mbbv
(\theta_{0})$.
By means of the Cram\'er--Wold device, it suffices to prove that
for each $v'_{1},v'_{2}\in\mbbr^{p_{\al}}$ and $v''_{1},v''_{2}\in
\mbbr
^{p_{\beta}}$,
%
%
%e4.24 #&#
%e4.25 #&#
%e4.26 #&#
\begin{eqnarray}
\sum_{j=1}^{n}\frac{1}{T_{n}}E_{0}^{j-1}
\bigl[\bigl(\tilde{\gamma}^{\al
}_{j}\bigr)^{\otimes2}
\bigr] \bigl[v'_{1},v'_{2}\bigr] &
\to^{p}&\mbbg_{\infty}^{\prime\al}\bigl[v'_{1},v'_{2}
\bigr], \label {clt_qv_aa}
\\
\mbbv_{\al\beta,n}\bigl[v'_{1},v''_{1}
\bigr]:= \sum_{j=1}^{n}
\frac{1}{T_{n}}E_{0}^{j-1} \bigl[\tilde{
\gamma}^{\al
}_{j}\otimes\tilde{\gamma}^{\beta}_{j}
\bigr] \bigl[v'_{1},v''_{1}
\bigr] &\to^{p}&\mbbv_{\al\beta}\bigl[v'_{1},v''_{1}
\bigr], \label{clt_qv_ba}
\\
\mbbv_{\beta\beta,n}\bigl[v''_{1},v''_{2}
\bigr]:= \sum_{j=1}^{n}
\frac{1}{T_{n}}E_{0}^{j-1} \bigl[\bigl(\tilde{\gamma
}^{\beta
}_{j}\bigr)^{\otimes2} \bigr]
\bigl[v''_{1},v''_{2}
\bigr] &\to^{p}&\mbbv_{\beta\beta}\bigl[v''_{1},v''_{2}
\bigr]. \label{clt_qv_bb}
\end{eqnarray}
First, (\ref{clt_qv_aa}) readily follows by noting
$E_{0}^{j-1}[m_{j}^{\otimes2}]=h_{n}V_{j-1}$ and
applying the ergodic theorem (\ref{LLN}). Next,
%
%
%e4.27 #&#
\begin{eqnarray}\label{clt_qv_ba-1}
&&\mbbv_{\al\beta,n}\bigl[v'_{1},v''_{1}
\bigr]\nonumber\\
&&\qquad=\frac{1}{n}\sum_{j=1}^{n}
\frac
{1}{h_{n}} E_{0}^{j-1} \bigl[M'_{j-1}[m_{j}]
\otimes M''_{j-1}\bigl[m_{j}^{\otimes
2}
\bigr] \bigr] \bigl[v'_{1},v''_{1}
\bigr]
\\
&&\qquad=\frac{1}{n}\sum_{j=1}^{n}
\frac{1}{h_{n}}\sum_{k,l,s} E_{0}^{j-1}
\bigl[m_{j}^{(k)}m_{j}^{(l)}m_{j}^{(s)}
\bigr] \bigl\{M_{j-1}^{\prime(\cdot s)}\otimes M_{j-1}^{\prime\prime(\cdot
kl)}
\bigr\}\bigl[v'_{1},v''_{1}
\bigr].\nonumber
\end{eqnarray}
For later use, we here note that, as $h\to0$,
\[
E \bigl[J_{h}^{(i_{1})}\cdots J_{h}^{(i_{m})}
\bigr]=\cases{ %
 h\nu_{i_{1}i_{2}i_{3}}(3),&\quad $m=3,$
\vspace*{2pt}\cr
h\nu_{i_{1}i_{2}i_{3}i_{4}}(4)+O\bigl(h^{2}\bigr), & \quad $m=4;$}
\]
this can be easily seen through the relation between the mixed moments
and cumulants of $J_{h}$,
where the latter can be computed as the values at $0$ of the partial
derivatives of the cumulant function
$u\mapsto\log E[\exp(iJ_{h}[u])]=h\int\{\exp(iu[z])-1-iu[z]\}\nu(dz)$.
In view of the expression
\[
m_{j}^{(k)}=\sum_{k'}b_{j-1}^{(kk')}
\D_{j}w^{(k')}+\sum_{k''}c_{j-1}^{(kk'')}
\D_{j}J^{(k'')}\
\]
together with the orthogonalities between the increments of $w$ and
$J$, we get
%
%
%e4.28 #&#
\begin{eqnarray} \label{clt_qv_ba-2}
E_{0}^{j-1} \bigl[m_{j}^{(k)}m_{j}^{(l)}m_{j}^{(s)}
\bigr] &=&\sum_{k',l',s'}c_{j-1}^{(kk')}c_{j-1}^{(ll')}c_{j-1}^{(ss')}
E \bigl[\D_{j}J^{(k')}\D_{j}J^{(l')}
\D_{j}J^{(s')} \bigr]
\nonumber
\\
&=&\sum_{k',l',s'}c_{j-1}^{(kk')}c_{j-1}^{(ll')}c_{j-1}^{(ss')}
E \bigl[J^{(k')}_{h_{n}}J^{(l')}_{h_{n}}J^{(s')}_{h_{n}}
\bigr]
\\
&=&h_{n}\sum_{k',l',s'}c_{j-1}^{(kk')}c_{j-1}^{(ll')}c_{j-1}^{(ss')}
\nu _{k'l's'}(3).\nonumber
\end{eqnarray}
(Since $E[J_{1}]=0$, the $3$rd mixed cumulants and the $3$rd mixed
moments of $J_{h_{n}}$ coincides.)
Substituting (\ref{clt_qv_ba-2}) in (\ref{clt_qv_ba-1}), we get (\ref
{clt_qv_ba})
\begin{eqnarray*}
&& \mbbv_{\al\beta,n}\bigl[v'_{1},v''_{1}
\bigr]
\\
&&\qquad=\frac{1}{n}\sum_{j=1}^{n}\sum
_{k,l,s}\sum_{k',l',s'}c_{j-1}^{(kk')}c_{j-1}^{(ll')}c_{j-1}^{(ss')}
\nu_{k'l's'}(3) \bigl\{M_{j-1}^{\prime(\cdot s)}\otimes
M_{j-1}^{\prime\prime(\cdot
kl)} \bigr\}\bigl[v'_{1},v''_{1}
\bigr]
\\
&&\qquad=\frac{1}{n}\sum_{j=1}^{n}\sum
_{k',l',s'}\nu_{k'l's'}(3) \bigl
\{M'_{j-1}\bigl[v'_{1},c^{(\cdot s')}_{j-1}
\bigr] \bigr\} \bigl\{M''_{j-1}
\bigl[v''_{1},c^{(\cdot k')}_{j-1},c^{(\cdot
l')}_{j-1}
\bigr] \bigr\}
\\
&&\qquad\to^{p}\mbbv_{\al\beta}\bigl[v'_{1},v''_{1}
\bigr].
\end{eqnarray*}
Finally, we look at $\mbbv_{\beta\beta,n}$. Direct computation gives
%
%
%e4.29 #&#
%e4.30 #&#
\begin{eqnarray}\label{clt_qv_ba-3}
&&\mbbv_{\beta\beta,n}\bigl[v''_{1},v''_{2}
\bigr] \nn
\\[-2pt]
&&\qquad=\frac{1}{n}\sum_{j=1}^{n}
\frac{1}{h_{n}}E_{0}^{j-1} \bigl[ \bigl(M''_{j-1}
\otimes M''_{j-1} \bigr)\bigl[
\bigl(v''_{1},m_{j}^{\otimes
2}
\bigr),\bigl(v''_{2},m_{j}^{\otimes2}
\bigr)\bigr] \bigr]
\nonumber
\\[-2pt]
&&\qquad{} -\frac{1}{n}\sum_{j=1}^{n}E_{0}^{j-1}
\bigl[ \bigl(d_{j-1}\otimes M''_{j-1}
\bigr) \bigl[v''_{1},\bigl(v''_{2},m_{j}^{\otimes2}
\bigr)\bigr] \bigr]
\nonumber
\\[-2pt]
&&\qquad{} -\frac{1}{n}\sum_{j=1}^{n}E_{0}^{j-1}
\bigl[ \bigl(d_{j-1}\otimes M''_{j-1}
\bigr) \bigl[v''_{2},\bigl(v''_{1},m_{j}^{\otimes2}
\bigr)\bigr] \bigr]
\nonumber
\\[-8pt]
\\[-8pt]
\nonumber
&&\qquad\quad{}+h_{n} \Biggl(\frac{1}{n}\sum
_{j=1}^{n}d_{j-1}^{\otimes
2}
\bigl[v''_{1},v''_{2}
\bigr] \Biggr)
\nonumber
\\[-2pt]
&&\qquad=\frac{1}{n}\sum_{j=1}^{n}
\frac{1}{h_{n}}E_{0}^{j-1} \bigl[\bigl\{M''_{j-1}
\bigl[v''_{1},m_{j}^{\otimes2}
\bigr]\bigr\} \bigl\{ M''_{j-1}
\bigl[v''_{2},m_{j}^{\otimes2}
\bigr]\bigr\} \bigr] +O_{p}(h_{n})
\nonumber\\[-2pt]
&&\qquad=\frac{1}{n}\sum_{j=1}^{n}
\frac{1}{h_{n}}\sum_{k,l,k',l'} M_{j-1}^{\prime\prime(\cdot kl)}
\bigl[v''_{1}\bigr]M_{j-1}^{\prime\prime
(\cdot
k'l')}
\bigl[v''_{2}\bigr] E_{0}^{j-1}
\bigl[m_{j}^{(k)}m_{j}^{(l)}m_{j}^{(k')}m_{j}^{(l')}
\bigr]\nonumber\\[-2pt]
&&\qquad\quad{}+O_{p}(h_{n}).\nonumber
\end{eqnarray}
Using the orthogonality as before
and noting the fact that $E[|w_{h_{n}}|^{4}]=O(h_{n}^{2})$, we get
%
%
%e4.31 #&#
%e4.32 #&#
\begin{eqnarray}\label{clt_qv_ba-4}
&& E_{0}^{j-1} \bigl[m_{j}^{(k)}m_{j}^{(l)}m_{j}^{(k')}m_{j}^{(l')}
\bigr] \nn
\\[-2pt]
&&\qquad=\sum_{s,t,s',t'}c_{j-1}^{(ks)}c_{j-1}^{(lt)}c_{j-1}^{(k's')}c_{j-1}^{(l't')}
E \bigl[J^{(s)}_{h_{n}}J^{(t)}_{h_{n}}J^{(s')}_{h_{n}}J^{(t')}_{h_{n}}
\bigr]+R_{j-1}h_{n}^{2}
\nonumber
\\[-8pt]
\\[-8pt]
\nonumber
&&\qquad=h_{n}\sum_{s,t,s',t'}c_{j-1}^{(ks)}c_{j-1}^{(lt)}c_{j-1}^{(k's')}c_{j-1}^{(l't')}
\bigl\{\nu_{sts't'}(4)+O(h_{n})\bigr\}+R_{j-1}h_{n}^{2}
\\[-2pt]
&&\qquad=h_{n}\sum_{s,t,s',t'}c_{j-1}^{(ks)}c_{j-1}^{(lt)}c_{j-1}^{(k's')}c_{j-1}^{(l't')}
\nu_{sts't'}(4)+R_{j-1}h_{n}^{2}.\nonumber
\end{eqnarray}
%
%(The mixed moments coming from the terms of
%the form ``$(\D_{j}w^{(s)})^{2}(\D_{j}J^{(t)})^{2}$'' are not null,
%but have no contribution.)
By putting (\ref{clt_qv_ba-3}) and (\ref{clt_qv_ba-4}) together, we get
(\ref{clt_qv_bb})
\begin{eqnarray*}
&& \mbbv_{\beta\beta,n}\bigl[v''_{1},v''_{2}
\bigr]
\\[-2pt]
&&\quad=\frac{1}{n}\sum_{j=1}^{n}\sum
_{s,t,s',t'}\nu_{sts't'}(4) \bigl
\{M''_{j-1}\bigl[v''_{1},c^{(\cdot s)}_{j-1},c^{(\cdot
t)}_{j-1}
\bigr] \bigr\} \bigl\{M''_{j-1}
\bigl[v''_{2},c^{(\cdot s')}_{j-1},c^{(\cdot
t')}_{j-1}
\bigr] \bigr\}+O_{p}(h_{n})
\\[-2pt]
&&\quad\to^{p}\mbbv_{\beta\beta}\bigl[v''_{1},v''_{2}
\bigr].
\end{eqnarray*}
The proof is thus complete.\vadjust{\goodbreak}
\end{pf}

%

%s4.1.3 #&#
\subsubsection{Verification of the conditions on the derivatives of
$\mbbg_{n}$}

Based on (\ref{G_al}) and (\ref{G_beta}), we derive the following
bilinear forms:
% to allow break between pages
%
%
%e4.33 #&#
%e4.34 #&#
%e4.35 #&#
%e4.36 #&#
%e4.37 #&#
\begin{eqnarray}
 \label{pdGa_al}\qquad\p_{\al}\mbbg^{\al}_{n}(\theta)&=&\sum
_{j=1}^{n}\p_{\al
}M'_{j-1}(
\theta )[\chi_{j}] -h_{n}\sum_{j=1}^{n}\,
\p_{\al}M'_{j-1}(\theta)\bigl[a_{j-1}(
\al)-a_{j-1}\bigr] \nn
\nonumber
\\[-8pt]
\\[-8pt]
\nonumber
&&{} -h_{n}\sum_{j=1}^{n}M'_{j-1}(
\theta)\p_{\al}a_{j-1}(\al),
\\
 \label{pdGa_beta}\p_{\beta}\mbbg^{\al}_{n}(\theta)&=&\sum
_{j=1}^{n}\,\p_{\beta
}M'_{j-1}(
\theta)[\chi_{j}]
\nonumber
\\[-8pt]
\\[-8pt]
\nonumber
 &&{}-h_{n}\sum_{j=1}^{n}
\p_{\beta}M'_{j-1}(\theta)\bigl[a_{j-1}(
\al)-a_{j-1}\bigr],
\\
\label{pdGb_al}\p_{\al}\mbbg^{\beta}_{n}(\theta)&=&
-2h_{n}\sum_{j=1}^{n} \bigl
\{M''_{j-1}(\beta)\p_{\al}a_{j-1}(
\al ) \bigr\}
\nonumber
\\[-8pt]
\\[-8pt]
\nonumber
&&{}\times \bigl[\chi_{j}-h_{n}\bigl\{a_{j-1}(
\al)-a_{j-1}\bigr\}\bigr],
\\
\label{pdGb_beta}
\p_{\beta}\mbbg^{\beta}_{n}(\theta)&=& \sum
_{j=1}^{n} \bigl\{\p_{\beta}M''_{j-1}(
\beta)\bigl[\chi_{j}^{\otimes
2}\bigr]-h_{n}
\p_{\beta}d_{j-1}(\beta) \bigr\} \nn
\\
&&{} -2h_{n}\sum_{j=1}^{n}
\,\p_{\beta}M''_{j-1}(\beta)\bigl[\chi
_{j},a_{j-1}(\al )-a_{j-1}\bigr]
\\
&&{} + h_{n}^{2}\sum_{j=1}^{n}\,
\p_{\beta}M''_{j-1}(\beta) \bigl[\bigl\{
a_{j-1}(\al )-a_{j-1}\bigr\}^{\otimes2} \bigr].
\nonumber
\end{eqnarray}
We can prove the following lemma in a similar way to the proof of Lemma~\ref{lem_check_me}.

%le4.7 #&#
\begin{lem}
For every $K>0$,
\[
\sup_{n}E_{0} \biggl[\sup_{\theta}
\biggl\llvert \frac{1}{T_{n}}\p_{\theta
}^{k}
\mbbg_{n}(\theta)\biggr\rrvert ^{K} \biggr] <\infty,\qquad
k=1,2,3.
\]
\end{lem}

Recall that the matrix $\mbbg_{\infty}'(\theta_{0})=\operatorname{ diag}
\{\mbbg_{\infty}^{\prime\al}(\theta_{0}),\mbbg_{\infty}^{\prime
\beta
}(\theta_{0})\}$
is given by (\ref{G'a}) and~(\ref{G'b}).
%
%le4.8 #&#
\begin{lem}\label{pG_ume}
For every $K>0$,
\[
\sup_{n\in\mbbn}E_{0} \biggl[ \biggl\llvert
\sqrt{T_{n}} \biggl(\frac{1}{T_{n}}\p_{\theta}\mbbg
_{n}(\theta_{0}) -\mbbg_{\infty}'(
\theta_{0}) \biggr)\biggr\rrvert ^{K} \biggr]<\infty.
\]
\end{lem}

\begin{pf}
First, concerning the off-diagonal parts, we have
\begin{eqnarray*}
\frac{1}{T_{n}}\p_{\beta}\mbbg_{n}^{\al} &=&
\frac{1}{\sqrt{T_{n}}}\sum_{j=1}^{n}
\frac{1}{\sqrt{T_{n}}} \p_{\beta}M'_{j-1}[
\chi_{j}]=O^{\ast}_{p} \biggl(\frac{1}{\sqrt {T_{n}}}
\biggr),
\\
\frac{1}{T_{n}}\p_{\al}\mbbg_{n}^{\beta} &=&-2
\frac{h_{n}}{\sqrt{T_{n}}}\sum_{j=1}^{n}
\frac{1}{\sqrt{T_{n}}} M''_{j-1} [
\p_{\al}a_{j-1},\chi_{j} ]=O^{\ast}_{p}
\biggl(\frac
{h_{n}}{\sqrt{T_{n}}} \biggr),
\end{eqnarray*}
where the moment estimates for the martingale terms
will be proved in an analogous way to the proof of Lemma~\ref{lem_check_me}.
Next, we observe
\begin{eqnarray*}
\frac{1}{T_{n}}\p_{\al}\mbbg_{n}^{\al}-
\mbbg_{\infty}^{\prime
\al} &=&\frac{1}{\sqrt{T_{n}}}\sum
_{j=1}^{n}\frac{1}{\sqrt{T_{n}}}\p _{\al
}M'_{j-1}[
\chi_{j}] -\frac{1}{n}\sum_{j=1}^{n}M'_{j-1}
\p_{\al}a_{j-1}-\mbbg_{\infty
}^{\prime
\al}(
\theta_{0})
\\
&=&O^{\ast}_{p} \biggl(\frac{1}{\sqrt{T_{n}}} \biggr)+
\frac{1}{\sqrt{T_{n}}} \Biggl\{\!\sqrt{T_{n}} \Biggl(\!-\frac{1}{n}\sum
_{j=1}^{n}M'_{j-1}\p
_{\al}a_{j-1} -\mbbg_{\infty}^{\prime\al}(
\theta_{0}) \!\Biggr)\! \Biggr\}
\nonumber
\\
&=&O^{\ast}_{p} \biggl(\frac{1}{\sqrt{T_{n}}} \biggr),
\end{eqnarray*}
where we used Lemma~\ref{lem_uem} for the last equality.
It remains to look at $T_{n}^{-1}\p_{\beta}\mbbg_{n}^{\beta}$.
Plugging in the identity $\chi_{j}=m_{j}+r'_{j}+h_{n}^{2}R_{j-1}$ and
making use of what we have seen in the first half of the proof of Lemma
\ref{lem_aux_clt}, we proceed as follows:
%
%
%e4.38 #&#
%e4.39 #&#
\begin{eqnarray}\label{pG_ume_1}
\frac{1}{T_{n}}\p_{\beta}\mbbg_{n}^{\beta} \nn
&=&\frac{1}{T_{n}}\sum_{j=1}^{n} \bigl(
\p_{\beta
}M''_{j-1}\bigl[
\bigl(m_{j}+r'_{j}\bigr)^{\otimes2}\bigr]
-h_{n}\p_{\beta}d_{j-1} \bigr)+O^{\ast}_{p}(h_{n})
\nonumber
\\
&=&\frac{1}{T_{n}}\sum_{j=1}^{n} \bigl(
\p_{\beta
}M''_{j-1}\bigl[m_{j}^{\otimes2}
\bigr] -h_{n}\p_{\beta}d_{j-1} \bigr)+O^{\ast}_{p}
(\sqrt {h_{n}} )
\nonumber
\\
&=&\frac{1}{\sqrt{T_{n}}} \Biggl\{\sum_{j=1}^{n}
\frac{1}{\sqrt{T_{n}}} \bigl(\p_{\beta}M''_{j-1}
\bigl[m_{j}^{\otimes2}\bigr] -E_{0}^{j-1}
\bigl[\p_{\beta}M''_{j-1}
\bigl[m_{j}^{\otimes2}\bigr] \bigr] \bigr) \Biggr\}
\nonumber
\\
&&{}+\frac{1}{T_{n}}\sum_{j=1}^{n}
\bigl(E_{0}^{j-1} \bigl[\p_{\beta}M''_{j-1}
\bigl[m_{j}^{\otimes2}\bigr] \bigr]-h_{n}
\p_{\beta}d_{j-1} \bigr) +O^{\ast}_{p} (
\sqrt{h_{n}} )
\\
&=&\frac{1}{T_{n}}\sum_{j=1}^{n}
\bigl(E_{0}^{j-1} \bigl[\p_{\beta}M''_{j-1}
\bigl[m_{j}^{\otimes2}\bigr] \bigr]-h_{n}
\p_{\beta}d_{j-1} \bigr) +O^{\ast}_{p} \biggl(
\frac{1}{\sqrt{T_{n}}} \biggr)
\nonumber
\\
&=&\frac{1}{n}\sum_{j=1}^{n} \bigl[\operatorname{trace} \bigl\{ \bigl(-\p_{\beta_{l}}\,\p_{\beta_{l'}}V_{j-1}^{-1}
\bigr)V_{j-1} \bigr\} -\p_{\beta_{l}}\,\p_{\beta_{l'}}
\log|V_{j-1}| \bigr]_{l,l'=1}^{p_{\beta}}\nonumber\\
&&{}+O^{\ast}_{p}
\biggl(\frac{1}{\sqrt {T_{n}}} \biggr). \nonumber
\end{eqnarray}
The $(l,l')$th component of the first term in (\ref{pG_ume_1}) tends in
probability to
\begin{eqnarray*}
&& \int \bigl[\operatorname{ trace} \bigl\{-\p_{\beta_{l}}\,\p_{\beta_{l'}}V^{-1}V(x,
\beta_{0}) \bigr\} -\p_{\beta_{l}}\,\p_{\beta_{l'}}\log|V|(x,
\beta_{0}) \bigr]\pi_{0}(dx)
\\
&&\qquad =-\int\operatorname{ trace} \bigl\{ \bigl(V^{-1}(\p_{\beta_{l}}V)V^{-1}(
\p_{\beta_{l'}}V) \bigr) (x,\beta _{0}) \bigr\}
\pi_{0}(dx).
\end{eqnarray*}
Accordingly, a reduced version of Lemma~\ref{lem_uem} with $\Theta=\{
\theta_{0}\}$ applies to conclude that
$T_{n}^{-1}\p_{\beta}\mbbg_{n}^{\beta}(\theta_{0})-\mbbg_{\infty
}^{\prime\beta}(\theta_{0})
=O^{\ast}_{p}(T_{n}^{-1/2})$. The proof is complete.
\end{pf}

%%%

%s4.2 #&#
\subsection{\texorpdfstring{Proof of Corollary \protect\ref{cor_main}}{Proof of Corollary 2.8}}\label{sec_proof2}

By Theorem~\ref{thmm_main}, we know that $\sqrt{T_{n}}(\hat{\al
}_{n}-\al
_{0})=O_{p}(1)$ and
$\sqrt{T_{n}}(\hat{\beta}_{n}-\beta_{0})=O_{p}(1)$.
It is easy to see from Taylor expansion that
$\hat{\mbbg}_{n}^{\prime\al}\to^{p}\mbbg_{\infty}^{\prime\al
}(\theta
_{0})$ and
$\hat{\mbbg}_{n}^{\prime\beta}\to^{p}\mbbg_{\infty}^{\prime
\beta}(\theta_{0})$.
Turning to $\hat{\mbbv}_{\al\beta,n}$ and $\hat{\mbbv}_{\beta
\beta,n}$,
we plug the expression
$\chi_{j}(\hat{\al}_{n})=\chi_{j}+\sqrt{h_{n}/n}R_{j-1}[\sqrt {T_{n}}(\hat{\al}_{n}-\al_{0})]$
into their definitions and then apply Taylor expansion
with respect to $\hat{\theta}_{n}$ around $\theta_{0}$ as before, to obtain
%
%
%e4.40 #&#
%e4.41 #&#
\begin{eqnarray}\label{cor_main_1}
\hat{\mbbv}_{\al\beta,n}\bigl[v'_{1},v''_{1}
\bigr] &=&-\sum_{j=1}^{n}\frac{1}{T_{n}}
\bigl(V_{j-1}^{-1}\otimes\p_{\beta}V_{j-1}^{-1}
\bigr) \bigl[ \bigl(\p_{\al}a_{j-1}\bigl[v'_{1}
\bigr],\chi_{j} \bigr),\bigl(v''_{1},
\chi _{j}^{\otimes2}\bigr) \bigr] \nn
\\
&&{} +O_{p} \biggl(\frac{1}{\sqrt{T_{n}}} \biggr),
\nonumber
\\[-8pt]
\\[-8pt]
\nonumber
\hat{\mbbv}_{\beta\beta,n}\bigl[v''_{1},v''_{2}
\bigr]&=& \sum_{j=1}^{n}\frac{1}{T_{n}}
\bigl(\p_{\beta}V_{j-1}^{-1} \otimes
\p_{\beta}V_{j-1}^{-1} \bigr) \bigl[
\bigl(v''_{1},\chi_{j}^{\otimes2}
\bigr),\bigl(v''_{2},\chi_{j}^{\otimes
2}
\bigr) \bigr]
\\
&&{} +O_{p} \biggl(\frac{1}{\sqrt{T_{n}}} \biggr).
\nonumber
\end{eqnarray}
We only show that $\hat{\mbbv}_{\al\beta,n}[v'_{1},v''_{2}]\to
^{p}\mbbv
_{\al\beta}[v'_{1},v''_{1}]$,
for the case of $\hat{\mbbv}_{\beta\beta,n}$ is similar.

Write $\sum_{j=1}^{n}\eta_{j}$ for the first term in the right-hand
side of (\ref{cor_main_1}).
We can show that
\[
\sum_{j=1}^{n}E_{0}^{j-1}[
\eta_{j}]\to^{p}\mbbv_{\al\beta}\bigl[v'_{1},v''_{1}
\bigr]
\]
in a similar manner to show the convergence of the quadratic characteristics
in the proof of Lemma~\ref{lem_aux_clt}.
Noting that $E_{0}^{j-1}[|\chi_{j}|^{q}]\le h_{n}R_{j-1}$ for every
$q\ge2$, we also have
\[
\sum_{j=1}^{n}E_{0} \bigl[
\bigl(\eta_{j}-E_{0}^{j-1}[\eta _{j}]
\bigr)^{2} \bigr] \lesssim\sum_{j=1}^{n}E_{0}
\bigl[\eta_{j}^{2}\bigr]\lesssim\frac
{1}{T_{n}}\to0.
\]
Applying the Lenglart domination property for the martingale $\sum_{j=1}^{n}(\eta_{j}-E_{0}^{j-1}[\eta_{j}])$
(cf. Jacod and Shiryaev~\cite{JacShi03}, I.3.30),
we conclude that $\sum_{j=1}^{n}\eta_{j}\to^{p}\mbbv_{\al\beta
}[v'_{1},v''_{1}]$,
hence $\hat{\mbbv}_{\al\beta,n}[v'_{1},v''_{1}]\to^{p}\mbbv_{\al
\beta
}[v'_{1},v''_{1}]$.

%%%

%s4.3 #&#
\subsection{\texorpdfstring{Proof of Theorem \protect\ref{thmm_AN}}{Proof of Theorem 2.9}}
\label{sec_proof3}

First, we mention an auxiliary estimate. Recall~(\ref{chi_def}) and
(\ref{zeta_def}):
$\chi_{j}:=\D_{j}X-h_{n}a_{j-1}(\al_{0})=m_{j}+(r_{j}+r'_{j})$.
Using Birkholder's inequality and then the Lipschitz continuity of the
coefficients, we see that
\[
E_{0} \bigl[\bigl|r_{j}+r'_{j}\bigr|^{q'}
\bigr]\lesssim\int_{j}E_{0}
\bigl[|X_{s}-X_{t_{j-1}}|^{q'} \bigr]\,ds \lesssim
h_{n}^{2}\|g\|_{\infty}^{q'}\lesssim
h_{n}^{2}
\]
for $q'\in[2,q]$, where $g$ is the one given in Lemma~\ref{lem_inc_ME}.
In this proof, $R$ denotes a generic essentially bounded function on
$\mbbr^{d}$
possibly depending on $n$ and~$\theta$.

By means of the classical $M$-estimation theory (e.g., van der Vaart
\cite{vdV98}, Chapter~5),
it is crucial to have the uniform convergence
%
%
%e4.42 #&#
\begin{equation}
\sup_{\theta\in\Theta}\biggl\llvert \frac{1}{T_{n}}
\mbbg_{n}(\theta )-\mbbg _{\infty}(\theta)\biggr\rrvert +\sup
_{\theta\in\Theta}\biggl\llvert \frac{1}{T_{n}}\p_{\theta}\mbbg
_{n}(\theta )-\mbbg'_{\infty}(\theta)\biggr\rrvert
\to^{p}0. \label{rem_add_rev1_eq1}
\end{equation}
Most key materials to prove this have been obtained
in the proof of Theorem~\ref{thmm_main}, so we only give a sketch.

Note that the variables $M'_{j-1}(\theta)$ and $M''_{j-1}(\beta)$
are now essentially bounded uniformly in $\theta$.
Substituting $\chi_{j}=m_{j}+h_{n}^{2}R_{j-1}$ in
the expressions (\ref{G_al}) and (\ref{G_beta}) about $\mbbg_{n}$,
and also
(\ref{pdGa_al}), (\ref{pdGa_beta}), (\ref{pdGb_al}) and (\ref
{pdGb_beta}) about $\p_{\theta}\mbbg_{n}$,
it is not difficult to deduce (\ref{rem_add_rev1_eq1}); as was in the
proof of Theorem~\ref{thmm_main},
for the estimate to be valid uniformly in $\theta$
we applied Sobolev inequality in part, where it was needed that
$E[|J_{1}|^{q}]<\infty$ for some $q>p$.

Now, the consistency of $\hat{\theta}_{n}$ follows from (\ref
{rem_add_rev1_eq1}):
$\hat{\theta}_{n}\to^{p}\theta_{0}$.
Since $P[\omega\dvtx\break  \mbbg_{n}(\hat{\theta}_{n}(\omega))=0]\to1$,
we may and do suppose that $\mbbg_{n}(\hat{\theta}_{n})=0$.
In view of (\ref{rem_add_rev1_eq1}) and the Taylor expansion
$0=T_{n}^{-1/2}\mbbg_{n}(\theta_{0})
+T_{n}^{-1}\,\p_{\theta}\mbbg_{n}(\tilde{\theta}_{n})[\sqrt {T_{n}}(\hat
{\theta}_{n}-\theta_{0})]$,
where the point $\tilde{\theta}_{n}$ lies on the segment connecting
$\hat{\theta}_{n}$ and $\theta_{0}$,
it suffices to have the central limit theorem (\ref{rev_clt}).
By close inspection of the proof of Lemma~\ref{lem_aux_clt},
we note that the present assumption [especially $q>(4\vee p)$ about the
moment order]
is enough to conclude (\ref{rev_clt}).
The proof is complete.

%%%%%
%%%%%

%s5 #&#
\section{A criterion for the exponential ergodicity in dimension one}\label{sec_proof_prop_ergo}

In this section, we set $d=r'=r''=1$ and suppress dependence on the
parameter from the notation
%
%
%e5.1 #&#
\begin{equation}
dX_{t}=a(X_{t})\,dt+b(X_{t})\,dW_{t}+c(X_{t-})\,dJ_{t}.
\label{SDE'}
\end{equation}
We here forget Assumptions~\ref{A_J} to~\ref{A_nond}, and instead
introduce the following set of conditions.

%as5.1 #&#
\begin{ass}\label{A_E1}
$(a,b,c)$ is of class $\mcc^{1}(\mbbr)$ and globally Lipschitz, and
$(b,c)$ is bounded.\vadjust{\goodbreak}
\end{ass}

%as5.2 #&#
\begin{ass}\label{A_E2}
Either one of the following conditions holds true:
\begin{longlist}[{(ii)}]
\item[{(i)}] $b(x')\ne0$ for some $x'$, $c(x'')\ne0$ for every $x''$,
and there exists a constant $\overline{\ep}>0$
such that $\nu(-\ep,0)\wedge\nu(0,\ep)>0$ for every $\ep\in
(0,\overline
{\ep})$;
\item[{(ii)}] $b\equiv0$, $c(x'')\ne0$ for every $x''$, and we
have the decomposition
\[
\nu=\nu_{\star}+\nu_{\natural}
\]
for two L\'evy measures $\nu_{\star}$ and $\nu_{\natural}$,
where the restriction of $\nu_{\star}$ to some open set of the form
$(-\overline{\ep},0)\cup(0,\overline{\ep})$
admits a continuously differentiable positive density $g_{\star}$.
\end{longlist}
\end{ass}

%as5.3 #&#
\begin{ass}\label{A_E3}
\begin{longlist}[{(ii)}]
\item[{(i)}] $E[J_{1}]=0$ and $\int_{|z|>1}|z|^{q}\nu(dz)<\infty$
for some $q\ge1$, and
\[
\limsup_{|x|\to\infty}\frac{a(x)}{x}<0.
\]
\item[{(ii)}] $E[J_{1}]=0$ and $\int_{|z|>1}\exp(q|z|)\nu
(dz)<\infty
$ for some $q>0$, and
\[
\limsup_{|x|\to\infty}\operatorname{ sgn}(x)a(x)<0.
\]
\end{longlist}
\end{ass}

The next proposition gives a pretty simple criterion for Assumption
\ref
{A_ergo}.

%pr5.4 #&#
\begin{prop}\label{prop_ergo}
The following holds true:
\begin{longlist}[(a)]
\item[(a)]
Suppose conditions~\ref{A_E1},~\ref{A_E2},~\ref{A_E3}\textup{(i)}, and that
$E[|X_{0}|^{q}]<\infty$.
Then, there exist a probability measure $\pi$ and a constant $a>0$
such that
(\ref{exp_ergo_def}) holds true for
a $\mcc^{2}$-function $g$ satisfying that $g(x)=1+|x|^{q}$ outside a
neighborhood of the origin.
Further, (\ref{g_moment_bound}) holds true for the $q$ given in~\ref{A_E3}\textup{(i)}.

\item[(b)] Suppose~\ref{A_E1},~\ref{A_E2},~\ref{A_E3}\textup{(ii)}, and that
$E[\exp(q|X_{0}|)]<\infty$.
Then, there exist a probability measure $\pi$ and constants $a,\ep>0$
such that (\ref{exp_ergo_def})
holds true for a $\mcc^{2}$-function $g$ satisfying that
$g(x)=1+\exp(\ep|x|)$ outside a neighborhood of the origin.
Further, (\ref{g_moment_bound}) holds true for arbitrary $q>0$.
\end{longlist}
\end{prop}

\begin{pf}
The Lipschitz continuity implies that the SDE (\ref{SDE'}) admits a
unique strong solution.
We consider the following conditions:
\begin{longlist}[(II)]
\item[(I)] there exists a constant $\D>0$ for which
every compact sets are petite for the Markov chain $(X_{j\D})_{j\in
\mbbzp}$;
\item[(II)] the exponential Lyapunov-drift criterion
%
%
%e5.2 #&#
\begin{equation}
\mca\vp\le-c\vp+d \label{edc}
\end{equation}
holds true for some constants $c,d>0$ and some $\vp\dvtx \mbbr\to\mbbrp$
belonging to the domain of $\mca$
such that $\lim_{|x|\to\infty}\vp(x)=\infty$, where $\mca$
denotes the
extended generator of $X$.
\end{longlist}
As in the proof of Masuda~\cite{Mas07}, the proof of Theorem 2.2,
in each of (a) and (b) the exponential ergodicity (\ref{exp_ergo_def})
follows from (I) and (II),
and the moment bound~(\ref{g_moment_bound}) from (II) alone.
In order to prove (I), we will first verify the \textit{Local Doeblin
(LD) condition}
(see Kulik~\cite{Kul09} for details); we note that the LD condition
implies (I) for any $\D>0$.
Then we will verify the drift condition (II)
with different choices of $\vp$ under Assumptions~\ref{A_E3}(i) and~\ref{A_E3}(ii).

\textit{Verification of (I): the LD condition.}

First, we verify the LD condition under Assumption~\ref{A_E2}(i).
Let $\Pi_{x}(A):=\nu(\{z\in\mbbr\dvtx  c(x)z\in A\})$,
and refer to Kulik's condition (\ref{eqS}) in the reduced form
{\renewcommand{\theequation}{$\mathrm{S}$}
%e5.3 #&#
\begin{eqnarray}\label{eqS}
&& \forall x\in\mbbr\ \forall v\in\{-1,1\}\ \exists\rho\in(-1,1)\ \forall\del>0
\dvtx \nn
\nonumber
\\[-8pt]
\\[-8pt]
\nonumber
&&\qquad{} \Pi_{x} \bigl(\bigl\{y\in\mbbr\dvtx yv\ge\rho|y|\bigr\}\cap\bigl\{y\in\mbbr
\dvtx |y|\le\del \bigr\} \bigr)>0. %\tag{$$}
\end{eqnarray}}
\hspace*{-2pt}Under Assumption~\ref{A_E2}(i), it follows form Kulik~\cite{Kul09}, Theorem
1.3, Proposition~A.2 and Proposition 4.7,
that the condition (\ref{eqS}) above implies the LD condition.
Simple manipulation shows that the last condition is equivalent to the
following:
\begin{eqnarray*}
&& \forall x\in\mbbr\ \forall\del>0\dvtx
\\
&&\qquad{} \nu \bigl(\bigl\{z\in\mbbr\dvtx 0\le c(x)z\le\del\bigr\} \bigr) \wedge \nu
\bigl(\bigl\{z\in\mbbr\dvtx -\del\le c(x)z\le0\bigr\} \bigr)>0.
\end{eqnarray*}
Since $\nu(\mbbr)>0$, it suffices to look at $x$ such that $c(x)\ne0$.
However, for such~$x$, the condition obviously holds true under
Assumption~\ref{A_E2}(i).

Next we verify the LD condition under Assumption~\ref{A_E2}(ii).
If $c$ is constant, then we can apply Kulik~\cite{Kul09}, Proposition
0.1,
to verify the LD condition. Therefore, we suppose that $\p_{x}c\not
\equiv0$ in what follows.
We smoothly truncate the support of $\nu_{\star}$ as follows:
pick any $\underline{\ep}\in(0,\overline{\ep})$, let $\psi\dvtx \mbbr
\to
[0,1]$ be given by\setcounter{footnote}{1}\footnote{The author owes Professor A.~M. Kulik for this clear-cut
choice of $\psi$.}
\[
\psi(z):=\cases{ %
\exp \bigl\{-(z-\underline{
\ep})^{-1}-(\overline{\ep }-z)^{-1} \bigr\}, & \quad $(\underline{
\ep}<z<\overline{\ep}),$
\vspace*{2pt}\cr
0, & \quad$(\mbox{otherwise})$ }
\]
and set
\[
\nu_{1}(dz):=\bigl\{\psi(z)+\psi(-z)\bigr\}\nu_{\star}(dz)=
\bigl\{\psi(z)+\psi (-z)\bigr\} g_{\star}(z)\,dz.
\]
Then we have the decomposition $\nu=\nu_{1}+\nu_{2}$,
where $\nu_{2}(dz):=[1-\{\psi(z)+\psi(-z)\}]\nu_{\star}(dz)+\nu
_{\natural}(dz)$ defines a L\'evy measure.
The function $z\mapsto\{\psi(z)+\psi(-z)\}g_{\star}(z)$ is smooth and
supported by
$[-\overline{\ep},-\underline{\ep}]\cup[\underline{\ep
},\overline{\ep}]$.
With this truncation in hand, we can apply Kulik~\cite{Kul09}, Proposition A.1,
which states that, when the diffusion part is absent,
the LD condition is implied by the conditions (\ref{eqS}) plus~(\ref{eqN}),
{\renewcommand{\theequation}{$\hat{\mathrm{N}}$}
%e5.4 #&#
\begin{equation}\label{eqN}
\exists x''\in\mbbr\ \exists t''>0
\dvtx P_{x''} [\hat {S}_{t''}=\mbbr ]>0, %\tag{}%\nn
\end{equation}}
\hspace*{-2pt}where
$\hat{S}_{t}:= \{u\mce^{t}_{\tau}c(X_{\tau-});u\in\mbbr,\tau\in
\mcd_{1}\cap(0,t) \}$,
with $\mcd_{1}$ and $(\mce^{t}_{s})_{0\le s\le t}$, respectively, denoting
the domain of the point process $N_{1}$ associated with $\nu_{1}$ and a
right-continuous solution to
\[
\mce^{t}_{s}=1+\int_{s}^{t}
\p_{x}a(X_{u})\mce^{u}_{s}\,du +\int
_{s}^{t}\p_{x}c(X_{u-})
\mce^{u-}_{s}\,dJ_{u}.
\]
As (\ref{eqS}) has been already verified in the previous paragraph, it
remains to prove (\ref{eqN});
obviously, if $\nu$ fulfils Assumption~\ref{A_E2}(ii), then it does
Assumption~\ref{A_E2}(i) too.
The stochastic-exponential formula leads to
\[
\mce^{t}_{s}=\exp(Y_{t}-Y_{s})
\prod_{s<u\le t}(1+\D Y_{u})\exp(-\D
Y_{u}),\qquad  s\le t,
\]
where $Y_{u}:=\int_{0}^{u}\p_{x}a(X_{v})\,dv+\int_{0}^{u}\p
_{x}c(X_{v-})\,dJ_{v}$.
We now introduce the two auxiliary sets
\begin{eqnarray*}
A'(t)&:=&\bigl\{\omega\in\Omega\dvtx \mcd_{1}\cap(0,t)\ne
\varnothing\bigr\},
\\
A''(t)&:=& \bigl\{\omega\in\Omega\dvtx N \bigl((0,t],
\bigl\{z\in\mbbr; |z|\ge \|\p _{x}c\|_{\infty}^{-1}
\bigr\} \bigr)=0 \bigr\},
\end{eqnarray*}
where $N(dt,dz)$ denotes the Poisson random measure associated with
$J$. According to the implications
\begin{eqnarray*}
\bigl\{|\D J_{u}|<\|\p_{x}c\|_{\infty}^{-1}, u
\in(0,t] \bigr\} &\subset& \bigl\{\bigl|\p_{x}c(X_{u-})\D
J_{u}\bigr|<1, u\in(0,t] \bigr\}
\\
&=& \bigl\{|\D Y_{u}|<1, u\in(0,t] \bigr\}
\\
&\subset& \bigl\{\mce^{t}_{s}\ne0,\ s\in[0,t] \bigr\},
\end{eqnarray*}
the process $(\mce^{t}_{s})_{0\le s\le t}$ stays positive a.s. on $A''(t)$.
Since $P[A'(t)\cap A''(t)]>0$ for every $t>0$ and $c$ is nonvanishing
on $\mbbr$,
we observe that for every $x\in\mbbr$ and $t>0$
\begin{eqnarray*}
P_{x}[\hat{S}_{t}=\mbbr]&\ge& P_{x} \bigl[\{
\hat{S}_{t}=\mbbr\}\cap A'(t)\cap A''(t)
\bigr]
\\
&\ge& P_{x} \bigl[\bigl\{\mce^{t}_{s}c(X_{s-})
\ne0\ \mbox{for some $s\in (0,t)$}\bigr\}\cap A'(t)\cap
A''(t) \bigr]
\\
&=&P_{x} \bigl[\bigl\{c(X_{s-})\ne0\ \mbox{for some $s
\in(0,t)$}\bigr\}\cap A'(t)\cap A''(t)
\bigr]
\\
&=&P_{x} \bigl[A'(t)\cap A''(t)
\bigr]>0,
\end{eqnarray*}
hence the LD condition.

\textit{Verification of (II): the drift condition.}
Now we turn to the verification of (\ref{edc}). For verification under
Assumption~\ref{A_E3}(i),
one can refer to Kulik~\cite{Kul09} and Masuda~\cite{Mas07,Mas08};
in this case, we may set $\vp(x)=|x|^{q}$ outside a sufficiently large
neighborhood of the origin.
We are left to showing (\ref{edc}) under Assumption~\ref{A_E3}(ii),
where, compared with Assumption~\ref{A_E3}(i), we impose a weaker
condition on the drift function $a$
while a stronger moment condition on $\nu$.
We will achieve the proof in a somewhat similar manner to the proof of
Masuda~\cite{Mas08}, Theorem 1.2.

Fix any $\ep\in(0,q\|c\|_{\infty}^{-1}\wedge1)$ and pick a $\vp
=\vp
_{\ep}\in\mcc^{2}(\mbbr)$ fulfilling:
\begin{itemize}
\item$\vp(x)=\exp(\ep|x|)$ for $|x|\ge\ep^{-1}$;
\item$\vp(x)\le\exp(\ep|x|)$ for every $x$;
\item$|\p_{x}^{2}\vp(x)|\le C\ep^{2}\vp(x)$ for every $x$.
\end{itemize}
We can write $\mca\vp=\mcg\vp+\mcj\vp$, where
\begin{eqnarray*}
\mcg\vp(x)&:=&\p_{x}\vp(x)a(x)+\frac{1}{2}
\p_{x}^{2}\vp(x)b^{2}(x),
\\
\mcj\vp(x)&:=&\int \bigl\{\vp\bigl(x+c(x)z\bigr)-\vp(x)-\p_{x}\vp
(x)c(x)z \bigr\}\nu(dz).
\end{eqnarray*}
According to the local boundedness of $x\mapsto\mca\vp(x)$,
we may and do concentrate on $x$ with $|x|$ large enough. Direct
algebra gives
\setcounter{equation}{2}
%e5.5 #&#
\begin{equation}
\mcg\vp(x)\le\ep\vp(x) \bigl\{\operatorname{ sgn}(x)a(x)+C\ep \bigr\}.
\label{prop_ergo_1}
\end{equation}
Further, by means of Taylor's theorem and the property of $\vp$,
%
%
%e5.6 #&#
\begin{eqnarray}\label{prop_ergo_2}
\bigl|\mcj\vp(x)\bigr|&\lesssim&\bigl|c(x)\bigr|^{2} \int|z|^{2} \Bigl(\sup
_{0\le s\le1}\bigl\llvert \p_{x}^{2}\vp
\bigl(x+sc(x)z\bigr)\bigr\rrvert \Bigr)\nu(dz)
\nonumber
\\
&\lesssim&\ep^{2}\exp\bigl(\ep|x|\bigr)\int|z|^{2}\exp \bigl(\ep\|c\|
_{\infty
}|z| \bigr)\nu(dz)
\\
&\lesssim&\ep^{2}\vp(x).\nonumber
\end{eqnarray}
By putting (\ref{prop_ergo_1}) and (\ref{prop_ergo_2}) together and by
taking $\ep$ small enough,
we can find a constant $c_{0}>0$ for which $\mca\vp(x)\le-c_{0}\vp(x)$
for every $|x|$ large enough.
The proof of Proposition~\ref{prop_ergo} is complete.
\end{pf}

%re5.5 #&#
\begin{rem}\label{rem_ergo_+1}
If the condition on $\nu$ in Assumption~\ref{A_E2}(i) fails to hold,
then $J$ is necessarily a compound-Poisson process.
In this case, we can utilize the criteria given in Masuda~\cite{Mas08}.
\end{rem}

%re5.6 #&#
\begin{rem}\label{rem_ergo_2}
By combining the results of the LD-condition argument and general
stability theory for Markov processes,
it is possible to formulate subexponential- and polynomial-ergodicity versions,
as well as the ergodicity version (without rate specification);
see, for example, Meyn and Tweedie~\cite{mt3} and Fort and Roberts
\cite
{ForRob05}.
Especially, as in Masuda~\cite{Mas08}, the conditions on $(a,b,c)$ in
Proposition~\ref{prop_ergo}
can be considerably relaxed in case of the ergodicity version,
because the Lyapunov condition required then becomes much weaker.
\end{rem}

%%%%%
%%%%%

%

%

\section*{Acknowledgments}
The author is grateful to two anonymous referees for careful reading
and for several valuable and constructive comments,
which led to substantial improvement of the paper.
He also thanks Professor A.~M. Kulik for his helpful comments
concerning the verification of the LD condition appearing in the proof
of Proposition~\ref{prop_ergo}.

%

%% AOS,AOAS: If there are supplements please fill:
%%\begin{supplement}[id=suppA]
%% \sname{Supplement A}
%% \stitle{Title}
%% \slink[doi]{10.1214/00-AOASXXXXSUPP}
%% \sdatatype{.pdf}"
%% \sdescription{Some text}
%%\end{supplement}

%
%
%accubitu suo, nardus mea dedit odorem suavitatis. Quoniam confortavit
%seras portarum tuarum, benedixit filiis tuis in te. Qui posuit fines
%tuos}

%%%%%
%%%%%

% imsref loaded by akundreckaite, 2013-06-20 07:42:52
%


%suskaldyti doi

\printaddresses

\end{document}